\pdfoutput=1
\documentclass[11pt, oneside]{article}

\usepackage{url}
\usepackage[title]{appendix}
\usepackage{amsmath,amstext,amssymb,amsfonts,amsthm}
\usepackage{stmaryrd}
\usepackage{xcolor}
\usepackage{graphicx}
\usepackage{tikz}
\usetikzlibrary{calc,angles,quotes,arrows.meta,decorations.pathmorphing,fadings,positioning}
\usepackage[american]{circuitikz}
\ctikzset{bipoles/length=0.7cm, resistors/zigs=4, resistors/thickness=0.9}
\usepackage{booktabs}
\usepackage{multirow}
\usepackage[labelfont=bf]{caption}
\usepackage{subcaption}
\usepackage{comment}
\usepackage{algorithm}
\usepackage[noend]{algpseudocode}

\newtheorem{rem}{Remark}[section]

\numberwithin{equation}{section}
\numberwithin{figure}{section}
\numberwithin{table}{section}
\allowdisplaybreaks[4]

\newcommand{\Gm}{\Gamma}

\title{A non-conforming finite difference discrete fracture model\\ based on an energy principle}
\author{Ziyao Xu\footnotemark[1]}
\date{}

\begin{document}

\maketitle
\renewcommand{\thefootnote}{\fnsymbol{footnote}}
\footnotetext[1]{Department of Mathematics and Statistics, Binghamton University, Binghamton, NY 13902, USA. E-mail: zxu24@binghamton.edu}

\begin{center}
\small
\begin{minipage}{0.9\textwidth}
\textbf{Abstract.}
We propose a finite difference discrete fracture model for single-phase flow in fractured porous media. 
The method is derived from a unified energy principle and is implemented on Cartesian grids that need not conform to the fractures. It introduces no additional degrees of freedom, modifies the underlying finite difference scheme only locally, handles both highly conductive fractures and low-permeability barriers, and naturally preserves symmetry and positive definiteness of the discrete system.
Barrier-induced pressure jumps and fracture-induced fluxes can be recovered through inexpensive local post-processing, yielding a sharper representation of the interface effects. 
Numerical experiments based on manufactured solutions and published $2$D and $3$D benchmark problems demonstrate the accuracy, effectiveness, and flexibility of the method for isolated fractures, complex networks, and anisotropic porous media.

\medskip
\textbf{Keywords.} finite difference method, discrete fracture model, conductive fractures, low-permeability barriers, non-conforming grid, energy principle

\end{minipage}
\end{center}
\setlength{\parindent}{2em}

\pagenumbering{arabic}

\section{Introduction}\label{sec:intro}

Fractures in geologic formations are common as a result of the deformation of crustal rocks or human activities.
Depending on the materials they contain and the degree of infilling, they may act as preferential flow pathways or as obstructions to subsurface flow. In this paper, the former are referred to as conductive fractures and the latter as blocking barriers.
Their accurate representation is important in applications involving groundwater flow, oil recovery, geothermal systems, and subsurface contaminant transport.

Because the aperture of a fracture is typically several orders of magnitude smaller than its tangential extent, resolving its full thickness may require prohibitively fine meshes. 
Discrete fracture models (DFMs) bypass this difficulty through a mixed-dimensional approach, treating the surrounding porous matrix in full dimension and the fractures as lower-dimensional objects.
Early works in this area focused primarily on geometry-conforming methods, in which the computational mesh is fitted to fractures.
Noorishad and Mehran \cite{noorishad1982upstream} represented fractures by $1$D line elements and the porous matrix by $2$D quadrilateral elements in a Petrov–Galerkin framework.
Similarly, Baca et al. \cite{baca1984modelling} used a mixed-dimensional finite element formulation for fluid flow and solute transport, coupling $1$D fracture and $2$D matrix elements by superposing their Galerkin contributions at shared nodes.
This superposition-based mixed-dimensional finite element (FE) framework was subsequently extended to multiphase flow \cite{kim2000finite,karimi2003numerical} and incorporated into multiscale finite element formulations \cite{zhang2013accurate}.
Finite volume (FV) DFMs form another major class among conforming discretizations. 
Representative vertex-centered formulations, often referred to as box methods, were developed for conductive fractures \cite{monteagudo2004control,reichenberger2006mixed} and later extended to include low-permeability barriers \cite{glaser2022comparison,xu2025extension} and to achieve higher-order accuracy \cite{liu2026high}. 
Cell-centered FV-DFMs commonly employ either two-point flux approximation (TPFA) \cite{karimi2004efficient,angot2009asymptotic} or multi-point flux approximation (MPFA) \cite{sandve2012efficient,ahmed2015control,glaser2017discrete}.
Other conforming DFMs have also been extensively studied. Examples include mixed finite element methods (MFEM) \cite{martin2005modeling,hoteit2005multicomponent,hoteit2008efficient,zidane2014efficient,frih2012modeling,boon2018robust}, discontinuous Galerkin (DG) methods \cite{antonietti2019discontinuous,chen2023discontinuous,liu2026interior}, and finite difference (FD) methods \cite{liu2018block,liu2020finite}, among others. 
Comparative overviews of representative DFM discretizations are provided in \cite{flemisch2018benchmarks,berre2021verification}.

Conforming DFMs can represent fractures directly and achieve satisfactory approximations, but they generally require fractures to coincide with element edges or faces. 
For complex networks, this requirement can make grid generation a substantial part of the simulation workflow, and may introduce poor-quality cells.
Non-conforming methods avoid this geometric constraint by allowing fractures and barriers to intersect an independently generated background grid.
A particularly influential class of non-conforming methods is the embedded discrete fracture model (EDFM). In the classical EDFM \cite{li2008efficient,moinfar2014development}, fractures are cut by the background grid used for the porous matrix, assigned independent degrees of freedom, and coupled to the matrix cells through transmissibility relations. 
This construction allows arbitrarily oriented conductive fractures to be embedded in a structured background mesh.
The projection-based EDFM (pEDFM) \cite{ctene2017projection} extends this framework to a broader range of conductivity contrasts, including low-permeability barriers, by modifying the matrix–matrix transmissibilities affected by the fracture projection and introducing additional matrix–fracture connections.
The continuous pEDFM (CPEDFM) \cite{rashid2024continuous} further constructs a connected stair-step projection of each fracture plane on the background cell faces, thereby preventing leakage of low-permeability barriers in three dimensions.
Other recent developments improve the near-fracture representation of matrix-fracture coupling and blocking-fracture effects. 
The local EDFM (LEDFM) \cite{losapio2023local} computes near-fracture transmissibilities from fine-scale local conforming problems and can represent both conductive and blocking fractures. The enriched EDFM (nEDFM) \cite{jiao2024enriched} introduces local enriched degrees of freedom per fracture cell, accommodating both conductive and blocking fractures.
These embedded formulations offer an effective balance between geometric flexibility and computational cost, but their construction is primarily based on explicit matrix-fracture connections and transmissibility-based coupling. The global system generally contains independent lower-dimensional fracture-pressure degrees of freedom, while specific extensions modify the transmissibilities or introduce additional enriched variables and equations.

Other routes to non-conforming DFMs include the reinterpreted DFM (RDFM) based on hybrid-dimensional Darcy's law \cite{xu2020hybrid,xu2023hybrid, fu2023hybridizable}, in which conductive fractures and blocking barriers are represented by Dirac delta measures. 
XFEM-based methods \cite{d2012mixed,schwenck2015dimensionally,flemisch2016review,del2017well} enrich the approximation spaces in elements cut by fractures to capture fracture-induced pressure discontinuities and impose the matrix-fracture coupling weakly. 
Lagrange-multiplier methods \cite{koppel2019lagrange,schadle20193d} couple independently meshed matrix and lower-dimensional fracture problems through Lagrange multipliers, whereas immersed finite element methods \cite{zhao2024discrete,zhao2026petrov} incorporate the coupling conditions directly into local basis functions on cut elements. 
These approaches avoid conforming bulk meshes through different mechanisms. 
Many introduce independent fracture variables, enriched local unknowns, or transmissibility-based coupling constructions, while others retain a single bulk system but require discontinuous or hybrid approximation spaces.
Beyond the porous media literature, finite difference methods have also been developed for general elliptic interface problems \cite{leveque1994immersed}. 
In particular, the pressure-jump condition associated with a low-permeability barrier is the classical imperfect-contact condition in elliptic interface problems \cite{dong2026gradient}. 
Second-order finite difference schemes that build this condition into the difference stencils have been developed based on Taylor expansions \cite{cao2022finite,cao2024finite}, and a variational mimetic finite difference method on body-fitted but non-matching meshes has recently been derived from an augmented Dirichlet functional \cite{lipnikov2026variational}. 
These methods address isolated diffusion interfaces rather than the conductive fractures, interface networks, and mixed-dimensional coupling considered here.

These developments motivate a finite difference formulation that handles conductive fractures and low-permeability barriers through local modifications while retaining the original pressure unknowns and algebraic structure. In this work, we develop such a non-conforming discrete fracture model from the energy principle of the continuous interface problem. The matrix flow, barrier transmission, and tangential fracture transport are discretized as contributions to a unified discrete energy, and the numerical scheme is obtained by taking its first variation.
For scalar matrix permeability, a barrier crossing a grid edge is represented by a local pressure-jump variable. Its elimination yields an effective conductance in which the matrix and barrier resistances combine in series.
For full-tensor matrix permeability, we instead extend the local-energy condensation technique developed in \cite{xu2026barrierfem} for a barrier-only unfitted $P^1$ finite element method: cell-local pressure offsets are introduced and eliminated from the energy of each barrier-cut cell, producing a closed-form modification of the local discrete operator.
Conductive fractures contribute local positive-semidefinite updates through their tangential-flow energy, with a common interpolated pressure imposed at fracture junctions. 
The resulting global system introduces no additional interface unknowns, remains symmetric positive definite, and coincides with the underlying five- or nine-point discretization away from the interfaces. 
Barrier-induced pressure jumps and fracture-induced fluxes can also be recovered through inexpensive local post-processing.
The principal contribution of this work is therefore a unified discrete energy construction that reduces fracture and barrier effects to explicit local modifications of the standard finite difference system. 
To the best of our knowledge, no previous finite difference DFM has incorporated both conductive fractures and low-permeability barriers on non-conforming meshes through closed-form local modifications.

The remainder of the paper is organized as follows.
Section~\ref{sec:continuous_model} introduces the interface model and derives the underlying continuous energy principle.
Section~\ref{sec:discretization}, which forms the main part of the paper, develops the standard five- and nine-point finite difference discretizations, the edgewise and cellwise barrier treatments, the fracture contributions, and the local recovery procedure within a unified discrete energy framework.
Section~\ref{sec:numerical_results} presents the $2$D and $3$D numerical experiments.
Section~\ref{sec:conclusion} summarizes the main conclusions and discusses possible directions for future work.
Appendix~\ref{app:three_dimensional_extension} outlines the extension of the method to three dimensions.

\section{Interface model and energy principle}
\label{sec:continuous_model}
This section presents the interface model for conductive fractures and blocking barriers and derives the energy principle underlying the finite difference discretization.

\subsection{Interface model}
Let $\Omega\subset\mathbb R^2$ be a porous domain containing two interfaces, $\Gamma_b$ and $\Gamma_f$, representing low-permeability barriers and highly conductive fractures, respectively.
The porous matrix occupies $\Omega_m:=\Omega\setminus(\Gamma_b\cup\Gamma_f)$.
Its permeability $\mathbf K_m(\mathbf x)$ is assumed to be a symmetric positive definite tensor and may be heterogeneous and anisotropic.

The Darcy velocity and mass-conservation law in the porous matrix are
\begin{equation}
\label{eq:matrix_flow}
\mathbf u=-\mathbf K_m\nabla p,
\qquad
\nabla\cdot\mathbf u=f,
\qquad
\mathbf x\in\Omega_m,
\end{equation}
where $p$ is the pressure and $f$ is the source term. Equivalently,
\begin{equation}
\label{eq:bulk_equation}
-\nabla\cdot(\mathbf K_m\nabla p)=f,
\qquad
\mathbf x\in\Omega_m.
\end{equation}
We prescribe the boundary conditions
\begin{equation}
\label{eq:boundary_conditions}
p=g_D
\quad\text{on }\Gamma_D,
\qquad
\mathbf u\cdot\mathbf n=g_N
\quad\text{on }\Gamma_N,
\end{equation}
where $\mathbf n$ is the outward unit normal to $\partial\Omega$, and
$\partial\Omega=\Gamma_D\cup\Gamma_N$.

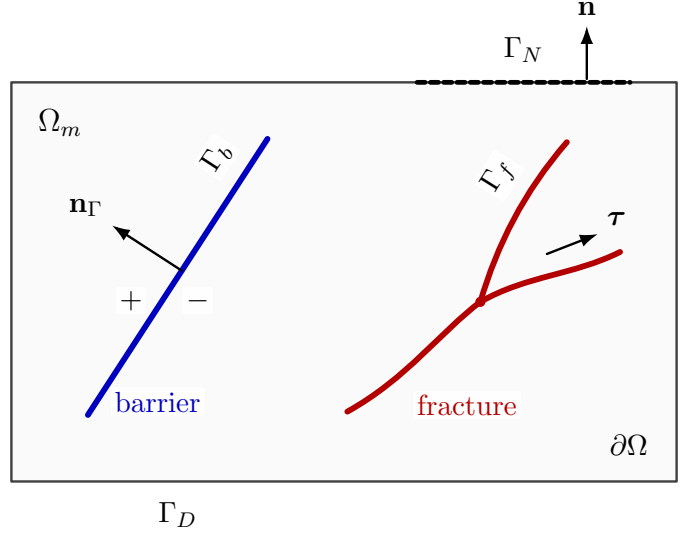
\begin{figure}[t]
\centering
\begin{tikzpicture}[
    x=0.88cm,
    y=0.88cm,
    line cap=round,
    line join=round,
    barrier/.style={
        draw=blue!75!black,
        line width=2.2pt
    },
    fracture/.style={
        draw=red!70!black,
        line width=2.2pt
    },
    vector/.style={
        -{Latex[length=2.4mm,width=1.7mm]},
        line width=0.9pt
    },
    labelbox/.style={
        fill=white,
        fill opacity=0.92,
        text opacity=1,
        inner sep=1.5pt,
        rounded corners=1pt
    }
]

\fill[gray!4] (0,0) rectangle (10,6);
\draw[black!75,line width=0.9pt] (0,0) rectangle (10,6);
\node[anchor=north west,font=\large] at (0.25,5.75) {$\Omega_m$};

\node[below=4pt] at (2.5,0) {$\Gamma_D$};

\draw[black,densely dashed,line width=1.8pt] (6.1,6) -- (9.3,6);
\node[above=4pt] at (7.7,6) {$\Gamma_N$};

\draw[vector] (8.65,6.05) -- (8.65,6.85)
    node[above] {$\mathbf n$};

\draw[barrier]
    (1.15,1.0) -- (3.85,5.15)
    node[pos=0.87,sloped,above=5pt,labelbox] {$\Gamma_b$};

\node[labelbox,text=blue!75!black]
    at (2.20,1.20)
    {barrier};

\coordinate (Bmid) at (2.56,3.17);

\draw[vector]
    (Bmid) -- ++(-1.05,0.68)
    node[above left=-1pt] {$\mathbf n_{\Gamma}$};

\node[labelbox] at (1.80,2.72) {$+$};
\node[labelbox] at (2.80,2.72) {$-$};

\draw[fracture]
    (5.05,1.05)
    .. controls (5.95,1.55) and (6.45,2.25) ..
    (7.05,2.70)
    .. controls (7.75,3.10) and (8.45,3.10) ..
    (9.15,3.45);

\draw[fracture]
    (7.05,2.70)
    .. controls (7.30,3.55) and (7.70,4.35) ..
    (8.35,5.10)
    node[pos=0.66,sloped,above=5pt,labelbox] {$\Gamma_f$};

\fill[red!70!black] (7.05,2.70) circle (1.9pt);

\node[labelbox,text=red!70!black]
    at (6.85,1.15)
    {fracture};

\draw[vector]
    (8.05,3.42) -- (8.82,3.72)
    node[above right=-1pt] {$\boldsymbol\tau$};
    
\node[anchor=south east] at (9.72,0.25) {$\partial\Omega$};

\end{tikzpicture}
\caption{Geometry of the interface model. 
The low-permeability barrier $\Gamma_b$ permits a pressure jump and resists normal flow,
whereas the highly conductive fracture $\Gamma_f$ supports tangential flow. 
The vector $\mathbf n_{\Gamma}$ fixes the orientation of the interface traces, and $\boldsymbol\tau$ denotes an interface tangent.}
\label{fig:interface_geometry}
\end{figure}

We next specify the interface conditions on the barriers and fractures, which represent two important limiting cases of the general interface model developed in \cite{martin2005modeling}.
On each segment of $\Gamma_b\cup\Gamma_f$, choose a unit normal $\mathbf n_{\Gamma}$. The side from which $\mathbf n_{\Gamma}$ points is denoted by $-$, and the side toward which it points by $+$. For any quantity $v$ with traces $v^-$ and $v^+$ on the two sides, we define $\llbracket v\rrbracket:=v^- -v^+$.
The geometry of the interface model and the orientation conventions are illustrated in Figure~\ref{fig:interface_geometry}.

Across $\Gamma_b$, the normal Darcy flux is continuous, while the pressure jump is proportional to the common normal flux:
\begin{equation}
\label{eq:barrier_interface_law}
\mathbf u^-\cdot\mathbf n_{\Gamma}
=
\mathbf u^+\cdot\mathbf n_{\Gamma},
\qquad
\llbracket p\rrbracket
=
\frac{a_b}{k_b}\mathbf u^-\cdot\mathbf n_{\Gamma}
\quad\text{on }\Gamma_b,
\end{equation}
where $a_b$ and $k_b$ are the aperture and permeability of the barrier, respectively.

Across $\Gamma_f$, the pressure is continuous, and the imbalance of the normal Darcy flux from the surrounding matrix is balanced by tangential flow along the fracture:
\begin{equation}
\label{eq:fracture_interface_law}
\llbracket p\rrbracket=0,
\qquad
-\partial_\tau\left(a_fk_f\partial_\tau p\right)
=
\llbracket\mathbf u\cdot\mathbf n_{\Gamma}\rrbracket
\quad\text{on }\Gamma_f,
\end{equation}
where $\partial_\tau$ denotes the tangential derivative, and $a_f$ and $k_f$ are the aperture and permeability of the fracture, respectively. A distributed source along the fracture can be added to the right-hand side of \eqref{eq:fracture_interface_law} without changing the discrete operator and is omitted here.
At a free fracture tip, the tangential fracture flux is set to zero. At a fracture junction, the pressure is shared by all incident branches, and the sum of their outgoing tangential fluxes vanishes. These conditions close the model on the fracture network.

\subsection{Energy principle}
To make the admissible pressure fields explicit, we introduce
\begin{equation}
\label{eq:energy_spaces}
V:=
\left\{
v\in H^1(\Omega_m):
\llbracket v\rrbracket=0 \text{ on }\Gamma_f,
\quad
v|_{\Gamma_f}\in H^1(\Gamma_f)
\right\},
\end{equation}
with $V_D=\{v\in V:\ v=g_D \text{ on }\Gm_D\}$ and the test space $V_0=\{v\in V:\ v=0\text{ on
}\Gm_D\}$.
Here $H^1(\Gamma_f)$ is understood branchwise, with a common value at each fracture junction. 
Since $\Omega_m$ is cut by the interfaces, functions in $H^1(\Omega_m)$ may have distinct traces across $\Gamma_b$, whereas continuity across $\Gamma_f$ is imposed explicitly in \eqref{eq:energy_spaces}.

We now derive the energy formulation underlying the interface model.
Let $p\in V_D$ and $v\in V_0$. 
Testing \eqref{eq:bulk_equation} against $v$ and integrating it by parts over the matrix subregions gives
\begin{align}
\label{eq:bulk_ibp}
\int_{\Omega_m}
\mathbf K_m\nabla p\cdot\nabla v\,d\mathbf x
={}&
\int_{\Omega_m}fv\,d\mathbf x
-
\int_{\Gamma_N}g_Nv\,ds
\notag\\
&-
\int_{\Gamma_b}
\left[
(\mathbf u^-\cdot\mathbf n_{\Gamma})v^-
-
(\mathbf u^+\cdot\mathbf n_{\Gamma})v^+
\right]\,ds
\notag\\
&-
\int_{\Gamma_f}
\left[
(\mathbf u^-\cdot\mathbf n_{\Gamma})v^-
-
(\mathbf u^+\cdot\mathbf n_{\Gamma})v^+
\right]\,ds .
\end{align}
On $\Gamma_b$, the normal-flux continuity in
\eqref{eq:barrier_interface_law} gives
\begin{equation}
\label{eq:barrier_weak_term}
(\mathbf u^-\cdot\mathbf n_{\Gamma})
\llbracket v\rrbracket
=
\frac{k_b}{a_b}
\llbracket p\rrbracket
\llbracket v\rrbracket .
\end{equation}
On $\Gamma_f$, both $p$ and $v$ are continuous, so the corresponding interface contribution in \eqref{eq:bulk_ibp} becomes $-\llbracket\mathbf u\cdot\mathbf n_{\Gamma}\rrbracket v$.
Using \eqref{eq:fracture_interface_law} and integrating by parts along the fracture network yields
\begin{equation}
\label{eq:fracture_weak_term}
-\int_{\Gamma_f}
\llbracket\mathbf u\cdot\mathbf n_{\Gamma}\rrbracket v\,ds
=
-\int_{\Gamma_f}
a_fk_f\,\partial_\tau p\,\partial_\tau v\,ds ,
\end{equation}
where the endpoint and junction terms vanish under the free-tip and junction conditions.

Substituting \eqref{eq:barrier_weak_term} and \eqref{eq:fracture_weak_term} into \eqref{eq:bulk_ibp} and moving the interface contributions to the left-hand side gives the symmetric weak formulation
\begin{align}
\label{eq:weak_form}
\int_{\Omega_m}
\mathbf K_m\nabla p\cdot\nabla v\,d\mathbf x
&+
\int_{\Gamma_b}
\frac{k_b}{a_b}
\llbracket p\rrbracket
\llbracket v\rrbracket\,ds
+
\int_{\Gamma_f}
a_fk_f\,\partial_\tau p\,\partial_\tau v\,ds
\notag\\
&=
\int_{\Omega_m}fv\,d\mathbf x
-
\int_{\Gamma_N}g_Nv\,ds
=: \ell(v),
\qquad
\forall v\in V_0.
\end{align}
The three terms on the left-hand side represent flow through the porous matrix, transmission across the barriers, and tangential transport along the fractures, respectively.

The weak problem \eqref{eq:weak_form} is the Euler equation of the
quadratic energy functional
\begin{align}
\label{eq:continuous_energy}
\mathcal J(q)
={}&
\frac{1}{2}
\int_{\Omega_m}
\mathbf K_m\nabla q\cdot\nabla q\,d\mathbf x
+
\frac{1}{2}
\int_{\Gamma_b}
\frac{k_b}{a_b}
\llbracket q\rrbracket^2\,ds
\notag\\
&+
\frac{1}{2}
\int_{\Gamma_f}
a_fk_f\left|\partial_\tau q\right|^2\,ds
-
\ell(q).
\end{align}
The functional \eqref{eq:continuous_energy} is a mixed-dimensional extension of the classical Dirichlet principle for interface problems. 
Related functionals can be traced back at least to \cite{lipton1996composites,lipton1997variational}, while the same variational principle \eqref{eq:weak_form} was presented in \cite{angot2009asymptotic}, where a cell-centered control volume TPFA discretization of the DFM was developed on conforming meshes.

Equivalently, the pressure is characterized by
\begin{equation}
\label{eq:continuous_minimization}
p
=
\operatorname*{arg\,min}_{q\in V_D}
\mathcal J(q).
\end{equation}
The functional \eqref{eq:continuous_energy} is the starting point for
the finite difference construction developed in the next section.

\section{Energy-based finite difference discretization}
\label{sec:discretization}

This section constructs the finite difference discrete fracture model by discretizing the three contributions to the continuous energy \eqref{eq:continuous_energy}. 
We first introduce the Cartesian grid and the discrete matrix energy, which recover the classical finite difference scheme in the absence of interfaces. The barrier treatment is then derived for scalar matrix permeability and extended to full-tensor permeability. 
Finally, the fracture energy is discretized directly along the fractures. 
The resulting global system contains only the original grid pressure unknowns.

\subsection{Cartesian grid and discrete matrix energy}
\label{sec:discrete_matrix_energy}

The discrete energy functionals provides a systematic framework for deriving finite difference schemes \cite{kuhn2021energy}. 
In the present work, the energy viewpoint is used first to recover the standard background discretization, and then to derive explicit non-conforming fracture and barrier treatments through local modifications.

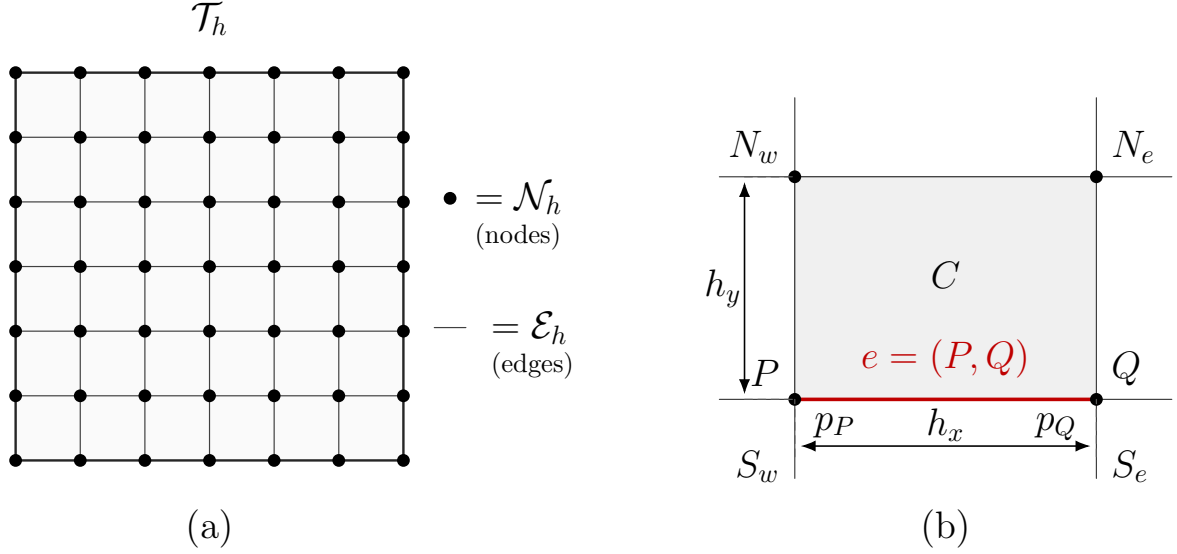
\begin{figure}[t]
\centering
\begin{tikzpicture}[
    x=0.95cm,
    y=0.95cm,
    line cap=round,
    line join=round,
    >=Latex,
    gridline/.style={black!75, line width=0.45pt},
    domainborder/.style={black!85, line width=0.9pt},
    nodept/.style={circle, fill=black, inner sep=1.7pt},
    rededge/.style={draw=red!75!black, line width=1.35pt},
    blueedge/.style={draw=blue!65!black, line width=1.35pt},
    dim/.style={{Latex[length=2mm]}-{Latex[length=2mm]}, line width=0.6pt},
    dashedguide/.style={densely dashed, line width=0.45pt, black!65},
    every node/.style={font=\small}
]

\begin{scope}[shift={(0,0)}]
\path[use as bounding box] (0,0.0) rectangle (8.1,8.0);
\fill[gray!4] (0.9,1.3) rectangle (6.3,6.7);
\draw[domainborder] (0.9,1.3) rectangle (6.3,6.7);

\foreach \x in {1.8,2.7,3.6,4.5,5.4}
    \draw[gridline] (\x,1.3) -- (\x,6.7);
\foreach \y in {2.2,3.1,4.0,4.9,5.8}
    \draw[gridline] (0.9,\y) -- (6.3,\y);

\foreach \x in {0.9,1.8,2.7,3.6,4.5,5.4,6.3}
    \foreach \y in {1.3,2.2,3.1,4.0,4.9,5.8,6.7}
        \node[nodept] at (\x,\y) {};

\node[font=\Large] at (3.6,7.45) {$\mathcal{T}_h$};

\node[nodept] at (6.95,4.95) {};
\node[anchor=west,font=\Large] at (7.15,4.95) {$=\mathcal{N}_h$};
\node[anchor=west] at (7.15,4.42) {(nodes)};

\draw[gridline] (6.73,3.15) -- (7.17,3.15);
\node[anchor=west,font=\Large] at (7.38,3.15) {$=\mathcal{E}_h$};
\node[anchor=west] at (7.38,2.62) {(edges)};

\node[font=\Large] at (3.6,0.35) {(a)};
\end{scope}

\begin{scope}[shift={(10.6,0)}]
\path[use as bounding box] (0,0.0) rectangle (7.4,8.0);

\coordinate (SW) at (1.15,2.15);
\coordinate (SE) at (5.35,2.15);
\coordinate (NW) at (1.15,5.25);
\coordinate (NE) at (5.35,5.25);

\draw[gridline] (0.10,2.15) -- (6.40,2.15);
\draw[gridline] (0.10,5.25) -- (6.40,5.25);
\draw[gridline] (1.15,1.05) -- (1.15,6.35);
\draw[gridline] (5.35,1.05) -- (5.35,6.35);

\fill[gray!12] (1.15,2.15) rectangle (5.35,5.25);
\draw[gridline] (1.15,2.15) rectangle (5.35,5.25);

\draw[rededge] (1.15,2.15) -- (5.35,2.15);

\node[nodept] at (1.15,2.15) {};
\node[nodept] at (5.35,2.15) {};
\node[nodept] at (1.15,5.25) {};
\node[nodept] at (5.35,5.25) {};

\node[anchor=south east,font=\Large] at (1.08,5.30) {$N_w$};
\node[anchor=south west,font=\Large] at (5.42,5.30) {$N_e$};

\node[anchor=north east,font=\Large] at (1.08,1.58) {$S_w$};
\node[anchor=north west,font=\Large] at (5.42,1.58) {$S_e$};

\node[anchor=south east,font=\Large] at (1.08,2.20) {$P$};
\node[anchor=south west,font=\Large] at (5.42,2.20) {$Q$};

\node[anchor=north west,font=\Large] at (1.28,2.1) {$p_P$};
\node[anchor=north east,font=\Large] at (5.22,2.1) {$p_Q$};

\node[font=\Large] at (3.25,3.88) {$C$};

\node[text=red!75!black,font=\Large] at (3.25,2.72) {$e=(P,Q)$};

\draw[dashedguide] (0.45,2.15) -- (1.15,2.15);
\draw[dashedguide] (0.45,5.25) -- (1.15,5.25);
\draw[dim] (0.45,2.22) -- (0.45,5.18);
\node[font=\Large] at (0.15,3.70) {$h_y$};

\draw[dashedguide] (1.15,1.35) -- (1.15,2.15);
\draw[dashedguide] (5.35,1.35) -- (5.35,2.15);
\draw[dim] (1.23,1.50) -- (5.27,1.50);
\node[font=\Large] at (3.25,1.8) {$h_x$};

\node[font=\Large] at (3.25,0.35) {(b)};
\end{scope}

\end{tikzpicture}
\caption{(a) A uniform Cartesian partition $\mathcal{T}_h$ with node set $\mathcal{N}_h$ and edge set $\mathcal{E}_h$. 
(b) A typical rectangular cell $C$ with side lengths $h_x$ and $h_y$. 
The horizontal edge $e=(P,Q)$ is highlighted in red; vertical edges are treated analogously.}
\label{fig:cartesian_grid_and_cell}
\end{figure}

Let $\mathcal T_h$ be a uniform Cartesian partition of $\Omega$ into rectangular cells with side lengths $h_x$ and $h_y$ in the $x$- and $y$-directions, respectively. 
Let $\mathcal N_h$ denote the set of grid nodes and $\mathcal E_h$ the set of horizontal and vertical grid edges.
For each $P\in\mathcal N_h$, let $p_P$ denote the approximation to $p(P)$, and write $\mathbf p=(p_P)_{P\in\mathcal N_h}$ for the vector of nodal pressures. 
The values at the Dirichlet boundary nodes are prescribed by $g_D$.
For a cell $C\in\mathcal T_h$, we label its southwest, southeast, northwest, and northeast vertices by $S_w$, $S_e$, $N_w$, and $N_e$, respectively, and set $\mathbf p_C=
\begin{pmatrix}
p_{S_w} & p_{S_e} & p_{N_w} & p_{N_e}
\end{pmatrix}^{\mathsf T}$.
The Cartesian partition and the local notation for a representative
cell $C$ are illustrated in Figure~\ref{fig:cartesian_grid_and_cell}.

In this subsection, we construct the interface-free discretization of the matrix contribution to the continuous energy
\eqref{eq:continuous_energy},
\begin{equation}
\label{eq:continuous_matrix_energy}
\mathcal E^m(p)
:=
\frac{1}{2}
\int_{\Omega_m}
\mathbf K_m\nabla p\cdot\nabla p\,d\mathbf x
=
\frac{1}{2}
\sum_{C\in\mathcal T_h}
\int_C
\mathbf K_m\nabla p\cdot\nabla p\,d\mathbf x.
\end{equation}
On each cell $C$, define the edge-based derivative approximations
\begin{equation}
\label{eq:cell_edge_derivatives}
\begin{aligned}
D_x^S p_C
&=
\frac{p_{S_e}-p_{S_w}}{h_x},
&
D_x^N p_C
&=
\frac{p_{N_e}-p_{N_w}}{h_x},
\\
D_y^W p_C
&=
\frac{p_{N_w}-p_{S_w}}{h_y},
&
D_y^E p_C
&=
\frac{p_{N_e}-p_{S_e}}{h_y}.
\end{aligned}
\end{equation}

\subsubsection{Edge-separable formulation for scalar permeability}
\label{sec:scalar_matrix_energy}

We first consider scalar matrix permeability, $\mathbf K_m=k_m\mathbf I$, and let $k_C$ be a representative value of $k_m$ on $C$. 
We approximate the cell averages of $(p_x)^2$ and $(p_y)^2$ by $\frac{1}{2}\left(D_x^S p_C\right)^2 +\frac{1}{2}\left(D_x^N p_C\right)^2$ and $\frac{1}{2}\left(D_y^W p_C\right)^2 + \frac{1}{2}\left(D_y^E p_C\right)^2$, respectively. 
This gives the following approximation of the local matrix energy $\frac12\int_C k_m|\nabla p|^2\,d\mathbf x$:
\begin{equation*}
\label{eq:scalar_cell_energy}
\mathcal E_{h,C}^{m}(\mathbf p_C)
=
\frac{k_C}{4}\frac{h_y}{h_x}
\left[
\left(p_{S_e}-p_{S_w}\right)^2
+
\left(p_{N_e}-p_{N_w}\right)^2
\right]
+
\frac{k_C}{4}\frac{h_x}{h_y}
\left[
\left(p_{N_w}-p_{S_w}\right)^2
+
\left(p_{N_e}-p_{S_e}\right)^2
\right].
\end{equation*}

Summing $\mathcal E_{h,C}^{m}(\mathbf p_C)$ over the cells $C\in\mathcal{T}_h$ and collecting the contributions associated with each grid edge yields the global edge-based representation
\begin{equation}
\label{eq:assembled_scalar_matrix_energy}
\mathcal E_h^{m}(\mathbf p)
=
\frac{1}{2}
\sum_{e=(P,Q)\in\mathcal E_h}
c_e\left(p_P-p_Q\right)^2
,
\end{equation}
where the edge conductance is
\begin{equation}
\label{eq:scalar_edge_conductance}
c_e
=
\begin{cases}
\displaystyle
k_e\frac{h_y}{h_x},
& \text{if $e$ is horizontal},\\[2mm]
\displaystyle
k_e\frac{h_x}{h_y},
& \text{if $e$ is vertical}.
\end{cases}
\end{equation}
Here $k_e>0$ denotes the permeability value associated with the edge.
For the cellwise construction above, $k_e$ is the arithmetic average of the permeability values in the two cells adjacent to an interior edge. 
A boundary edge belongs to a single cell and receives half the value of that cell, $k_e=k_C/2$.
More generally, $k_e$ may be replaced by any positive and consistent edge approximation of $k_m$.

The first variation of \eqref{eq:assembled_scalar_matrix_energy} with respect to an interior nodal value is
\begin{equation}
\label{eq:classical_five_point_scheme}
\frac{\partial\mathcal E_h^m}{\partial p_P}
=
\sum_{Q:\,(P,Q)\in\mathcal E_h}
c_{PQ}\left(p_P-p_Q\right).
\end{equation}
This is the classical
\emph{five-point finite difference} operator (up to a sign). 
The edge-separable form will be used in Section~\ref{sec:scalar_barrier_discretization} for scalar permeability: when a barrier cuts an edge, only the corresponding edge conductance $c_e$ is modified.

\subsubsection{Cellwise formulation for full-tensor permeability}
\label{sec:tensor_matrix_energy}

We next consider a symmetric positive definite permeability tensor.
Let
\begin{equation}
\label{eq:cell_permeability_tensor}
\mathbf K_C
=
\begin{pmatrix}
k_{11,C} & k_{12,C}\\
k_{12,C} & k_{22,C}
\end{pmatrix}
\end{equation}
be a representative value of $\mathbf K_m$ on the cell $C$.
We use the same approximations for the cell averages of $(p_x)^2$ and
$(p_y)^2$ as in the scalar case. For the mixed term, we approximate the
cell average of $p_xp_y$ by $\frac{1}{4}
\left(D_x^S p_C+D_x^N p_C\right)
\left(D_y^W p_C+D_y^E p_C\right)$.
The resulting approximation of the local matrix energy
$\frac12\int_C\mathbf K_m\nabla p\cdot\nabla p\,d\mathbf x$ is:
\begin{equation}
\label{eq:tensor_cell_matrix_energy}
\mathcal E_{h,C}^{m}(\mathbf p_C)
=
\frac{1}{2}
\mathbf p_C^{\mathsf T}S_C\mathbf p_C,
\end{equation}
where
\begin{equation}
\label{eq:tensor_cell_matrix}
S_C
=
\frac{1}{2}
\begin{pmatrix}
\alpha_C+\beta_C+\gamma_C
&
-\alpha_C
&
-\gamma_C
&
-\beta_C
\\
-\alpha_C
&
\alpha_C-\beta_C+\gamma_C
&
\beta_C
&
-\gamma_C
\\
-\gamma_C
&
\beta_C
&
\alpha_C-\beta_C+\gamma_C
&
-\alpha_C
\\
-\beta_C
&
-\gamma_C
&
-\alpha_C
&
\alpha_C+\beta_C+\gamma_C
\end{pmatrix},
\end{equation}
with 
\begin{equation*}
\label{eq:scaled_tensor_coefficients}
\alpha_C
=
k_{11,C}\frac{h_y}{h_x},
\qquad
\beta_C
=
k_{12,C},
\qquad
\gamma_C
=
k_{22,C}\frac{h_x}{h_y}.
\end{equation*}

The matrix $S_C$ is symmetric positive semidefinite and satisfies $S_C\mathbf 1=\mathbf 0$, as required for a discrete diffusion energy. 
The global matrix energy is therefore
\begin{equation}
\label{eq:global_tensor_matrix_energy}
\mathcal E_h^{m}(\mathbf p)
=
\frac{1}{2}
\sum_{C\in\mathcal T_h}
\mathbf p_C^{\mathsf T}S_C\mathbf p_C.
\end{equation}

Its first variation at an interior node $P$ is
\begin{equation}
\label{eq:assembled_tensor_scheme}
\frac{\partial\mathcal E_h^m}{\partial p_P}
=
\sum_{\substack{C\in\mathcal T_h\\P\in C}}
\left(S_C\mathbf p_C\right)_P,
\end{equation}
where $\left(S_C\mathbf p_C\right)_P$ denotes the component corresponding to the vertex $P$. 
For a constant full tensor, \eqref{eq:assembled_tensor_scheme} is the standard \emph{nine-point finite difference} operator (up to a sign).

When $k_{12,C}\neq0$, the energy cannot in general be written as a sum of grid-edge contributions.
The complete four-node matrix $S_C$ must then be retained as the basic local object. 
This cellwise form will be used in Section~\ref{sec:tensor_barrier_discretization} for full-tensor permeability: 
when a barrier cuts a cell, only the corresponding local cell matrix $S_C$ is modified.

\subsection{Barrier discretization}
\label{sec:barrier_discretization}

The barrier term in the continuous energy \eqref{eq:continuous_energy},
\[
\frac12
\int_{\Gamma_b}
\frac{k_b}{a_b}\llbracket p\rrbracket^2\,ds,
\]
represents the transmission energy associated with the pressure jump across $\Gamma_b$. 
On a non-conforming Cartesian grid, the nodal pressures do not directly provide separate traces on the two sides of a barrier. 
We therefore introduce local pressure-jump variables on the grid objects (edges for scalar permeability and cells for tensor permeability) intersected by $\Gamma_b$ and construct a discrete energy involving both the nodal pressures and these local variables. 
The jump variables are then eliminated locally by energy minimization, leaving only the original nodal pressures in the global system.

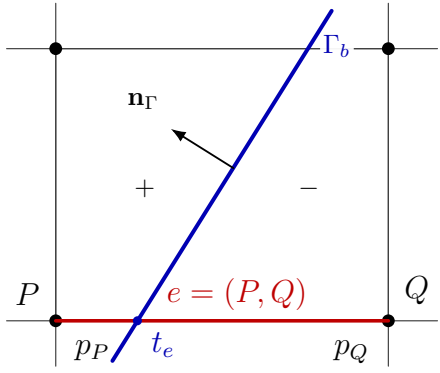
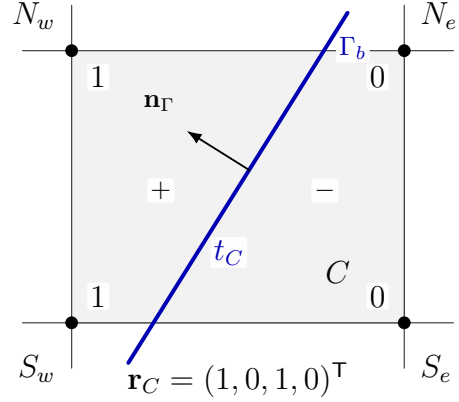
\begin{figure}[t]
\centering

\begin{subfigure}[t]{0.46\textwidth}
\centering
\begin{tikzpicture}[
    x=1.0cm,
    y=1.0cm,
    line cap=round,
    line join=round,
    >=Latex,
    gridline/.style={black!75, line width=0.45pt},
    cellfill/.style={fill=gray!10},
    nodept/.style={circle, fill=black, inner sep=1.7pt},
    edgehl/.style={draw=red!75!black, line width=1.35pt},
    barrier/.style={draw=blue!70!black, line width=1.5pt},
    vector/.style={
        -{Latex[length=2mm,width=1.5mm]},
        line width=0.7pt
    },
    labelbox/.style={
        fill=white,
        fill opacity=0.95,
        text opacity=1,
        inner sep=1.2pt,
        rounded corners=1pt
    },
    every node/.style={font=\small}
]

\path[use as bounding box] (0,0) rectangle (6.8,5.9);

\coordinate (SW) at (1.2,1.25);
\coordinate (SE) at (5.6,1.25);
\coordinate (NW) at (1.2,4.85);
\coordinate (NE) at (5.6,4.85);

\draw[gridline] (0.55,1.25) -- (6.25,1.25);
\draw[gridline] (0.55,4.85) -- (6.25,4.85);
\draw[gridline] (1.2,0.65) -- (1.2,5.45);
\draw[gridline] (5.6,0.65) -- (5.6,5.45);

\draw[gridline] (SW) rectangle (NE);

\draw[barrier] (1.95,0.72) -- (4.85,5.35);
\node[labelbox,text=blue!70!black]
    at (4.92,4.88) {$\Gamma_b$};

\draw[vector]
    (3.55,3.27) -- ++(-0.83,0.52);
\node[anchor=south east]
    at (2.70,3.92) {$\mathbf n_{\Gamma}$};
    
\node[labelbox] at (2.38,3.0) {$+$};
\node[labelbox] at (4.55,3.0) {$-$};

\node[nodept] at (SW) {};
\node[nodept] at (SE) {};
\node[nodept] at (NW) {};
\node[nodept] at (NE) {};

\draw[edgehl] (SW) -- (SE);

\node[anchor=south east,font=\large]
    at (1.12,1.33) {$P$};
\node[anchor=south west,font=\large]
    at (5.68,1.33) {$Q$};

\node[anchor=north west,font=\large]
    at (1.32,1.13) {$p_P$};
\node[anchor=north east,font=\large]
    at (5.48,1.13) {$p_Q$};

\node[font=\large,text=red!75!black]
    at (3.60,1.6) {$e=(P,Q)$};

\coordinate (X) at (2.28,1.25);
\fill[blue!70!black] (X) circle (1.7pt);

\node[labelbox,text=blue!70!black,font=\large]
    at (2.62,0.93) {$t_e$};

\end{tikzpicture}
\caption{Edgewise construction for scalar matrix permeability.}
\label{fig:scalar_cut_edge}
\end{subfigure}
\hfill
\begin{subfigure}[t]{0.46\textwidth}
\centering
\begin{tikzpicture}[
    x=1.0cm,
    y=1.0cm,
    line cap=round,
    line join=round,
    >=Latex,
    gridline/.style={black!75, line width=0.45pt},
    cellfill/.style={fill=gray!10},
    nodept/.style={circle, fill=black, inner sep=1.7pt},
    barrier/.style={draw=blue!70!black, line width=1.5pt},
    vector/.style={
        -{Latex[length=2mm,width=1.5mm]},
        line width=0.7pt
    },
    labelbox/.style={
        fill=white,
        fill opacity=0.95,
        text opacity=1,
        inner sep=1.2pt,
        rounded corners=1pt
    },
    every node/.style={font=\small}
]

\path[use as bounding box] (0,0) rectangle (6.8,5.9);

\coordinate (SW) at (1.2,1.25);
\coordinate (SE) at (5.6,1.25);
\coordinate (NW) at (1.2,4.85);
\coordinate (NE) at (5.6,4.85);

\draw[gridline] (0.55,1.25) -- (6.25,1.25);
\draw[gridline] (0.55,4.85) -- (6.25,4.85);
\draw[gridline] (1.2,0.65) -- (1.2,5.45);
\draw[gridline] (5.6,0.65) -- (5.6,5.45);

\fill[cellfill] (SW) rectangle (NE);
\draw[gridline] (SW) rectangle (NE);

\draw[barrier] (1.95,0.72) -- (4.85,5.35);
\node[labelbox,text=blue!70!black]
    at (4.92,4.88) {$\Gamma_b$};

\draw[vector]
    (3.55,3.27) -- ++(-0.83,0.52);
\node[anchor=south east]
    at (2.70,3.92) {$\mathbf n_{\Gamma}$};
    
\node[labelbox] at (2.38,3.0) {$+$};
\node[labelbox] at (4.55,3.0) {$-$};

\node[nodept] at (SW) {};
\node[nodept] at (SE) {};
\node[nodept] at (NW) {};
\node[nodept] at (NE) {};

\node[labelbox,font=\large] at (1.55,1.60) {$1$};
\node[labelbox,font=\large] at (5.25,1.60) {$0$};
\node[labelbox,font=\large] at (1.55,4.50) {$1$};
\node[labelbox,font=\large] at (5.25,4.50) {$0$};

\node[labelbox,text=blue!70!black,font=\large]
    at (3.28,2.18) {$t_C$};

\node[font=\large] at (4.72,1.92) {$C$};

\node[font=\large]
    at (3.40,0.5)
    {$\mathbf r_C=(1,0,1,0)^{\mathsf T}$};

\node[anchor=north east,font=\large] at (1.12,0.95) {$S_w$};
\node[anchor=north west,font=\large] at (5.68,0.95) {$S_e$};
\node[anchor=south east,font=\large] at (1.12,5.00) {$N_w$};
\node[anchor=south west,font=\large] at (5.68,5.00) {$N_e$};
\end{tikzpicture}
\caption{Cellwise construction for full-tensor matrix permeability.}
\label{fig:tensor_cut_cell}
\end{subfigure}
\caption{Local barrier constructions used in the energy-based finite difference discretization. Left: in the scalar case, a barrier crossing a grid edge introduces an edge-based jump variable. Right: in the full-tensor case, the local construction is carried out on the whole cut cell.}
\label{fig:barrier_local_constructions}
\end{figure}

For scalar matrix permeability, the discrete matrix energy \eqref{eq:assembled_scalar_matrix_energy} separates into edge contributions, so the barrier treatment is performed on the cut edges.
For a full permeability tensor, the coordinate directions are coupled within each cell, so the construction is instead performed at the level of the cut cells. 
The two local constructions are illustrated in Figure~\ref{fig:barrier_local_constructions} and derived below.

\subsubsection{Edgewise construction for scalar permeability}
\label{sec:scalar_barrier_discretization}

Let
\begin{equation}
\label{eq:barrier_cut_edges}
\mathcal E_h^b
=
\left\{
e\in\mathcal E_h:
e^\circ\cap\Gamma_b\neq\varnothing
\right\}
\end{equation}
denote the set of grid edges cut by the barriers. 
We first consider an edge $e=(P,Q)\in\mathcal E_h^b$ whose interior is intersected once by a barrier branch at $\mathbf x_e$.
The case of multiple intersections by barrier networks will be discussed later.

Choose an ordering of $e=(P,Q)$ and regard the edge as oriented from $P$ to $Q$. 
We define the corresponding pressure drop by $\delta_e\mathbf p:=p_P-p_Q$.
The ordering is arbitrary and is introduced only to fix the signs of the local quantities below.

Let $t_e$ be a local variable approximating the signed pressure jump encountered when the barrier is crossed from $P$ to $Q$. 
The total pressure drop along the edge is then decomposed into the pressure variation through the porous matrix and the jump across the barrier: $\delta_e\mathbf p
=
\left(\delta_e\mathbf p-t_e\right)+t_e$.
Accordingly, $\delta_e\mathbf p-t_e$ represents the part of the edge pressure drop associated with the porous matrix.

To relate the sign of $t_e$ to the interface convention introduced in Section~\ref{sec:continuous_model}, let $\mathbf x_P$ and $\mathbf x_Q$ denote the coordinates of the endpoints of $e$ and define $\sigma_e
=
\operatorname{sign}
\left(
(\mathbf x_Q-\mathbf x_P)\cdot
\mathbf n_\Gamma(\mathbf x_e)
\right)$.
Thus, $t_e$ approximates $\sigma_e\llbracket p\rrbracket(\mathbf x_e)$: 
it has the sign of $\llbracket p\rrbracket$ when the edge crosses the barrier from the $-$ side to the $+$ side and the opposite sign when the crossing direction is reversed. 
Reversing the ordering of the edge changes the signs of both $\delta_e\mathbf p$ and $t_e$ and therefore leaves the resulting energy unchanged.

Let $\omega_e>0$ be the quadrature weight associated with the barrier-edge intersection $\mathbf x_e$. 
We approximate the barrier energy by
\begin{equation}
\label{eq:barrier_line_quadrature}
\int_{\Gamma_b}
\frac{k_b}{a_b}\llbracket p\rrbracket^2\,ds
\approx
\sum_{e\in\mathcal E_h^b}
\omega_e
\frac{k_{b,e}}{a_{b,e}}
t_e^2,
\end{equation}
where $k_{b,e}$ and $a_{b,e}$ are representative values of the barrier permeability and aperture at $\mathbf x_e$.

In this work, we use the projected crossing weight
\begin{equation}
\label{eq:projected_barrier_weight}
\omega_e
=
\begin{cases}
|\mathbf{n}_{\Gamma,x}(\mathbf x_e)|h_y,
& \text{if $e$ is horizontal},\\
|\mathbf{n}_{\Gamma,y}(\mathbf x_e)|h_x,
& \text{if $e$ is vertical}.
\end{cases}
\end{equation}
Its form follows from
the decomposition $ds
=
|\mathbf{n}_{\Gamma,x}|\,|dy|
+
|\mathbf{n}_{\Gamma,y}|\,|dx|$.
On an edge lying on the domain boundary, whose control face is half a cell wide, the weight $\omega_e$ is halved, in the same way as the boundary value of the edge conductance $c_e$.

Define the discrete barrier conductance associated with the cut edge by
\begin{equation}
\label{eq:discrete_edge_barrier_conductance}
d_e
:=
\omega_e\frac{k_{b,e}}{a_{b,e}}.
\end{equation}
For an uncut edge, the matrix energy contribution is still $\frac12c_e(\delta_e\mathbf p)^2$ following \eqref{eq:assembled_scalar_matrix_energy}, where $c_e$ is the edge conductance defined in \eqref{eq:scalar_edge_conductance}. 
If the edge is cut by a barrier, we replace this term by the augmented local energy
\begin{equation}
\label{eq:augmented_scalar_barrier_energy}
\mathcal E_{h,e}^{m,b}(\mathbf p,t_e)
=
\frac{1}{2}
c_e\left(\delta_e\mathbf p-t_e\right)^2
+
\frac{1}{2}d_et_e^2.
\end{equation}
The first term is the matrix energy associated with the jump-corrected pressure difference, whereas the second term is the discrete barrier energy.

The discrete solution is obtained by minimizing the sum of these edge energies, together with the load contribution, over the nodal pressures $\mathbf{p}$ and the local jump variables $\{t_e\}_{e\in\mathcal E_h^b}$. 
Since each $t_e$ appears only in the augmented energy of its own cut edge, the minimization with respect to $t_e$ decouples edge by edge. We may therefore first eliminate $t_e$ by minimizing \eqref{eq:augmented_scalar_barrier_energy} for fixed $\mathbf p$, which gives 
\begin{equation*}
\label{eq:optimal_edge_jump}
t_e^\star
=
\frac{c_e}{c_e+d_e}\,
\delta_e\mathbf p.
\end{equation*}
Substituting it into \eqref{eq:augmented_scalar_barrier_energy} yields the effective energy for the cut edge,
\begin{equation}
\label{eq:condensed_scalar_barrier_energy}
\widetilde{\mathcal E}_{h,e}^{m,b}(\mathbf p)
=
\frac{1}{2}
\widetilde c_e
\left(\delta_e\mathbf p\right)^2,
\end{equation}
where
\begin{equation}
\label{eq:effective_barrier_edge_conductance}
\widetilde c_e
=
\frac{c_ed_e}{c_e+d_e}
=
\left(
\frac{1}{c_e}
+
\frac{1}{d_e}
\right)^{-1}.
\end{equation}
Thus, the matrix resistance $1/c_e$ and the discrete barrier resistance $1/d_e$ combine in series.

Define the effective edge conductance by
\begin{equation}
\label{eq:scalar_effective_edge_conductance}
\widehat c_e
=
\begin{cases}
c_e,
& e\notin\mathcal E_h^b,\\
\widetilde c_e,
& e\in\mathcal E_h^b.
\end{cases}
\end{equation}
The effective matrix-barrier energy is therefore
\begin{equation}
\label{eq:global_scalar_barrier_energy}
\widetilde{\mathcal E}_h^{m,b}(\mathbf p)
=
\frac12
\sum_{e=(P,Q)\in\mathcal E_h}
\widehat c_e
\left(p_P-p_Q\right)^2.
\end{equation}

Its first variation retains the five-point form:
\begin{equation}
\label{eq:scalar_barrier_nodal_balance}
\frac{\partial\widetilde{\mathcal E}_h^{m,b}}
{\partial p_P}
=
\sum_{Q:\,(P,Q)\in\mathcal E_h}
\widehat c_{PQ}\left(p_P-p_Q\right).
\end{equation} 

Because only the conductances of the cut edges are changed, the method introduces no additional global unknowns and preserves the original five-point sparsity pattern. 
Moreover, since $\widehat c_e>0$, its off-diagonal entries are nonpositive and its row sums vanish before Dirichlet conditions are imposed. 
If the edge graph contains a Dirichlet node, the reduced matrix is a symmetric nonsingular $M$-matrix. 
The corresponding discrete comparison and maximum-principle properties are therefore retained.

\begin{rem}[Arc-length weights]
\label{rem:arc_length_lumping}
An alternative quadrature for \eqref{eq:barrier_line_quadrature} is obtained by  assigning an arc-length portion $\omega_e^{\rm arc}>0$ to each intersection $\mathbf{x}_{e}$. 
The portions can be chosen using the arc-length midpoints between consecutive intersections.
Replacing $\omega_e$ by $\omega_e^{\rm arc}$ in \eqref{eq:discrete_edge_barrier_conductance} leaves all subsequent elimination and condensation formulas unchanged.
It performs similarly to the projected weights in our numerical tests.
We use the projected weights \eqref{eq:projected_barrier_weight} because they are computed independently at each barrier-edge intersection and do not require non-local information.
\end{rem}

\subsubsection{Cellwise construction for full-tensor permeability}
\label{sec:tensor_barrier_discretization}

For a full-tensor matrix permeability, the mixed-derivative contribution couples all four vertices of a Cartesian cell. 
Consequently, the cell energy
\[
\mathcal E_{h,C}^{m}(\mathbf p_C)
=
\frac12\mathbf p_C^{\mathsf T}S_C\mathbf p_C
\]
cannot, in general, be decomposed into independent edge energies. 
A barrier must therefore be incorporated directly into the complete energy of each cut cell rather than through separate modifications of the cut edges.

We first consider a cell $C$ crossed by a single barrier branch such that $\Gamma_{b,C}:=\Gamma_b\cap C$ is a nondegenerate chord, with length $\ell_C:=|\Gamma_{b,C}|$, separating the vertices of $C$ into two nonempty groups. 
Let $\mathbf r_C\in\{0,1\}^4$ be the corresponding side-indicator vector, with entries equal to zero on one side of the barrier and one on the other; see Figure \ref{fig:tensor_cut_cell}. 
The choice of sides is arbitrary.

Introduce a cell-local jump variable $t_C$ representing the signed pressure offset between the two vertex groups. The vector $\mathbf p_C-\mathbf r_Ct_C$ then represents the nodal pressure after this offset has been removed.
For the barrier contribution, define
\begin{equation}
\label{eq:tensor_cell_barrier_conductance}
d_C
:=
\ell_C\frac{k_{b,C}}{a_{b,C}}\approx
\int_{\Gamma_{b,C}}
\frac{k_b}{a_b}\,ds,
\end{equation}
where $k_{b,C}$ and $a_{b,C}$ are representative barrier parameters on $\Gamma_{b,C}$. 

The augmented energy of the barrier-cut cell is defined by
\begin{equation}
\label{eq:augmented_tensor_barrier_energy}
\mathcal E_{h,C}^{m,b}(\mathbf p_C,t_C)
=
\frac12
(\mathbf p_C-\mathbf r_Ct_C)^{\mathsf T}
S_C
(\mathbf p_C-\mathbf r_Ct_C)
+
\frac12d_Ct_C^2.
\end{equation}
The first term is the matrix energy of the jump-corrected nodal pressure, whereas the second is the discrete transmission energy of the barrier.

The global discrete energy is minimized over the nodal pressures and all cell-local jump variables. 
Since $t_C$ appears only in \eqref{eq:augmented_tensor_barrier_energy}, its elimination decouples cell by cell. 
For fixed $\mathbf p_C$, the stationarity condition is
\begin{equation}
\label{eq:tensor_cell_jump_stationarity}
\frac{\partial \mathcal{E}^{m,b}_{h,C}}{\partial t_{C}}=-\mathbf r_C^{\mathsf T}S_C
(\mathbf p_C-\mathbf r_Ct_C)
+
d_Ct_C
=
0,
\end{equation}
which gives
\begin{equation}
\label{eq:optimal_tensor_cell_jump}
t_C^\star
=
\frac{\mathbf r_C^{\mathsf T}S_C\mathbf p_C}
{\mathbf r_C^{\mathsf T}S_C\mathbf r_C+d_C}.
\end{equation}
Substituting \eqref{eq:optimal_tensor_cell_jump} into \eqref{eq:augmented_tensor_barrier_energy} yields an effective energy on the cut cell,
\begin{equation}
\label{eq:condensed_tensor_barrier_energy}
\widetilde{\mathcal E}_{h,C}^{m,b}(\mathbf p_C)
=
\frac12
\mathbf p_C^{\mathsf T}\widetilde S_C\mathbf p_C,
\end{equation}
where
\begin{equation}
\label{eq:modified_tensor_cell_matrix}
\widetilde S_C
=
S_C
-
\frac{
(S_C\mathbf r_C)(S_C\mathbf r_C)^{\mathsf T}
}{
\mathbf r_C^{\mathsf T}S_C\mathbf r_C+d_C
}.
\end{equation}
Thus, the effect of the barrier is a symmetric rank-one reduction of the original cell matrix.

This construction is independent of which vertex group is assigned the value one. 
Indeed, replacing $\mathbf r_C$ by $\mathbf 1-\mathbf r_C$ leaves
\eqref{eq:modified_tensor_cell_matrix} unchanged, since $S_C\mathbf 1=\mathbf0$.
Moreover, the same formulation applies to both the $1$--$3$ and $2$--$2$ vertex partitions induced by the barrier chord.

The modified local cell matrix is symmetric positive semidefinite. 
This follows directly from
\begin{equation*}
\label{eq:tensor_condensed_energy_minimum}
\frac12
\mathbf p_C^{\mathsf T}\widetilde S_C\mathbf p_C
=
\min_{t_C\in\mathbb R}
\mathcal E_{h,C}^{m,b}(\mathbf p_C,t_C)
\geq 0.
\end{equation*}
Moreover, $\widetilde S_C\mathbf1=\mathbf0$, so the constant-pressure null mode of the original cell energy is preserved. 

Let $\mathcal T_h^b$ denote the set of cells cut by barriers and define the effective local four node matrix
\begin{equation}
\label{eq:tensor_effective_cell_matrix}
\widehat S_C
=
\begin{cases}
S_C,
& C\notin\mathcal T_h^b,\\
\widetilde S_C,
& C\in\mathcal T_h^b.
\end{cases}
\end{equation}
The effective matrix-barrier energy is therefore
\begin{equation}
\label{eq:global_tensor_barrier_energy}
\widetilde{\mathcal E}_h^{m,b}(\mathbf p)
=
\frac12
\sum_{C\in\mathcal T_h}
\mathbf p_C^{\mathsf T}\widehat S_C\mathbf p_C.
\end{equation}
Its first variation retains the nine-point form
\begin{equation}
\label{eq:tensor_barrier_nodal_balance}
\frac{\partial\widetilde{\mathcal E}_h^{m,b}}
{\partial p_P}
=
\sum_{\substack{C\in\mathcal T_h\\P\in C}}
\left(\widehat S_C\mathbf p_C\right)_P.
\end{equation}
Hence only the local matrices of barrier-cut cells are changed.
No additional global unknowns are introduced, and the original global sparsity pattern is preserved. 

Unlike the scalar edgewise method, the cellwise construction does not generally yield a graph-Laplacian or $M$-matrix structure; such sign properties are not generally available even for the classical nine-point finite difference operator. 
The construction does preserve the symmetry and positive semidefiniteness of every local contribution and, under suitable Dirichlet constraints, the positive definiteness of the assembled system.

Even when the permeability tensor becomes scalar, the present cellwise construction is not generally identical to the edgewise method in Section~\ref{sec:scalar_barrier_discretization}. 
We use the edgewise method whenever the matrix permeability is a scalar because it retains the stronger five-point graph-Laplacian and $M$-matrix structure, and use the cellwise method when full-tensor coupling makes an edgewise decomposition unavailable.

\subsubsection{Multiple barrier segments and junction cells}
\label{sec:multiple_barrier_segments}

For the scalar edgewise construction, a barrier network requires no separate junction treatment. 
Each barrier-edge intersection contributes locally to the resistance of the corresponding edge. 
More precisely, let $\mathcal A_e$ denote the set of barrier crossings of an edge $e$, and define
\[
d_{e,\alpha}
=
\omega_{e,\alpha}
\frac{k_{b,e,\alpha}}{a_{b,e,\alpha}},
\qquad
\alpha\in\mathcal A_e.
\]
Introducing one signed jump variable $t_{e,\alpha}$ for each crossing gives the augmented edge energy
\begin{equation}
\label{eq:multiple_barrier_edge_energy}
\mathcal E_{h,e}^{m,b}
=
\frac12c_e
\left(
\delta_e\mathbf p
-
\sum_{\alpha\in\mathcal A_e}t_{e,\alpha}
\right)^2
+
\frac12
\sum_{\alpha\in\mathcal A_e}
d_{e,\alpha}t_{e,\alpha}^2.
\end{equation}
Eliminating these local variables yields
\begin{equation}
\label{eq:multiple_barrier_edge_conductance}
\widetilde c_e
=
\left(
\frac1{c_e}
+
\sum_{\alpha\in\mathcal A_e}
\frac1{d_{e,\alpha}}
\right)^{-1}.
\end{equation}
Thus, multiple barrier crossings of the same edge contribute resistances in series, while crossings on different edges are treated independently. 
The single-crossing formula \eqref{eq:effective_barrier_edge_conductance} is recovered when $|\mathcal A_e|=1$.

For the cellwise full-tensor construction, however, multiple barrier segments within the same cell are coupled through the common cell energy and must be treated jointly. We now describe this multiregion extension.

Let $C_0,C_1,\ldots,C_m$ denote the connected components of $C\setminus\Gamma_b$. 
We choose $C_0$ as a reference region and introduce the vector of regionwise pressure offsets $\mathbf t_C
=
(t_{C,1},\ldots,t_{C,m})^{\mathsf T}$,
with the offset in $C_0$ fixed to zero.
Let $E_C\in\{0,1\}^{4\times m}$ be the vertex-region incidence matrix defined by
\begin{equation}
\label{eq:vertex_region_incidence}
(E_C)_{ij}
=
\begin{cases}
1,
& \text{if the $i$th vertex of $C$ belongs to $C_j$},\\
0,
& \text{otherwise},
\end{cases}
\qquad j=1,\ldots,m.
\end{equation}
Rows corresponding to vertices in the reference region $C_0$ are therefore zero. 
The vector $\mathbf p_C-E_C\mathbf t_C$ represents the nodal pressure after the regionwise offsets have been removed.

Suppose that the barrier network inside $C$ consists of local arms $\{\gamma_{C,\alpha}\}_{\alpha\in\mathcal A_C}$ and $\gamma_{C,\alpha}$ separates the regions $C_{i_\alpha}$ and $C_{j_\alpha}$. 
Define
\begin{equation}
\label{eq:barrier_arm_conductance}
d_{C,\alpha}
:=|\gamma_{C,\alpha}|
\frac{k_{b,\alpha}}{a_{b,\alpha}}\approx
\int_{\gamma_{C,\alpha}}
\frac{k_b}{a_b}\,ds
.
\end{equation}
Let $\boldsymbol e_0=\mathbf0\in\mathbb R^m$ and let $\boldsymbol e_j=(0,\cdots,1,\cdots,0)^{\mathsf T}$ be the $j$th coordinate vector for $j=1,\ldots,m$. 
The offset difference across the arm is represented, up to orientation, by $(\boldsymbol e_{i_\alpha}-\boldsymbol e_{j_\alpha})^{\mathsf T} \mathbf t_C$.

The total barrier energy in the cell can then be written as
\[
\frac12
\sum_{\alpha\in\mathcal A_C}
d_{C,\alpha}
\left(
(\boldsymbol e_{i_\alpha}
-
\boldsymbol e_{j_\alpha})^{\mathsf T}\mathbf t_C
\right)^2
=
\frac12\mathbf t_C^{\mathsf T}P_C\mathbf t_C,
\]
where
\begin{equation}
\label{eq:barrier_penalty_matrix}
P_C
=
\sum_{\alpha\in\mathcal A_C}
d_{C,\alpha}
(\boldsymbol e_{i_\alpha}
-
\boldsymbol e_{j_\alpha})(\boldsymbol e_{i_\alpha}
-
\boldsymbol e_{j_\alpha})^{\mathsf T}.
\end{equation}

The augmented cell energy is
\begin{equation}
\label{eq:augmented_multiregion_barrier_energy}
\mathcal E_{h,C}^{m,b}
(\mathbf p_C,\mathbf t_C)
=
\frac12
(\mathbf p_C-E_C\mathbf t_C)^{\mathsf T}
S_C
(\mathbf p_C-E_C\mathbf t_C)
+
\frac12
\mathbf t_C^{\mathsf T}P_C\mathbf t_C.
\end{equation}
The local stationarity equation for the region offsets is
\begin{equation}
\label{eq:multiregion_jump_system}
\left(
E_C^{\mathsf T}S_CE_C+P_C
\right)\mathbf t_C
=
E_C^{\mathsf T}S_C\mathbf p_C.
\end{equation}
The matrix $P_C$ on the left-hand side is positive
definite, hence
\begin{equation}
\label{eq:optimal_multiregion_jump}
\mathbf t_C^\star
=
\left(
E_C^{\mathsf T}S_CE_C+P_C
\right)^{-1}
E_C^{\mathsf T}S_C\mathbf p_C.
\end{equation}
Substitution into \eqref{eq:augmented_multiregion_barrier_energy} gives the effective cell matrix
\begin{equation}
\label{eq:modified_multiregion_cell_matrix}
\widetilde S_C
=
S_C
-
S_CE_C
\left(
E_C^{\mathsf T}S_CE_C+P_C
\right)^{-1}
E_C^{\mathsf T}S_C.
\end{equation}
The single-chord formula
\eqref{eq:modified_tensor_cell_matrix} is recovered when $m=1$,
$E_C=\mathbf r_C$, and $P_C=d_C$.

As in the single-chord case, $\widetilde S_C$ is symmetric positive semidefinite and satisfies $\widetilde S_C\mathbf1=\mathbf0$.
All subregion offsets are eliminated locally, so the global unknowns remain the original nodal pressures and the global sparsity pattern is unchanged. 
An isolated free tip that does not separate the cell into distinct connected regions introduces no additional offset variable and is ignored in the treatment.
We assume that barriers do not pass exactly through grid vertices or overlap cell edges; such degenerate incidences can be resolved by a small geometric perturbation.

The barrier treatment therefore modifies the matrix energy only through local contributions and introduces no additional global degrees of freedom.
We next discretize the tangential energy carried by the conductive fractures.

\subsection{Fracture discretization}
\label{sec:fracture_discretization}

Unlike a barrier, a conductive fracture preserves pressure continuity and contributes an additional tangential-flow energy,
\[
\frac12
\int_{\Gamma_f}
a_fk_f|\partial_\tau p|^2\,ds.
\]
Its discrete contribution can therefore be added directly to the matrix energy, without introducing local variables.

\begin{figure}[t]
\centering

\begin{subfigure}[t]{0.46\textwidth}
\centering
\begin{tikzpicture}[
    x=1.0cm,
    y=1.0cm,
    line cap=round,
    line join=round,
    >=Latex,
    gridline/.style={black!75, line width=0.45pt},
    cellfill/.style={fill=gray!10},
    nodept/.style={circle, fill=black, inner sep=1.7pt},
    fracline/.style={draw=red!75!black, line width=1.5pt},
    fracpt/.style={circle, fill=red!75!black, inner sep=1.5pt},
    vector/.style={
        -{Latex[length=2mm,width=1.5mm]},
        line width=0.7pt
    },
    labelbox/.style={
        fill=white,
        fill opacity=0.95,
        text opacity=1,
        inner sep=1.2pt,
        rounded corners=1pt
    },
    every node/.style={font=\small}
]

\path[use as bounding box] (0,0) rectangle (6.8,5.9);

\coordinate (SW) at (1.2,1.25);
\coordinate (SE) at (5.6,1.25);
\coordinate (NW) at (1.2,4.85);
\coordinate (NE) at (5.6,4.85);

\draw[gridline] (0.55,1.25) -- (6.25,1.25);
\draw[gridline] (0.55,4.85) -- (6.25,4.85);
\draw[gridline] (1.2,0.65) -- (1.2,5.45);
\draw[gridline] (5.6,0.65) -- (5.6,5.45);

\fill[cellfill] (SW) rectangle (NE);
\draw[gridline] (SW) rectangle (NE);

\node[font=\large] at (4.95,1.75) {$C$};

\node[nodept] at (SW) {};
\node[nodept] at (SE) {};
\node[nodept] at (NW) {};
\node[nodept] at (NE) {};

\node[anchor=north east,font=\large] at (1.12,0.95) {$S_w$};
\node[anchor=north west,font=\large] at (5.68,0.95) {$S_e$};
\node[anchor=south east,font=\large] at (1.12,5.00) {$N_w$};
\node[anchor=south west,font=\large] at (5.68,5.00) {$N_e$};

\coordinate (A) at (1.20,2.15);
\coordinate (B) at (4.70,4.85);
\draw[fracline] (A) -- (B);
\node[fracpt] at (A) {};
\node[fracpt] at (B) {};

\node[labelbox,anchor=east] at (1.03,2.15)
    {$\mathbf a_{C,\alpha}$};
\node[labelbox,anchor=south west] at (4.78,4.93)
    {$\mathbf b_{C,\alpha}$};

\node[labelbox] at (2.00,2.0)
    {$p_h(\mathbf a_{C,\alpha})$};
\node[labelbox] at (4.8,4.1)
    {$p_h(\mathbf b_{C,\alpha})$};

\node[text=red!75!black,font=\large] at (3.35,3.2)
    {$\gamma_{C,\alpha}$};

\coordinate (TauStart) at ($(A)!0.36!(B)+(-0.30,0.39)$);
\coordinate (TauEnd)   at ($(A)!0.64!(B)+(-0.30,0.39)$);

\draw[vector]
    (TauStart) -- (TauEnd)
    node[midway, above left, xshift=-1pt, yshift=1pt] {$\tau$};

\end{tikzpicture}
\caption{A fracture chord crossing a cell.}
\label{fig:fracture_single_chord}
\end{subfigure}
\hfill
\begin{subfigure}[t]{0.46\textwidth}
\centering
\begin{tikzpicture}[
    x=1.0cm,
    y=1.0cm,
    line cap=round,
    line join=round,
    >=Latex,
    gridline/.style={black!75, line width=0.45pt},
    cellfill/.style={fill=gray!10},
    nodept/.style={circle, fill=black, inner sep=1.7pt},
    fracline/.style={draw=red!75!black, line width=1.5pt},
    fracpt/.style={circle, fill=red!75!black, inner sep=1.5pt},
    vector/.style={
        -{Latex[length=2mm,width=1.5mm]},
        line width=0.7pt
    },
    labelbox/.style={
        fill=white,
        fill opacity=0.95,
        text opacity=1,
        inner sep=1.2pt,
        rounded corners=1pt
    },
    every node/.style={font=\small}
]

\path[use as bounding box] (0,0) rectangle (6.8,5.9);

\coordinate (SW) at (1.2,1.25);
\coordinate (SE) at (5.6,1.25);
\coordinate (NW) at (1.2,4.85);
\coordinate (NE) at (5.6,4.85);

\draw[gridline] (0.55,1.25) -- (6.25,1.25);
\draw[gridline] (0.55,4.85) -- (6.25,4.85);
\draw[gridline] (1.2,0.65) -- (1.2,5.45);
\draw[gridline] (5.6,0.65) -- (5.6,5.45);

\fill[cellfill] (SW) rectangle (NE);
\draw[gridline] (SW) rectangle (NE);

\node[font=\large] at (4.95,1.75) {$C$};

\node[nodept] at (SW) {};
\node[nodept] at (SE) {};
\node[nodept] at (NW) {};
\node[nodept] at (NE) {};

\node[anchor=north east,font=\large] at (1.12,0.95) {$S_w$};
\node[anchor=north west,font=\large] at (5.68,0.95) {$S_e$};
\node[anchor=south east,font=\large] at (1.12,5.00) {$N_w$};
\node[anchor=south west,font=\large] at (5.68,5.00) {$N_e$};

\coordinate (J) at (3.45,2.95);

\coordinate (L) at (1.20,2.10);
\coordinate (T) at (3.10,4.85);
\coordinate (R) at (5.60,3.85);

\draw[fracline] (L) -- (J);
\draw[fracline] (T) -- (J);
\draw[fracline] (R) -- (J);

\node[fracpt] at (J) {};
\node[fracpt] at (L) {};
\node[fracpt] at (T) {};
\node[fracpt] at (R) {};

\node[labelbox,anchor=west] at (3.62,2.78)
    {$\mathbf x_J$};
\node[labelbox] at (4.30,2.35)
    {$p_h(\mathbf x_J)$};

\node[text=red!75!black,font=\large] at (2.15,2.65)
    {$\gamma_{C,1}$};
\node[text=red!75!black,font=\large] at (2.80,4.10)
    {$\gamma_{C,2}$};
\node[text=red!75!black,font=\large] at (4.95,3.20)
    {$\gamma_{C,3}$};

\coordinate (TauOneStart) at ($(L)!0.18!(J)+(0.10,-0.34)$);
\coordinate (TauOneEnd)   at ($(L)!0.58!(J)+(0.10,-0.34)$);

\draw[vector]
    (TauOneStart) -- (TauOneEnd)
    node[midway, below right, xshift=1pt, yshift=-1pt] {$\tau$};

\coordinate (TauTwoStart) at ($(T)!0.15!(J)+(0.38,0.02)$);
\coordinate (TauTwoEnd)   at ($(T)!0.52!(J)+(0.38,0.02)$);

\draw[vector]
    (TauTwoStart) -- (TauTwoEnd)
    node[midway, right, xshift=1pt] {$\tau$};

\coordinate (TauThreeStart) at ($(R)!0.16!(J)+(0.10,0.34)$);
\coordinate (TauThreeEnd)   at ($(R)!0.56!(J)+(0.10,0.34)$);

\draw[vector]
    (TauThreeStart) -- (TauThreeEnd)
    node[midway, above, yshift=1pt] {$\tau$};

\end{tikzpicture}
\caption{A fracture junction inside a cell.}
\label{fig:fracture_junction_cell}
\end{subfigure}

\caption{Local fracture constructions used in the energy-based finite difference discretization. Left: a fracture segment $\gamma_{C,\alpha}=\overline{\mathbf a_{C,\alpha}\mathbf b_{C,\alpha}}$ on a cell $C$. 
Its contribution is determined by the interpolated pressure difference between the two endpoints. 
Right: several fracture arms meeting at a junction point $\mathbf x_J$ inside a cell. 
All incident arms share the same interpolated junction pressure $p_h(\mathbf x_J)$, and their tangential energies are added locally.}
\label{fig:fracture_local_constructions}
\end{figure}
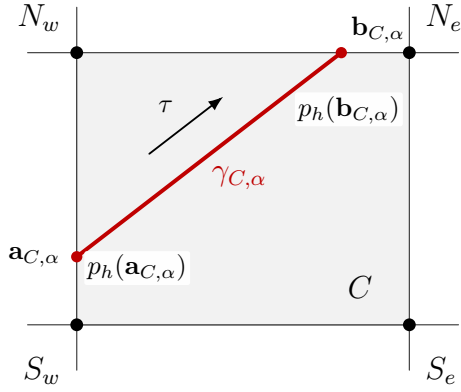
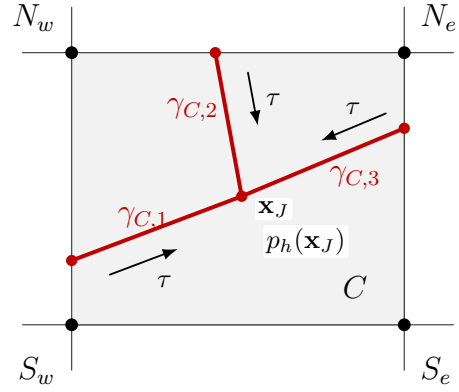

Let $\mathcal A_C^f$ denote the set of fracture segments on a cell $C$.
Consider a fracture segment $\gamma_{C,\alpha} = \overline{\mathbf a_{C,\alpha}\mathbf b_{C,\alpha}} \subset \Gamma_f\cap C$  of length $\ell_{C,\alpha} = \left|\mathbf b_{C,\alpha}-\mathbf a_{C,\alpha}\right|$ for $\alpha\in\mathcal{A}_C^f$.
The local constructions for a fracture chord and a fracture junction are illustrated in Figure~\ref{fig:fracture_local_constructions}.
Let $\boldsymbol\lambda_C(\mathbf x)\in\mathbb R^4$ denote the vector of bilinear Lagrange basis functions on $C$. The pressure at a point $\mathbf x\in C$ is approximated by $p_h(\mathbf x) = \boldsymbol\lambda_C(\mathbf x)^{\mathsf T}\mathbf p_C$.
Define
\[
\mathbf g_{C,\alpha}
:=
\boldsymbol\lambda_C(\mathbf b_{C,\alpha})
-
\boldsymbol\lambda_C(\mathbf a_{C,\alpha}).
\]
Then
\[
p_h(\mathbf b_{C,\alpha})
-
p_h(\mathbf a_{C,\alpha})
=
\mathbf g_{C,\alpha}^{\mathsf T}\mathbf p_C.
\]

Approximating the tangential derivative by the pressure difference along the segment gives 
\[
\partial_\tau p
\approx
\frac{\mathbf g_{C,\alpha}^{\mathsf T}\mathbf p_C}
{\ell_{C,\alpha}}.
\]
Let
\begin{equation}
\label{eq:effective_fracture_conductivity}
\kappa_{C,\alpha}
:=
\frac{1}{\ell_{C,\alpha}}
\int_{\gamma_{C,\alpha}}a_fk_f\,ds
\end{equation}
denote the average fracture conductivity on the segment. 
The corresponding discrete fracture energy is
\begin{equation}
\label{eq:local_fracture_energy}
\mathcal E_{h,C,\alpha}^{f}(\mathbf p_C)
=
\frac12
\frac{\kappa_{C,\alpha}}{\ell_{C,\alpha}}
\left(
\mathbf g_{C,\alpha}^{\mathsf T}\mathbf p_C
\right)^2.
\end{equation}
Equivalently, the segment contributes the rank-one positive
semidefinite matrix
\begin{equation}
\label{eq:local_fracture_matrix}
F_{C,\alpha}
=
\frac{\kappa_{C,\alpha}}{\ell_{C,\alpha}}
\mathbf g_{C,\alpha}\mathbf g_{C,\alpha}^{\mathsf T}
\end{equation}
to the four-node cell matrix $S_C$.
If a cell contains several fracture segments, their energies are added:
\begin{equation}
\label{eq:cell_fracture_energy}
\mathcal E_{h,C}^{f}(\mathbf p_C)
=
\sum_{\alpha\in\mathcal A_C^f}
\mathcal E_{h,C,\alpha}^{f}(\mathbf p_C)
=
\frac12\mathbf p_C^{\mathsf T}F_C\mathbf p_C,
\end{equation}
where
\begin{equation}
\label{eq:cell_fracture_matrix}
F_C
=
\sum_{\alpha\in\mathcal A_C^f}
\frac{\kappa_{C,\alpha}}{\ell_{C,\alpha}}
\mathbf g_{C,\alpha}\mathbf g_{C,\alpha}^{\mathsf T}.
\end{equation}
The same construction applies at fracture junctions. Each fracture is split at its junction points, so that every incident arm has the junction as one of its segment endpoints. 
A free fracture tip (terminal fracture segment) is discarded.

The nodal basis functions form a partition of unity, hence $\mathbf g_{C,\alpha}^{\mathsf T}\mathbf1=0$.
Consequently, $F_C\mathbf1=\mathbf0$, and the fracture contribution preserves the constant-pressure mode.
Moreover, $F_C$ is symmetric positive semidefinite, so adding the fracture energy preserves the symmetry and positive semidefiniteness of every local contribution.
The fracture discretization introduces no additional global unknowns and modifies only the four-node matrix of the intersected cells. 
In the scalar case, these blocks may add diagonal couplings to the five-point matrix stencil; for the full-tensor discretization, they remain within the existing nine-point cell stencil.

Let $\mathcal T_h^f$ denote the set of cells intersected by conductive fractures. 
The global fracture energy is
\begin{equation}
\label{eq:global_fracture_energy}
\mathcal E_h^f(\mathbf p)
=
\frac12
\sum_{C\in\mathcal T_h^f}
\mathbf p_C^{\mathsf T}F_C\mathbf p_C.
\end{equation}
We next combine these contributions into the global discrete problem and describe the recovery of the eliminated barrier jumps.

\subsection{Global discrete problem and local recovery}
\label{sec:global_scheme}

The local barrier and fracture contributions derived above are now assembled into the global nodal problem, after which the eliminated barrier jump variables can be recovered by inexpensive local post-processing when needed.

\subsubsection{Global discrete problem}
\label{sec:global_discrete_problem}

We now combine the barrier-modified matrix energy and the fracture energy into a single discrete problem. 
Let $\widetilde{\mathcal E}_h^{m,b}$ denote the matrix energy after the local barrier modifications. 
In the scalar edgewise formulation,
\[
\widetilde{\mathcal E}_h^{m,b}(\mathbf p)
=
\frac12
\sum_{e=(P,Q)\in\mathcal E_h}
\widehat c_e
(p_P-p_Q)^2,
\]
where $\widehat c_e=c_e$ on an uncut edge and
$\widehat c_e=\widetilde c_e$ on a cut edge; see \eqref{eq:effective_barrier_edge_conductance} and \eqref{eq:multiple_barrier_edge_conductance}. In the cellwise
full-tensor formulation,
\[
\widetilde{\mathcal E}_h^{m,b}(\mathbf p)
=
\frac12
\sum_{C\in\mathcal T_h}
\mathbf p_C^{\mathsf T}
\widehat S_C
\mathbf p_C,
\]
where $\widehat S_C=S_C$ for an uncut cell and $\widehat S_C=\widetilde S_C$ for a barrier-cut cell; see \eqref{eq:modified_tensor_cell_matrix} and \eqref{eq:modified_multiregion_cell_matrix}.

Let $\ell_h(\mathbf p)=\mathbf b_h^{\mathsf T}\mathbf p$ denote the discrete load functional, including the volume source and prescribed
Neumann fluxes. 
The total discrete energy is
\begin{equation}
\label{eq:total_discrete_energy}
\mathcal J_h(\mathbf p)
=
\widetilde{\mathcal E}_h^{m,b}(\mathbf p)
+
\mathcal E_h^f(\mathbf p)
-
\ell_h(\mathbf p).
\end{equation}
If $\mathcal N_{h,D}$ denotes the set of Dirichlet nodes, define
\[
V_{h,D}
=
\left\{
\mathbf p\in\mathbb R^{|\mathcal N_h|}:
p_P=g_D(\mathbf x_P)
\ \text{for every }P\in\mathcal N_{h,D}
\right\}.
\]
The finite difference solution is the minimizer
\begin{equation}
\label{eq:global_discrete_minimization}
\mathbf p_h
=
\operatorname*{argmin}_{\mathbf p\in V_{h,D}}
\mathcal J_h(\mathbf p).
\end{equation}

The stationarity condition for \eqref{eq:global_discrete_minimization} gives a finite difference equation at every free grid node. In the scalar formulation, it takes the explicit form
\begin{equation}
\label{eq:global_scalar_fd_scheme}
\sum_{Q:\,(P,Q)\in\mathcal E_h}
\widehat c_{PQ}
\left(p_P-p_Q\right)
+
\sum_{\substack{C\in\mathcal T_h^f\\P\in C}}
\left(F_C\mathbf p_C\right)_P
=
(\mathbf b_h)_P,
\qquad
P\in\mathcal N_h\setminus\mathcal N_{h,D}.
\end{equation}
Thus, relative to the interface-free five-point scheme, the conductance $c_{PQ}$ is replaced by the effective conductance $\widehat c_{PQ}$ only on barrier-cut edges, while every fracture-cut cell contributes the additional local term $F_C\mathbf p_C$. All other remain unchanged.

For the full-tensor formulation, the corresponding nodal equation is
\begin{equation}
\label{eq:global_tensor_fd_scheme}
\sum_{\substack{C\in\mathcal T_h\\P\in C}}
\left[
\left(\widehat S_C+F_C\right)\mathbf p_C
\right]_P
=
(\mathbf b_h)_P,
\qquad
P\in\mathcal N_h\setminus\mathcal N_{h,D},
\end{equation}
where $F_C=\mathbf0$ if $C$ is not intersected by a conductive fracture. 
Hence the interface-free cell matrix $S_C$ is replaced by $\widehat S_C$ only in barrier-cut cells, and the fracture matrix $F_C$ is added only in fracture-cut cells. The scheme is therefore identical to the standard Cartesian discretization away from the interfaces.

After assembly, both formulations can be written compactly as
\begin{equation}
\label{eq:global_discrete_system}
A_h\mathbf p_h=\mathbf b_h.
\end{equation}
If the nodal unknowns are divided into free (F) and Dirichlet (D) components,
then
\[
(A_h)_{FF}\mathbf p_F
=
(\mathbf b_h)_F
-
(A_h)_{FD}\mathbf p_D,
\qquad
\mathbf p_D
=
\bigl(g_D(\mathbf x_P)\bigr)_{P\in\mathcal N_{h,D}}.
\]
The matrix $A_h$ is symmetric. 
If the problem contains at least one Dirichlet node, the reduced matrix $(A_h)_{FF}$ is positive definite.

All barrier jump variables have been eliminated before global
assembly, and the fracture contributions use only the original nodal
pressures. Consequently, the global unknowns remain exactly those of
the interface-free finite difference scheme. In the scalar
barrier-only case, the five-point graph structure is preserved.
Fracture contributions may additionally couple diagonally opposite
vertices of an intersected cell, producing at most a nine-point
stencil, while the full-tensor discretization already possesses this
nine-point structure.

\subsubsection{Local recovery}
\label{sec:local_recovery}

The nodal vector $\mathbf p_h$ is the primary finite difference solution. 
Quantities associated with the interfaces can be recovered after the global solve.
These recovery procedures are optional post-processing steps and do not alter the global nodal solution $\mathbf p_h$.

For the scalar edgewise barrier formulation, let $\mathcal A_e$ denote the set of barrier crossings of an oriented edge $e=(P,Q)$. 
The effective conductance of the edge is
\[
\widetilde c_e
=
\left(
\frac{1}{c_e}
+
\sum_{\beta\in\mathcal A_e}
\frac{1}{d_{e,\beta}}
\right)^{-1},
\]
and the corresponding oriented matrix--barrier flux is
\[
q_e
=
\widetilde c_e\,\delta_e\mathbf p_h.
\]
The jump variable associated with crossing
$\alpha\in\mathcal A_e$ is therefore recovered as
\begin{equation}
\label{eq:recovered_edge_barrier_jump}
t_{e,\alpha}^\star
=
\frac{q_e}{d_{e,\alpha}}
=
\frac{\widetilde c_e}{d_{e,\alpha}}
\,\delta_e\mathbf p_h.
\end{equation}
With the
orientation convention introduced in
Section~\ref{sec:scalar_barrier_discretization}, the physical pressure
jump at the crossing is approximated by $\llbracket p\rrbracket(\mathbf x_{e,\alpha})
\approx
\sigma_{e,\alpha}t_{e,\alpha}^\star$.

For a full-tensor cell containing a single barrier chord, the
eliminated jump variable is recovered from
\begin{equation}
\label{eq:recovered_single_cell_barrier_jump}
t_C^\star
=
\frac{
\mathbf r_C^{\mathsf T}S_C\mathbf p_{h,C}
}{
\mathbf r_C^{\mathsf T}S_C\mathbf r_C+d_C
},
\end{equation}
where $\mathbf p_{h,C}$ is the restriction of the nodal solution to the four vertices of $C$.
The recovered $t_C^\star$ is the pressure offset of the vertex group with $(\mathbf r_C)_i=1$ relative to the group with $(\mathbf r_C)_i=0$. 
The physical jump is $\llbracket p\rrbracket\approx\sigma_Ct_C^\star$, where $\sigma_C=1$ if the group with $(\mathbf r_C)_i=1$ lies on the $-$ side of the chord and $\sigma_C=-1$ otherwise.
For a cell divided into several subregions, the regionwise pressure offsets are
recovered by
\begin{equation}
\label{eq:recovered_multiregion_barrier_offsets}
\mathbf t_C^\star
=
\left(
E_C^{\mathsf T}S_CE_C+P_C
\right)^{-1}
E_C^{\mathsf T}S_C\mathbf p_{h,C}.
\end{equation}
These formulas are precisely the local back-substitutions associated
with the barrier variables eliminated during the construction of the
condensed cell matrices.

The tangential flux along each conductive fracture segment is recovered
directly from the nodal pressure. Orient
$\gamma_{C,\alpha}$ from $\mathbf a_{C,\alpha}$ to
$\mathbf b_{C,\alpha}$ and let $\boldsymbol\tau_{C,\alpha}
=
\frac{
\mathbf b_{C,\alpha}-\mathbf a_{C,\alpha}
}{
\ell_{C,\alpha}
}.$
The segmentwise tangential pressure gradient is approximated by $\partial_\tau p
\approx
\frac{
\mathbf g_{C,\alpha}^{\mathsf T}\mathbf p_{h,C}
}{
\ell_{C,\alpha}
}.$
Consequently, the tangential fracture flux in the direction
$\boldsymbol\tau_{C,\alpha}$ is recovered as
\begin{equation}
\label{eq:recovered_fracture_flux}
q_{C,\alpha}^{f}
=
-
\frac{\kappa_{C,\alpha}}{\ell_{C,\alpha}}
\mathbf g_{C,\alpha}^{\mathsf T}\mathbf p_{h,C}.
\end{equation}
A positive value of $q_{C,\alpha}^{f}$ represents flow from
$\mathbf a_{C,\alpha}$ toward $\mathbf b_{C,\alpha}$. At a fracture
junction, the incident arms may be oriented outward from the junction
so that their recovered outgoing fluxes are reported with a common sign
convention.

The nodal solution $\mathbf p_h$ can be used to define a standard bilinear interpolant on each Cartesian cell and hence a continuous piecewise bilinear function over the computational grid.
However, this interpolant cannot represent the pressure jumps across $\Gamma_b$. 
To recover these discontinuities, we apply the cellwise local recovery (even when $\mathbf p_h$ is obtained from the scalar edgewise formulation) to obtain the region offsets $\mathbf t_C^\star$, and define the reconstruction
\begin{equation}
\label{eq:broken_pressure_reconstruction}
p_h|_{C_j}
=
I_C (\mathbf p_{h,C}-E_C\mathbf t_C^\star)+t_{C,j}^\star,
\qquad
j=0,\ldots,m,
\end{equation}
where $I_C$ denotes the bilinear interpolation operator on $C$.
Thus, the offsets are removed before interpolation and restored separately on each subregion.
The resulting broken piecewise bilinear field ${p}_{h}$ interpolates the nodal solution $\mathbf{p}_h$ at all grid points and explicitly represents the recovered pressure jumps across the barriers.
This construction is entirely local and requires no additional global solve.

\section{Numerical results}
\label{sec:numerical_results}

We now examine the performance of the proposed schemes in a sequence of two- and three-dimensional numerical experiments.  
Throughout this section the computational grids are uniformly Cartesian, and $N$ stands for the number of cells in each coordinate direction.
Problems with scalar matrix permeability are solved with the edgewise formulation of Sections~\ref{sec:scalar_matrix_energy} and~\ref{sec:scalar_barrier_discretization}, and problems with a full permeability tensor are solved with the cellwise formulation of Sections~\ref{sec:tensor_matrix_energy} and~\ref{sec:tensor_barrier_discretization}.
The pressure in barrier-cut cells is reconstructed using \eqref{eq:broken_pressure_reconstruction} so that the pressure jump appears exactly at the interface rather than being smeared across the cut cells.

\subsection{Example 1: convergence study}
\label{sec:example1}

The first experiment measures the convergence of the method on four
single-interface problems whose exact solutions are available in closed form \cite{xu2020hybrid, liu2026high}.
Cases~(a) and~(b) concern a homogeneous isotropic matrix, $\mathbf K_m=\mathbf I$;
cases~(c) and~(d) repeat the same two interface types with the heterogeneous
anisotropic field
\begin{equation}
\label{eq:example1_tensor}
\mathbf K_m(x,y)=k(x,y)\,\mathbf A,
\qquad
k(x,y)=2+\sin(\pi x)\sin(\pi y),
\qquad
\mathbf A=
\begin{pmatrix}
\alpha & \beta\\
\beta & \gamma
\end{pmatrix},
\end{equation}
where $\alpha=1-\eta/4$, $\beta=\tfrac{\sqrt3}{4}\eta$, $\gamma=1-3\eta/4$, and $\eta=0.9$.  
The eigenvalues of $\mathbf A$ are $1$ and $1-\eta=0.1$, with eigen-directions rotated by $\pi/6$ relative to the coordinate axes, so the grid is far from aligned with the anisotropy.

\emph{Case (a): fracture, scalar matrix permeability.}
On $\Omega=(-1,1)^2$, a conductive fracture with $a_fk_f=2$ occupies the segment $y=0$.  
The exact pressure and source are
\begin{equation*}
\label{eq:example1_case_a}
p(x,y)=\sin(x)\,e^{|y|},\qquad f=0.
\end{equation*}

\emph{Case (b): barrier, scalar matrix permeability.}
On $\Omega=(0,1)^2$, a low-permeability barrier with $a_b/k_b=1$ occupies the segment $x=\tfrac12$.
The exact pressure and source are
\begin{equation*}
\label{eq:example1_case_b}
p(x,y)=\sin(x)\sin(y)+\cos(\tfrac12)\sin(y)\,\chi_{\{x>1/2\}},
\quad
f(x,y)=2\sin(x)\sin(y)+\cos(\tfrac12)\sin(y)\,\chi_{\{x>1/2\}},
\end{equation*}
where $\chi$ denotes the indicator function. 

\emph{Case (c): fracture, tensor matrix permeability.}
The geometry and fracture conductivity are those of case~(a), combined with the permeability field \eqref{eq:example1_tensor}.  Setting
$\lambda=1/(2\gamma)$, the exact pressure is
\begin{equation*}
\label{eq:example1_case_c}
p(x,y)=\sin(x)\,e^{\lambda|y|}.
\end{equation*}

\emph{Case (d): barrier, tensor matrix permeability.}
The geometry and barrier resistance are those of case~(b), combined with the field \eqref{eq:example1_tensor}.  
Setting $J(y)=k(\tfrac12,y)
\left[
\alpha\cos(\tfrac12)\sin(y)+\beta\sin(\tfrac12)\cos(y)
\right]$,
the exact pressure is
\begin{equation*}
\label{eq:example1_case_d}
p(x,y)=\sin(x)\sin(y)
+\Big[J(y)-\frac{\beta}{\alpha}\big(x-\tfrac12\big)J'(y)\Big]
\chi_{\{x>1/2\}}.
\end{equation*}

\begin{figure}[htbp!]
\centering
\begin{subfigure}[b]{0.4\textwidth}
\includegraphics[width=\textwidth]{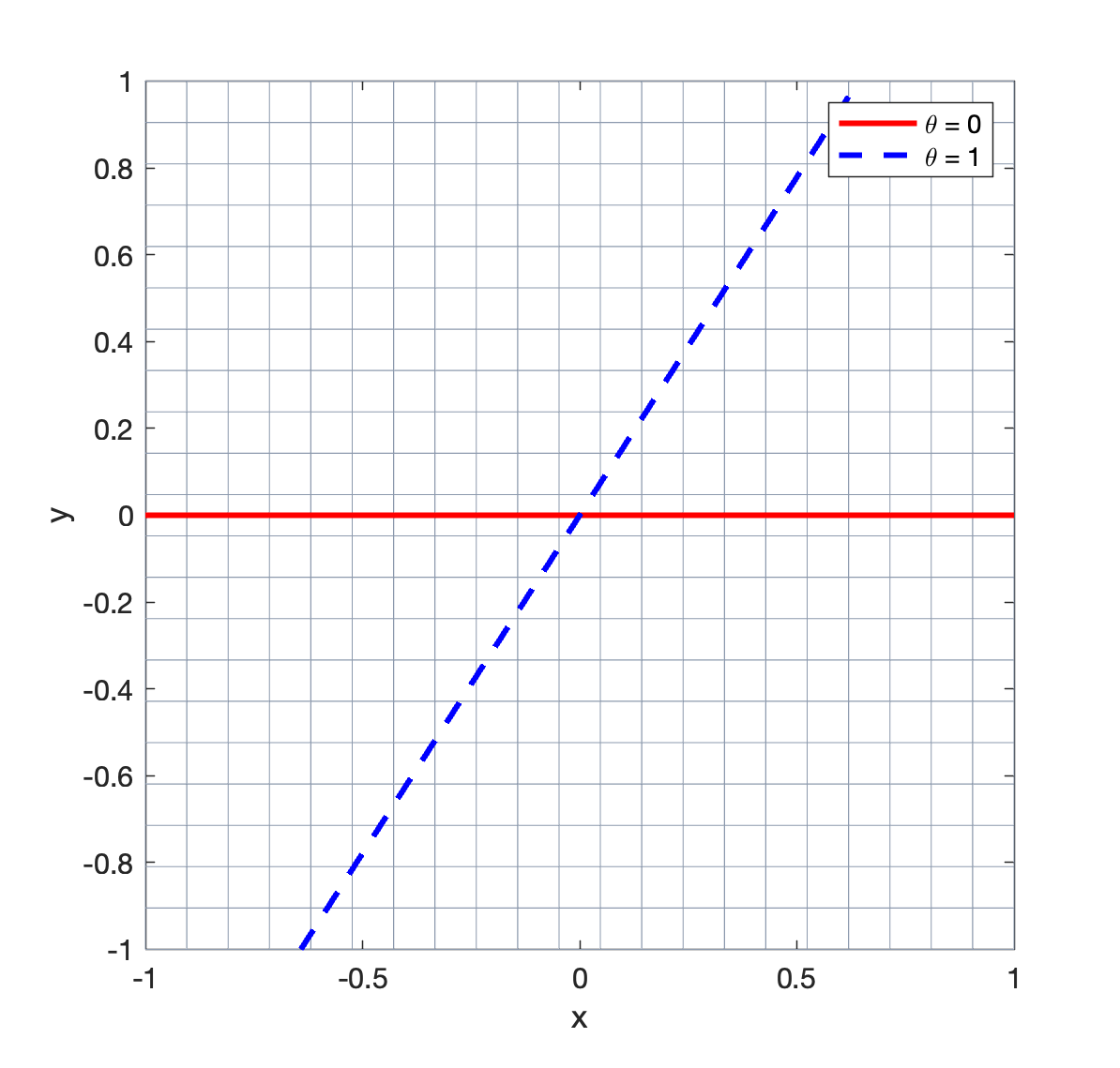}
\caption{fracture, cases (a) and (c)}
\end{subfigure}
\begin{subfigure}[b]{0.4\textwidth}
\includegraphics[width=\textwidth]{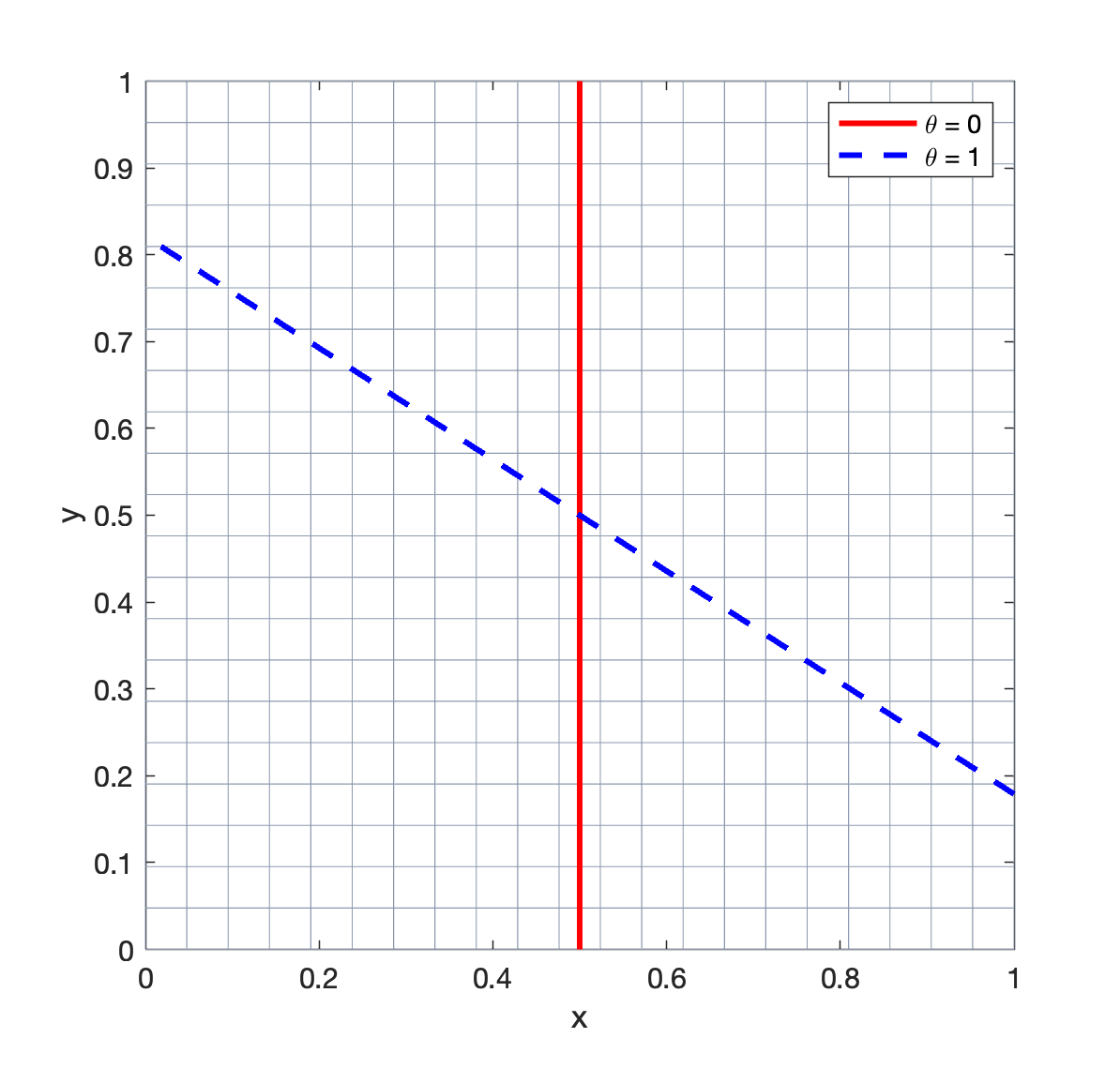}
\caption{barrier, cases (b) and (d)}
\end{subfigure}
\caption{Example 1, interface configurations on the coarsest grid ($N=21$).
Solid line: axis-aligned interface ($\theta=0$). Dashed line: interface rotated by one radian about the domain center ($\theta=1$).  The grids are identical in
both configurations; only the problem is rotated.}
\label{fig:example1_config}
\end{figure}

In cases~(c) and~(d) the source $f=-\nabla\cdot(\mathbf K_m\nabla p)$ is evaluated from its closed-form expression on each side of the interface. 
The formulas are lengthy and are omitted here.  
In all four cases the Dirichlet data on the whole of $\partial\Omega$ is taken from the exact solution.
Every case is run in two geometric configurations, shown in Figure~\ref{fig:example1_config}.  
In the axis-aligned configuration ($\theta=0$), the interface is the horizontal or vertical segment given above.
Since the cell counts $N$ are odd, the interface then passes midway between two rows or columns of grid lines.
In the rotated configuration ($\theta=1$), the entire problem, including the interface, exact solution, source, and $\mathbf K_m$, is rotated by one radian about the center of the domain.  
The rotated problem admits the rotated exact solution.
The grids stay fixed and the interface therefore crosses the cells obliquely.

\begin{table}[htbp]\centering\small
\setlength{\tabcolsep}{4.5pt}
\caption{Example 1, discrete $\ell^2$ pressure error and observed order for the
four cases, in the axis-aligned ($\theta=0$) and rotated ($\theta=1$)
configurations.  Upper block: homogeneous isotropic matrix, edgewise five-point scheme.  
Lower block: heterogeneous anisotropic matrix, cellwise nine-point
scheme.}
\label{tab:example1}
\begin{tabular}{r rr rr @{\hskip 1.8em} rr rr}
\toprule
 & \multicolumn{4}{c}{(a) fracture, $\mathbf K_m=\mathbf I$}
 & \multicolumn{4}{c}{(b) barrier, $\mathbf K_m=\mathbf I$}\\
\cmidrule(r){2-5}\cmidrule(l){6-9}
 & \multicolumn{2}{c}{$\theta=0$} & \multicolumn{2}{c}{$\theta=1$}
 & \multicolumn{2}{c}{$\theta=0$} & \multicolumn{2}{c}{$\theta=1$}\\
$N$ & error & order & error & order & error & order & error & order\\
\midrule
21   & 1.30e-3 & --   & 5.36e-3 & --      & 5.35e-6 & --   & 2.26e-3 & --   \\
41   & 6.46e-4 & 1.05 & 2.23e-3 & 1.31    & 1.41e-6 & 2.00 & 9.77e-4 & 1.26 \\
81   & 3.22e-4 & 1.02 & 9.33e-4 & 1.28    & 3.61e-7 & 2.00 & 4.08e-4 & 1.28 \\
161  & 1.61e-4 & 1.01 & 9.35e-4 & $-$0.00 & 9.13e-8 & 2.00 & 1.96e-4 & 1.07 \\
321  & 8.03e-5 & 1.01 & 2.47e-4 & 1.93    & 2.30e-8  & 2.00 & 9.02e-5 & 1.12 \\
641  & 4.02e-5 & 1.00 & 1.87e-4 & 0.40    & 5.76e-9  & 2.00 & 4.16e-5 & 1.12 \\
1281 & 2.01e-5 & 1.00 & 1.14e-4 & 0.71    & 1.45e-9  & 2.00 & 2.03e-5 & 1.04 \\
2561 & 1.00e-5 & 1.00 & 4.79e-5 & 1.26    & 3.77e-10 & 1.94 & 9.91e-6 & 1.03 \\
\midrule
 & \multicolumn{4}{c}{(c) fracture, $\mathbf K_m=k(x,y)\mathbf A$}
 & \multicolumn{4}{c}{(d) barrier, $\mathbf K_m=k(x,y)\mathbf A$}\\
\cmidrule(r){2-5}\cmidrule(l){6-9}
 & \multicolumn{2}{c}{$\theta=0$} & \multicolumn{2}{c}{$\theta=1$}
 & \multicolumn{2}{c}{$\theta=0$} & \multicolumn{2}{c}{$\theta=1$}\\
$N$ & error & order & error & order & error & order & error & order\\
\midrule
21   & 7.63e-3 & --   & 1.08e-2 & --   & 1.14e-2 & --   & 7.74e-3 & --   \\
41   & 2.00e-3 & 2.00 & 4.67e-3 & 1.26 & 6.12e-3 & 0.93 & 3.04e-3 & 1.40 \\
81   & 5.78e-4 & 1.82 & 2.03e-3 & 1.22 & 3.19e-3 & 0.96 & 1.16e-3 & 1.41 \\
161  & 2.27e-4 & 1.36 & 1.89e-3 & 0.11 & 1.63e-3 & 0.98 & 4.51e-4 & 1.38 \\
321  & 1.12e-4 & 1.03 & 6.00e-4 & 1.66 & 8.24e-4 & 0.99 & 1.70e-4 & 1.41 \\
641  & 5.73e-5 & 0.96 & 4.15e-4 & 0.53 & 4.15e-4 & 0.99 & 6.59e-5 & 1.37 \\
1281 & 2.93e-5 & 0.97 & 2.44e-4 & 0.77 & 2.08e-4 & 1.00 & 2.77e-5 & 1.25 \\
2561 & 1.48e-5 & 0.98 & 1.07e-4 & 1.18 & 1.04e-4 & 1.00 & 1.25e-5 & 1.15 \\
\bottomrule
\end{tabular}
\end{table}
Errors of the nodal pressure are measured in the $\ell^2$ norm
\begin{equation}
\label{eq:discrete_l2_norm}
\|v\|_{\ell^2(\mathcal N_h)}
=
\Big(
\sum_{P\in\mathcal N_h}
h_xh_y\,v_P^2
\Big)^{1/2}.
\end{equation}
Table~\ref{tab:example1} collects the discrete $\ell^2$ pressure errors and the observed orders of all cases on the grids $N=21,41,\ldots,2561$.
The axis-aligned barrier under the edgewise scheme (case~(b), $\theta=0$) converges at second order.
The axis-aligned barrier under the cellwise scheme (case~(d), $\theta=0$) is first order instead, because it uses a single constant jump per cut cell to approximate the varying jump.  
All remaining configurations converge at first order as an unfitted treatment.
In the rotated fracture cases, the per-step orders oscillate, since refinement relocates the fracture-gridline crossings. 
For $\theta=1$, cases~(a) and~(c) share the same crossing pattern, which is why their irregular steps coincide, and least-squares fits over all eight refinement levels give overall rates of $0.94$ and $0.92$.

\subsection{Example 2: a single immersed fracture}
\label{sec:example2}

The second experiment \cite{angot2009asymptotic, xu2025extension, liu2026interior} places a single interface in the unit square $\Omega=(0,1)^2$ with $\mathbf K_m=\mathbf I$ and examines the two conductivity extremes.
In the conductive regime, the fracture has the tangential conductivity $a_fk_f=10^{5}$, and the flow is driven across the domain from top to bottom: $p=1$ on the upper boundary, $p=0$ on the lower boundary, and no flow through the other two sides.  
In the blocking regime, the barrier has the resistance $a_b/k_b=10^{5}$, and the flow is driven from right to left: $p=1$ on the right, $p=0$ on the left, and no flow through the top and bottom.
Each regime is computed for two placements of the interface: a vertical segment from $(\tfrac12,\tfrac12)$ to $(\tfrac12,1)$, whose lower end is a free tip inside the domain, and a slanted segment from $(\tfrac14,\tfrac34)$ to $(\tfrac34,\tfrac14)$, both ends of which are free tips.  
The grid is uniform with $N=49$ cells per direction, so the vertical segment runs midway between two grid columns.  
The slanted segment lying on the diagonal is displaced by $h/1000$ in both coordinate (a geometric perturbation of about $2\times10^{-5}$) to avoid the crossing with grid nodes.
Free fracture tips are discarded on their terminating cells.
The results are compared with reference profiles computed by a box discrete fracture model \cite{xu2025extension} on conforming grids of roughly $23{,}000$ cells, sampled along $y=0.75$ for the vertical placement and along $y=0.5$ for the slanted one.

\begin{figure}[htbp!]
\centering
\begin{subfigure}[b]{0.40\textwidth}
\includegraphics[width=\textwidth]{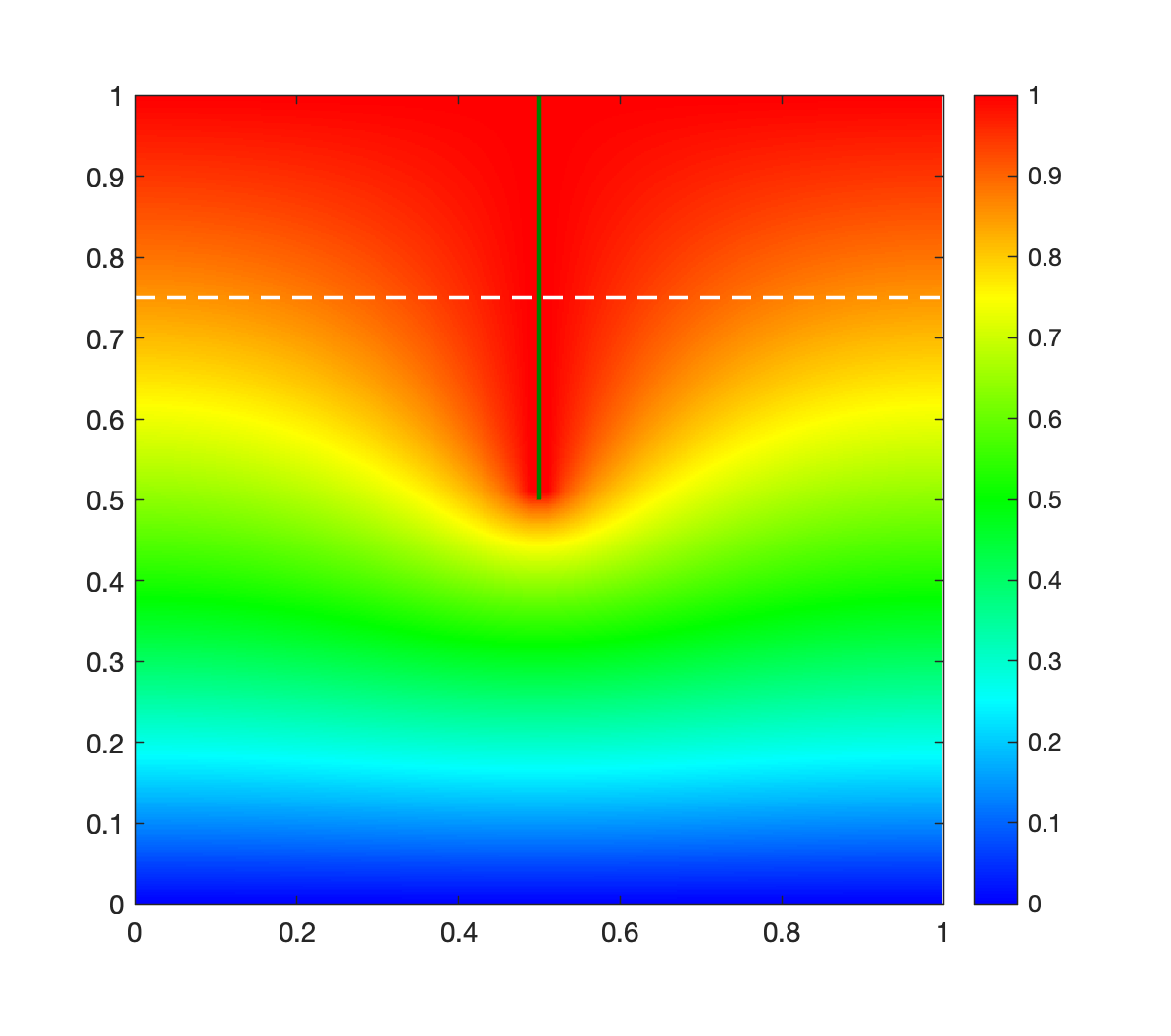}
\caption{vertical: pressure}
\end{subfigure}
\begin{subfigure}[b]{0.40\textwidth}
\includegraphics[width=\textwidth]{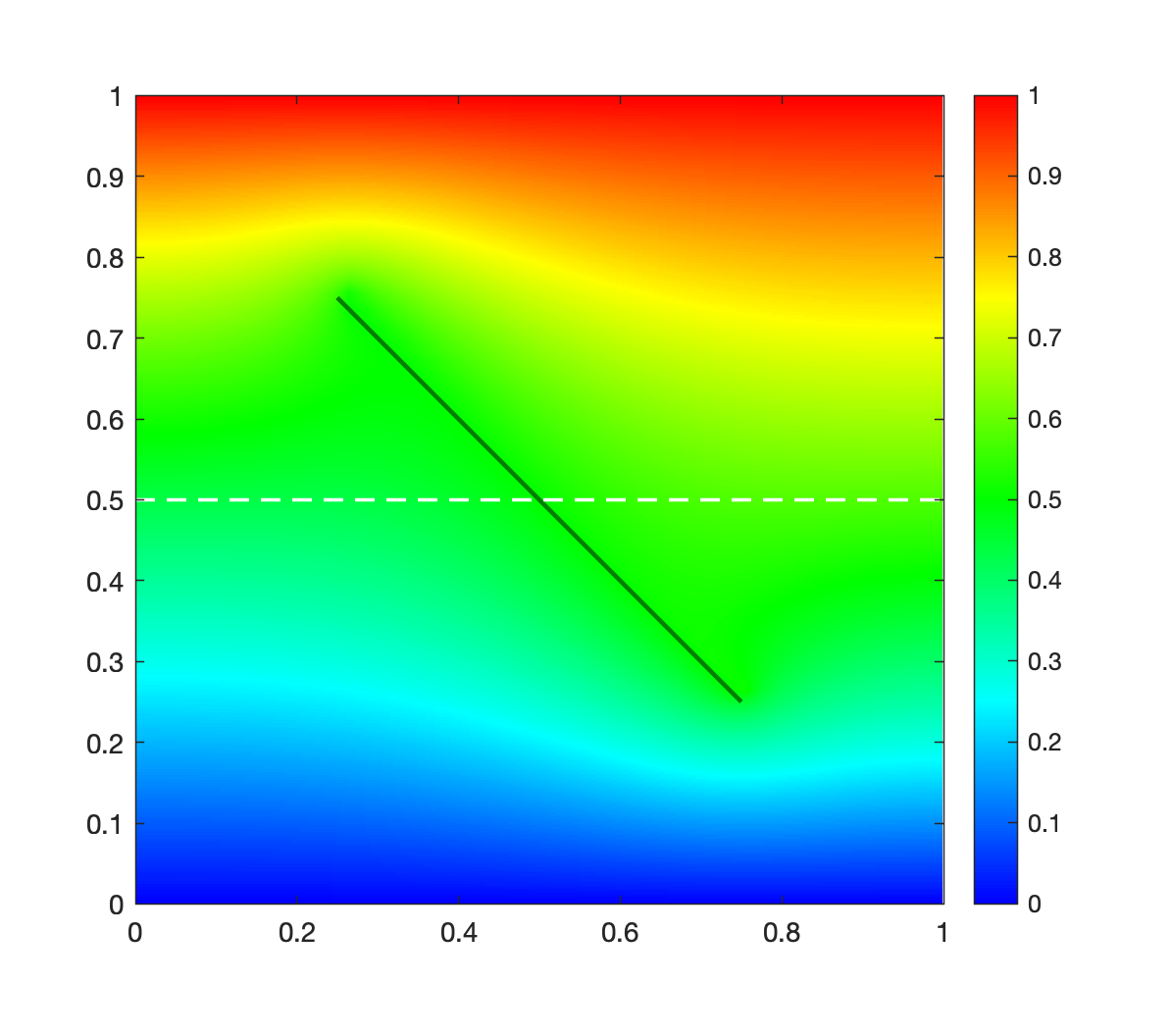}
\caption{slanted: pressure}
\end{subfigure}\\
\begin{subfigure}[b]{0.40\textwidth}
\includegraphics[width=\textwidth]{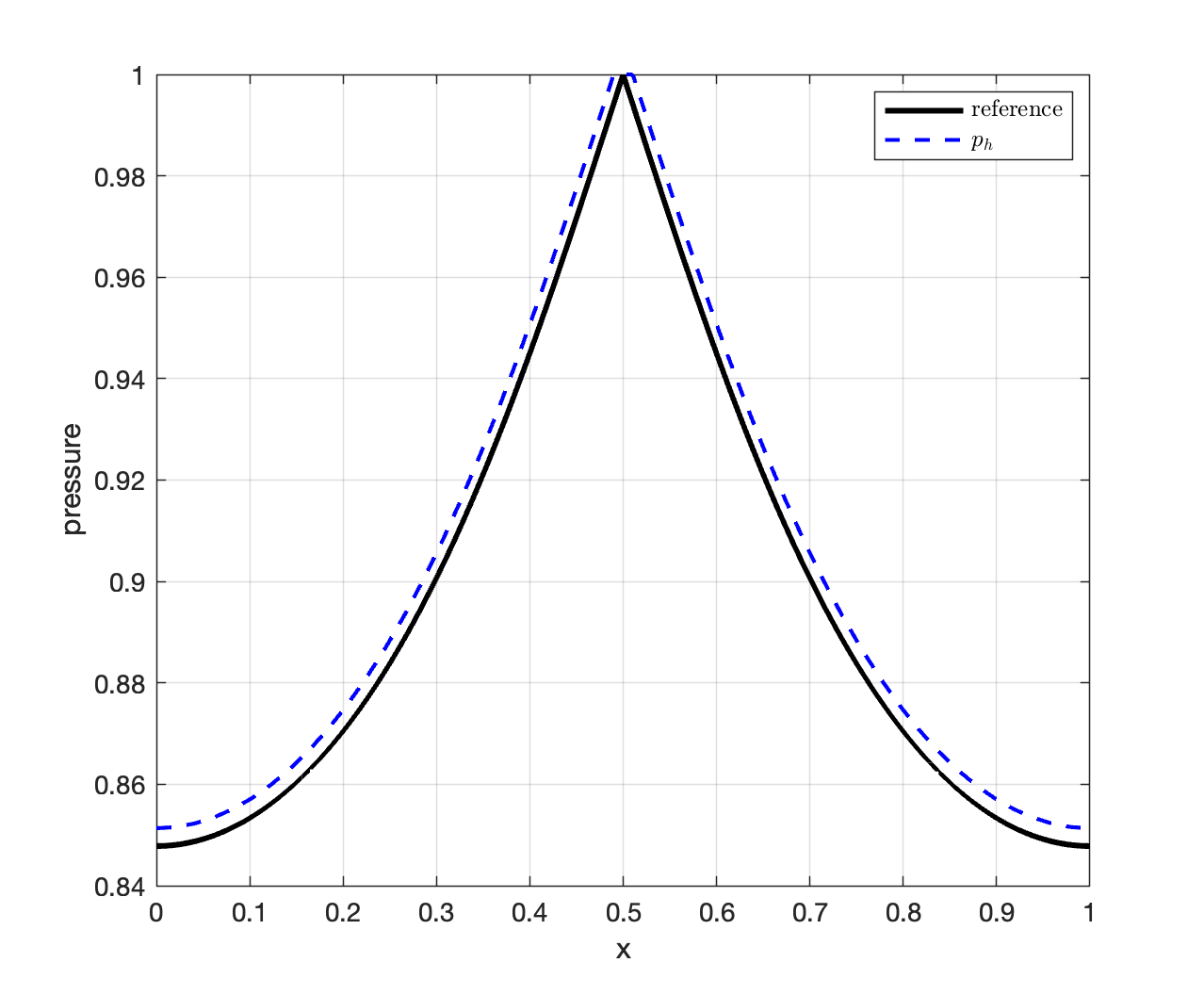}
\caption{vertical: slice $y=0.75$}
\end{subfigure}
\begin{subfigure}[b]{0.40\textwidth}
\includegraphics[width=\textwidth]{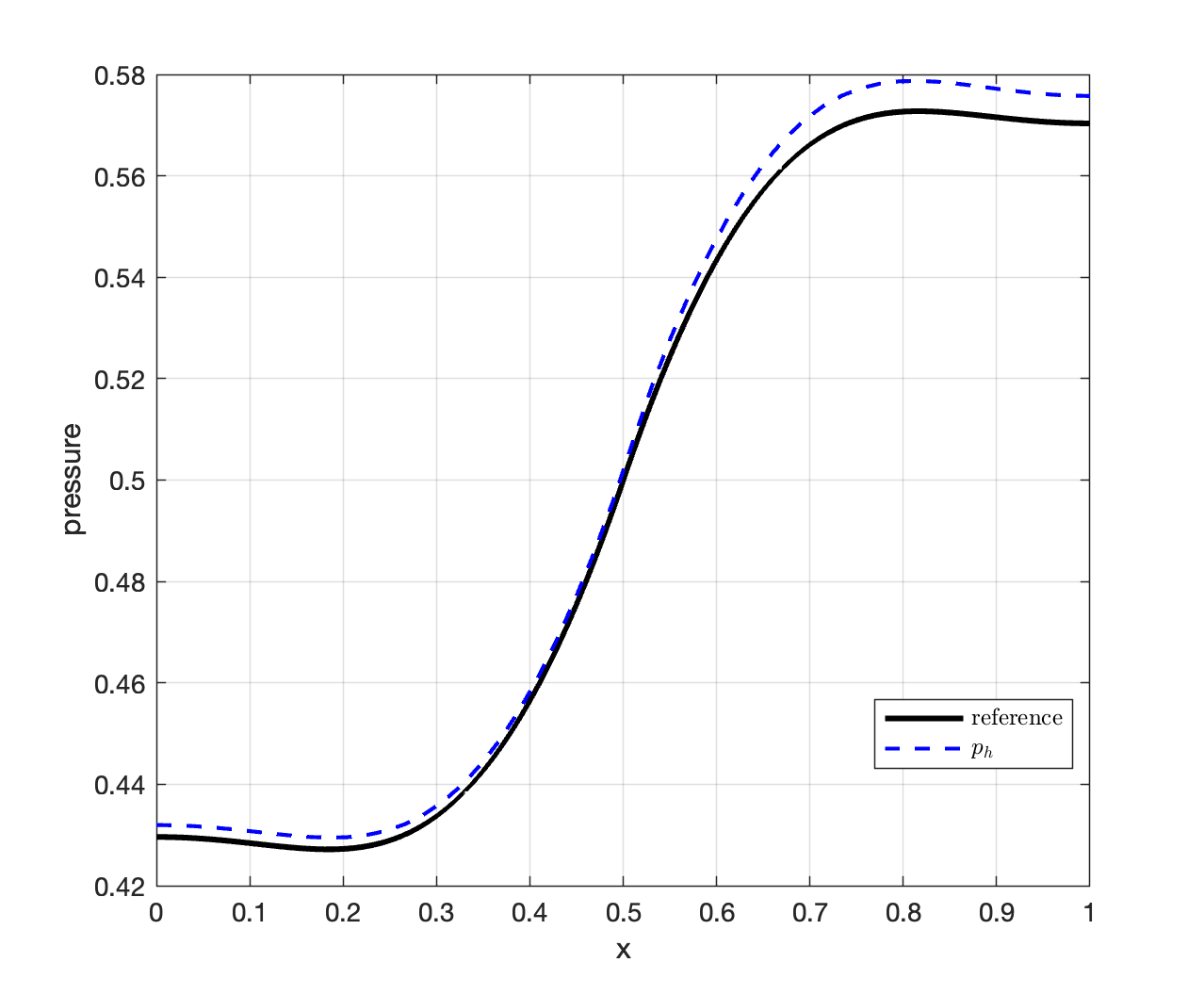}
\caption{slanted: slice $y=0.5$}
\end{subfigure}
\caption{Example 2, conductive regime ($a_fk_f=10^{5}$) on the $N=49$ grid.  
First row: pressure field with the fracture drawn in green and
the sampling slice dashed.  
Second row: pressure profile along the slice, against the conforming box-DFM reference.}
\label{fig:example2_fracture}
\end{figure}

\begin{figure}[htbp!]
\centering
\begin{subfigure}[b]{0.40\textwidth}
\includegraphics[width=\textwidth]{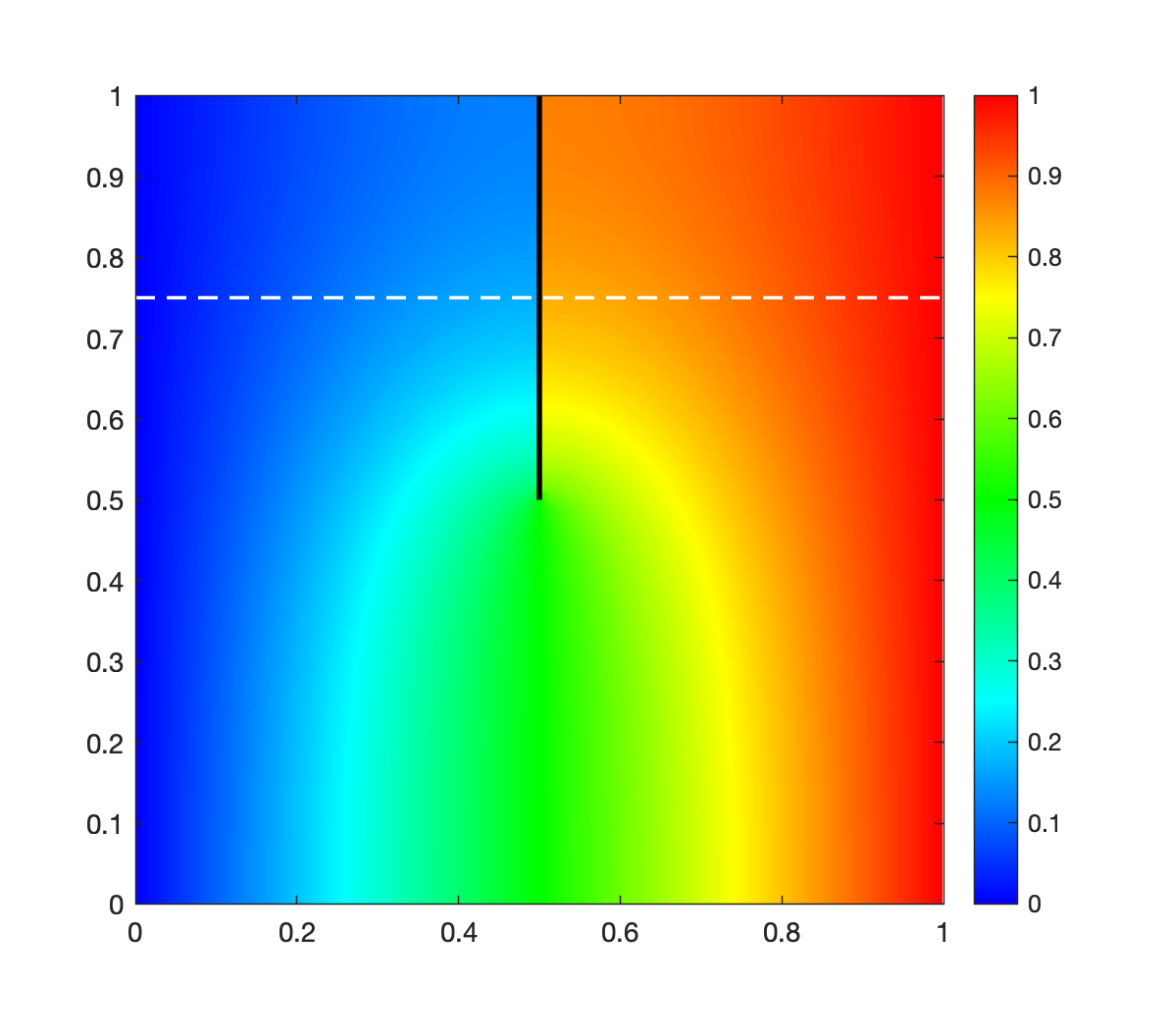}
\caption{vertical: pressure}
\end{subfigure}
\begin{subfigure}[b]{0.40\textwidth}
\includegraphics[width=\textwidth]{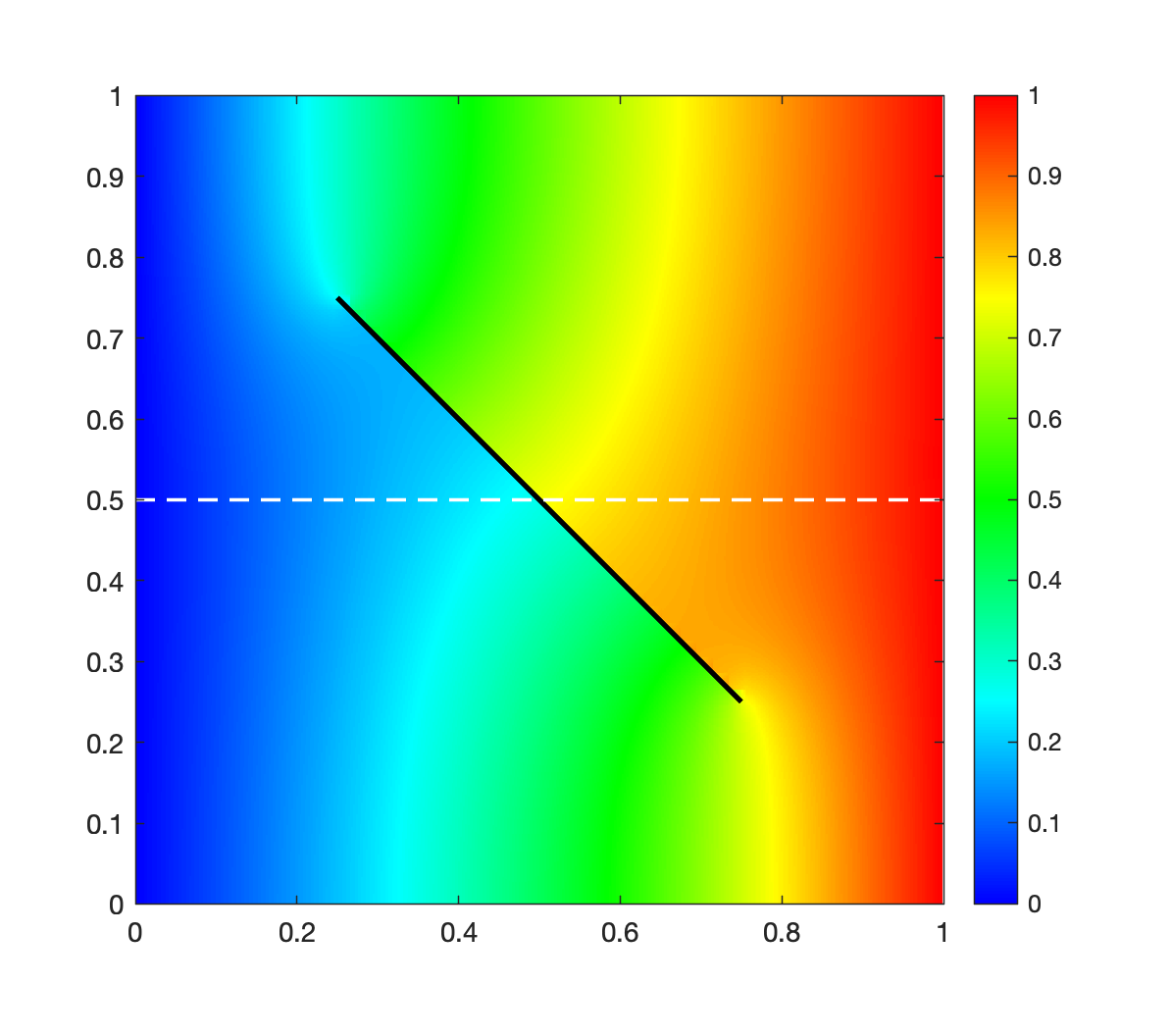}
\caption{slanted: pressure}
\end{subfigure}\\
\begin{subfigure}[b]{0.40\textwidth}
\includegraphics[width=\textwidth]{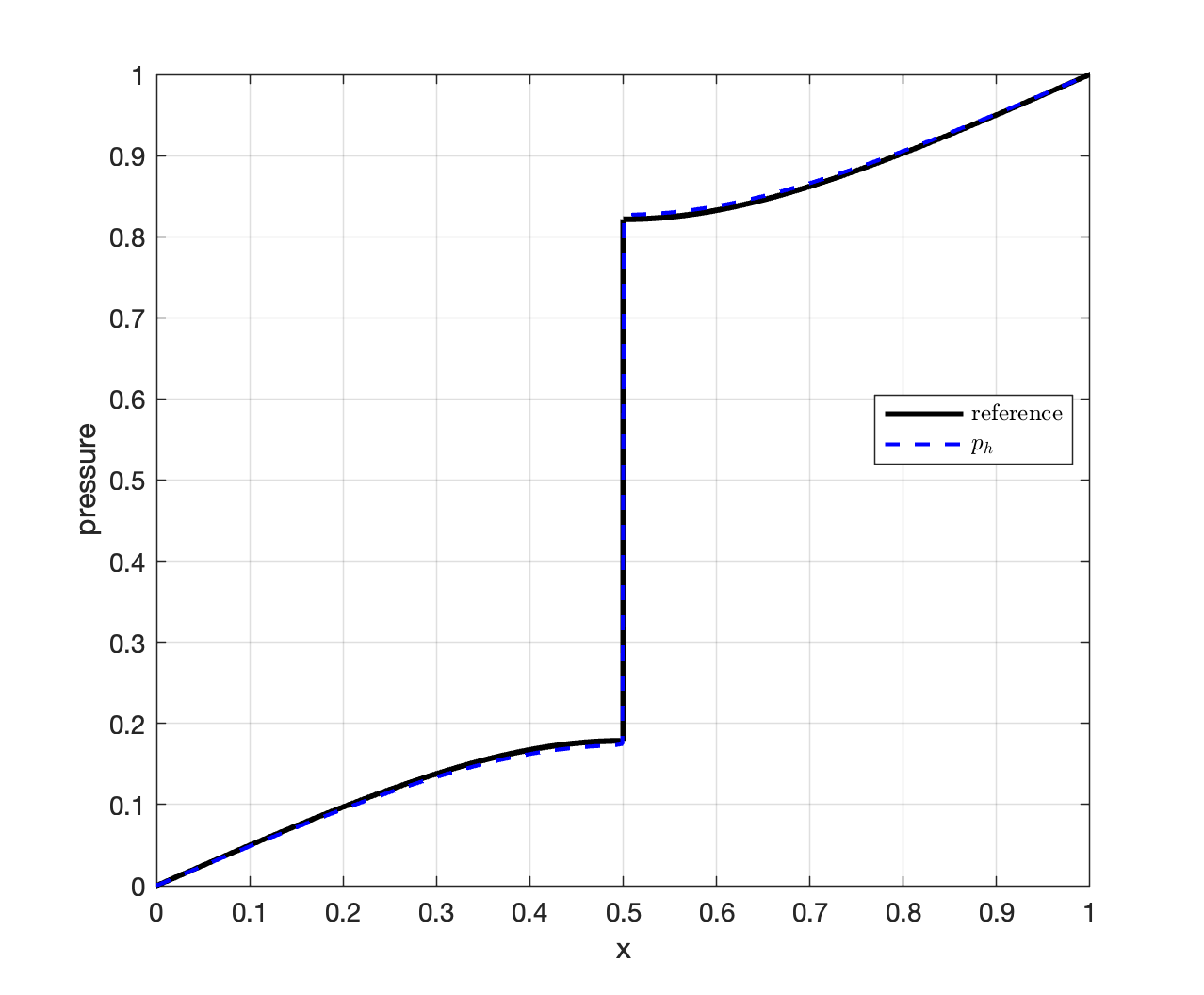}
\caption{vertical: slice $y=0.75$}
\end{subfigure}
\begin{subfigure}[b]{0.40\textwidth}
\includegraphics[width=\textwidth]{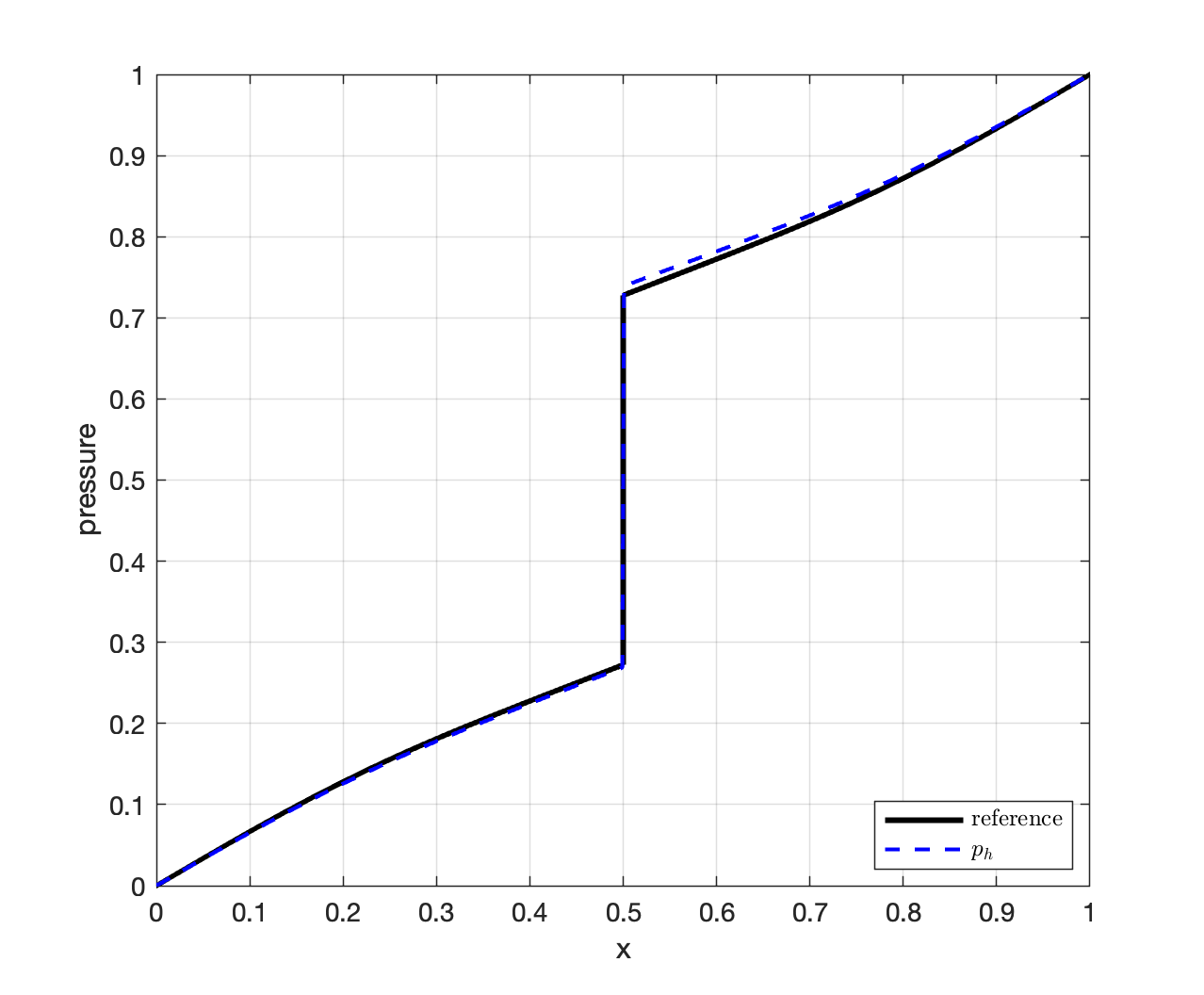}
\caption{slanted: slice $y=0.5$}
\end{subfigure}
\caption{Example 2, blocking regime ($a_b/k_b=10^{5}$) on the $N=49$ grid.  
First row: pressure field with the barrier drawn in black and the sampling slice dashed.  
Second row: pressure profile along the slice, against the conforming box-DFM reference.}
\label{fig:example2_barrier}
\end{figure}

Figure~\ref{fig:example2_fracture} shows the conductive regime, where the fracture short-circuits the pressure drop along its path.
Figure~\ref{fig:example2_barrier} shows the blocking regime, where the barrier splits the fields and the recovered pressure jumps across the interface. 
Both profiles follow the fine-grid reference closely on the slice.

\subsection{Example 3: a regular fracture network}
\label{sec:example3}

We next solve the regular-network benchmark of \cite{flemisch2018benchmarks}.
The unit square with $\mathbf K_m=\mathbf I$ contains six axis-aligned fractures: a vertical and a horizontal segment spanning the whole domain through $0.5$, a half-length horizontal and vertical pair through $0.75$, and a quarter-length pair through $0.625$.  
All fractures share the aperture $a=10^{-4}$.  
In case~(a) they are conductive fractures with $k_f=10^{4}$, hence $a_fk_f=1$; in case~(b) they are low-permeability barriers with $k_b=10^{-4}$, hence $a_b/k_b=1$.  
A unit influx $\mathbf u\cdot\mathbf n=-1$ enters through the left boundary, the right boundary carries the Dirichlet value $p=1$, and the top and bottom are impermeable.  
The computation uses a uniform grid with $N=49$ cells per direction, on which the coordinates $0.5$, $0.625$ and $0.75$ all fall between grid lines. 
Junctions of conductive fractures are coupled through their shared interpolated pressures as described in Section~\ref{sec:fracture_discretization}.  
The reference is a mimetic finite difference solution of the benchmark, sampled along $y=0.7$ and $x=0.5$ for the fractures and along the
segment from $(0,0.1)$ to $(0.9,1)$ for the barriers.

\begin{figure}[htbp!]
\centering
\begin{subfigure}[b]{0.4\textwidth}
\includegraphics[width=\textwidth]{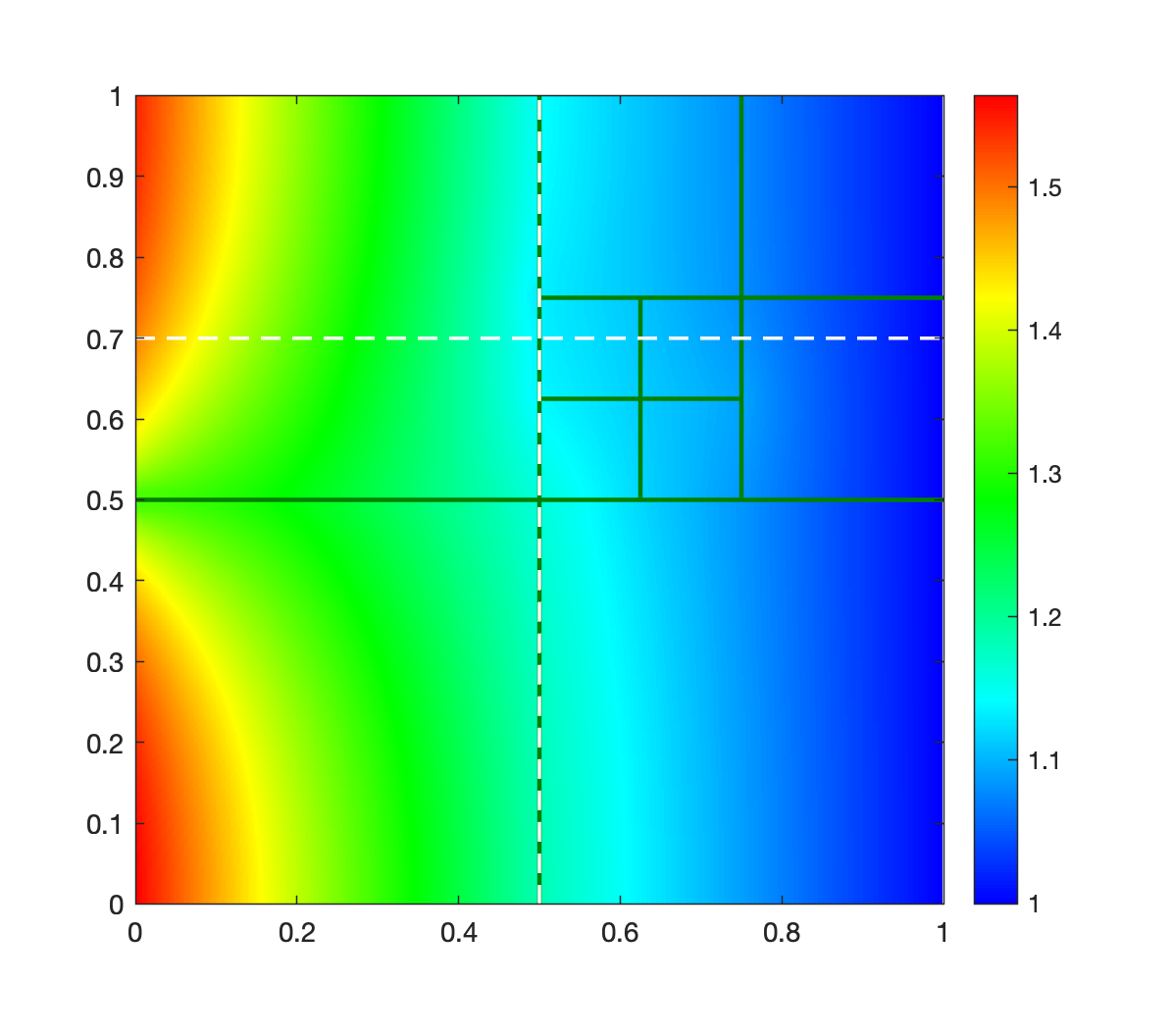}
\caption{fracture network}
\end{subfigure}
\begin{subfigure}[b]{0.4\textwidth}
\includegraphics[width=\textwidth]{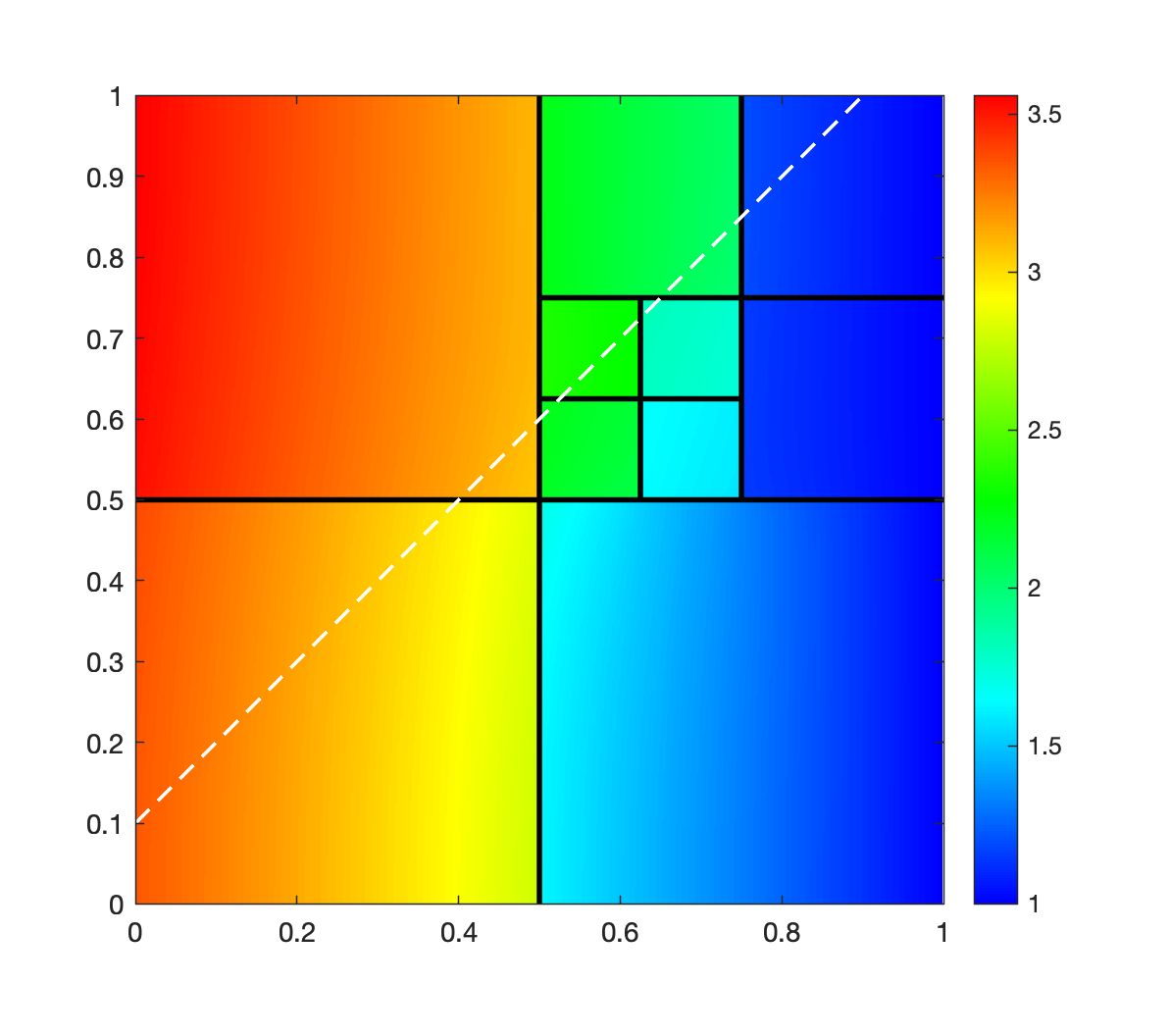}
\caption{barrier network}
\end{subfigure}
\caption{Example 3, pressure fields on the $N=49$ grid.
Left: the conductive fracture network.
Right: the blocking barrier network. 
Dashed lines mark the sampling slices.}
\label{fig:example3_fields}
\end{figure}

\begin{figure}[htbp!]
\centering
\begin{subfigure}[b]{0.32\textwidth}
\includegraphics[width=\textwidth]{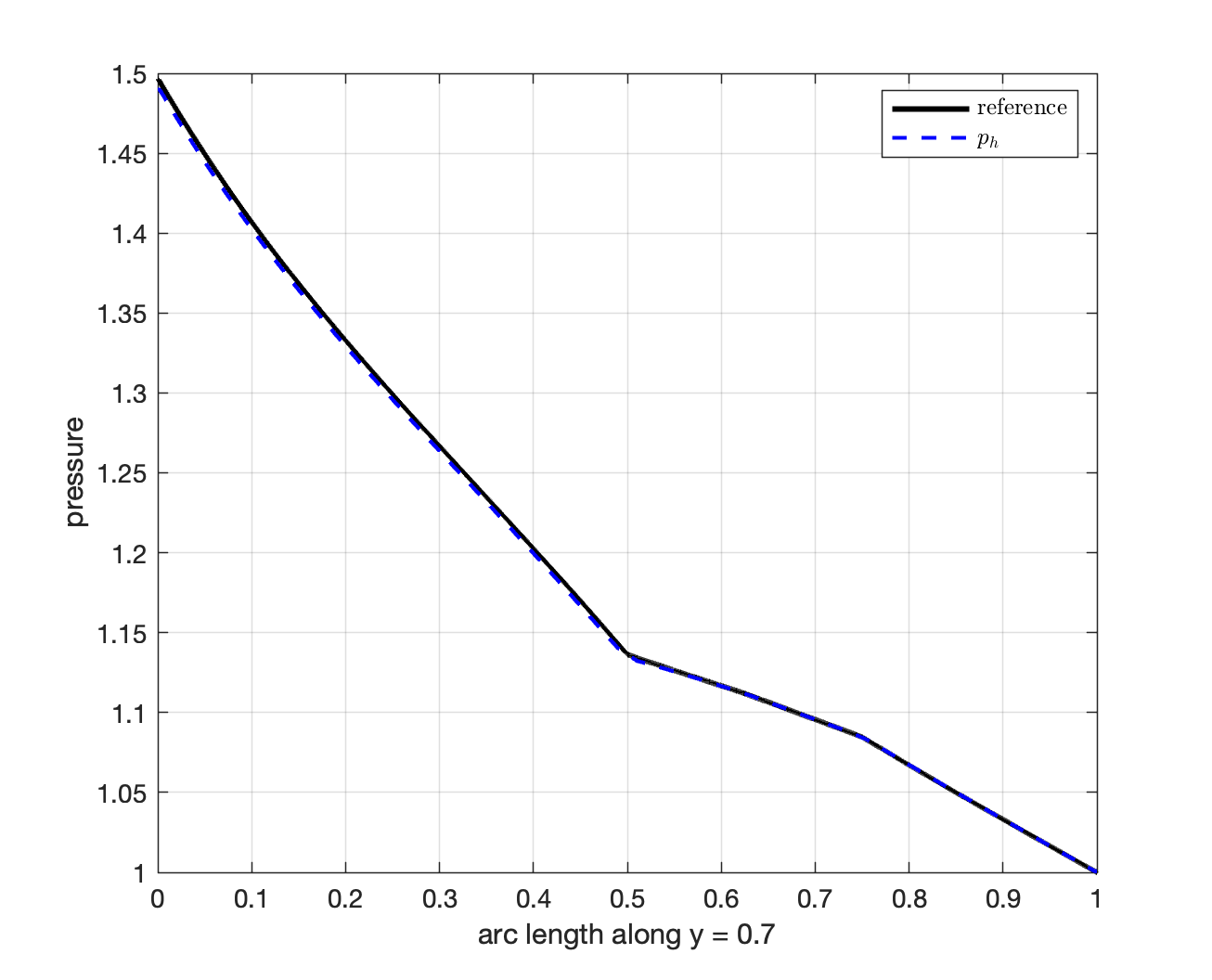}
\caption{fractures, slice $y=0.7$}
\end{subfigure}
\hfill
\begin{subfigure}[b]{0.32\textwidth}
\includegraphics[width=\textwidth]{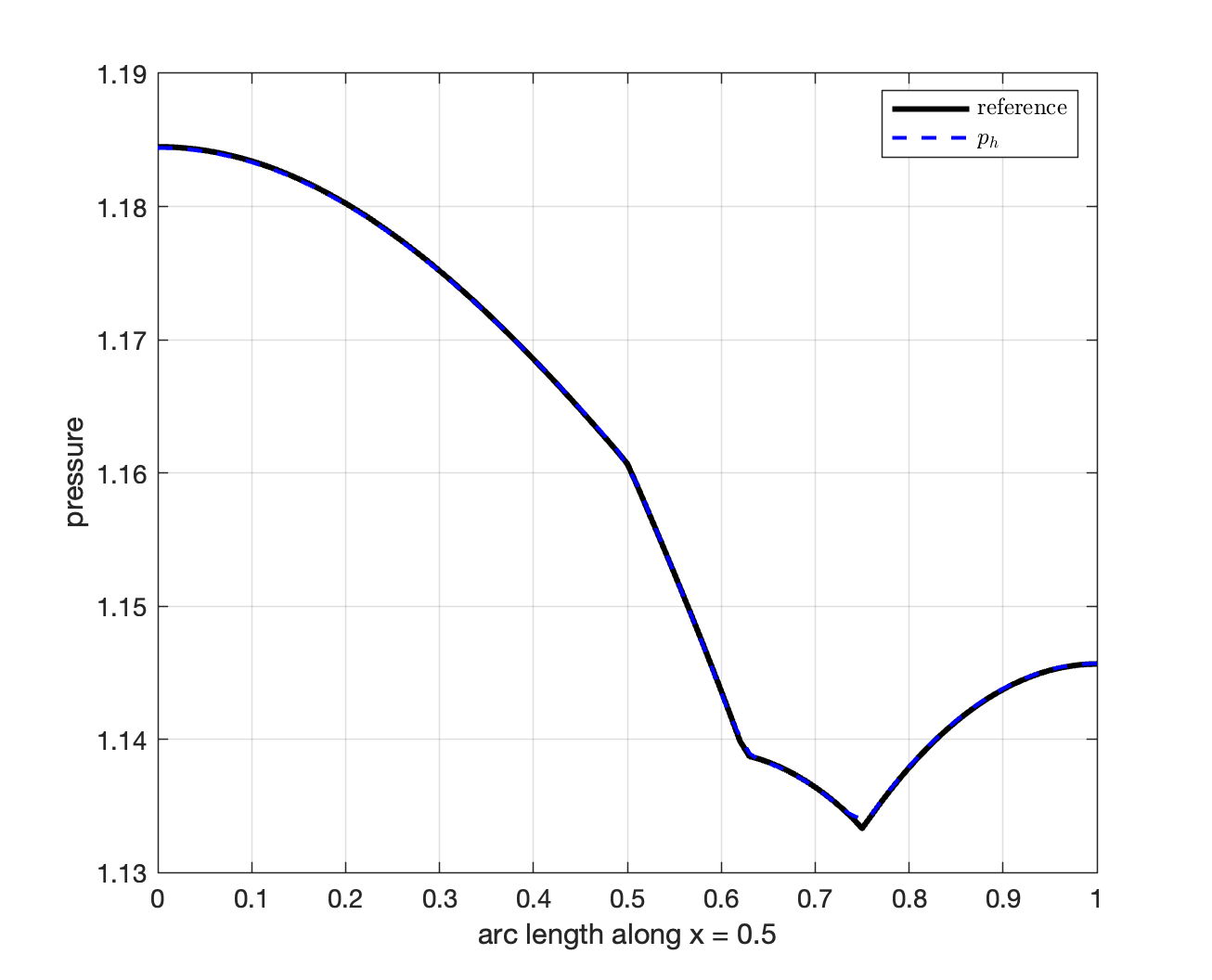}
\caption{fractures, slice $x=0.5$}
\end{subfigure}
\hfill
\begin{subfigure}[b]{0.32\textwidth}
\includegraphics[width=\textwidth]{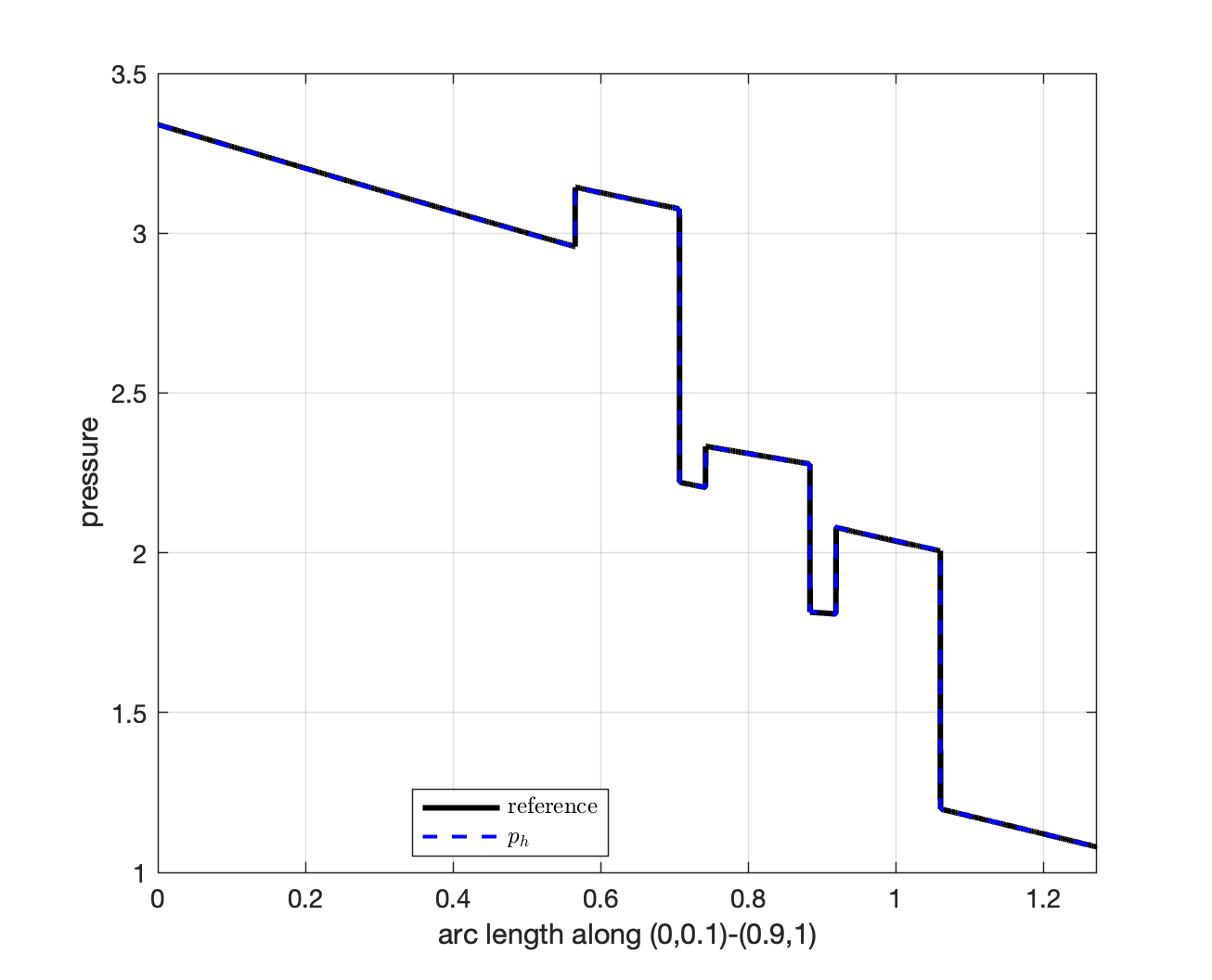}
\caption{barriers, slice $(0,0,1)-(0.9,1)$}
\end{subfigure}
\caption{Example 3, pressure profiles along the slices on the $N=49$ grid, against the mimetic finite difference reference.}
\label{fig:example3_slices}
\end{figure}

Figure~\ref{fig:example3_fields} shows the two pressure fields.
The conductive fractures act as a drainage system while the low-permeability barriers instead block the flow entering from the left.  Figure~\ref{fig:example3_slices} compares the pressure along the three reference slices, which shows excellent agreement.  

\subsection{Example 4: a complex network with mixed interface types}
\label{sec:example4}

The fourth experiment is the mixed-network benchmark of \cite{flemisch2018benchmarks}, in which both interface types appear simultaneously and intersect each other.  
The unit square with $\mathbf K_m=\mathbf I$ contains ten oblique line segments whose coordinates are listed in Appendix~C of \cite{flemisch2018benchmarks}: 
eight of them are conductive fractures with $k_f=10^4$ and two are blocking barriers with $k_b=10^{-4}$.  
All segments share the aperture $a=10^{-4}$, so $a_fk_f=1$ and $a_b/k_b=1$.  
Two flow configurations are computed on the same geometry: a predominantly vertical flow, driven by $p=4$ on the top and $p=1$ on the bottom with impermeable lateral sides, and a predominantly horizontal flow, driven by $p=4$ on the left and $p=1$ on the right with impermeable top and bottom.
The grid is uniform with $N=99$ cells per direction.  
At the fracture-barrier intersection, we let the barrier dominate: in a cell crossed by a barrier, no fracture contribution is assembled.  
This treatment plays a role similar to the harmonic averaging of the two interface permeabilities employed by the methods compared in \cite{flemisch2018benchmarks}.

\begin{figure}[htbp!]
\centering
\begin{subfigure}[b]{0.4\textwidth}
\includegraphics[width=\textwidth]{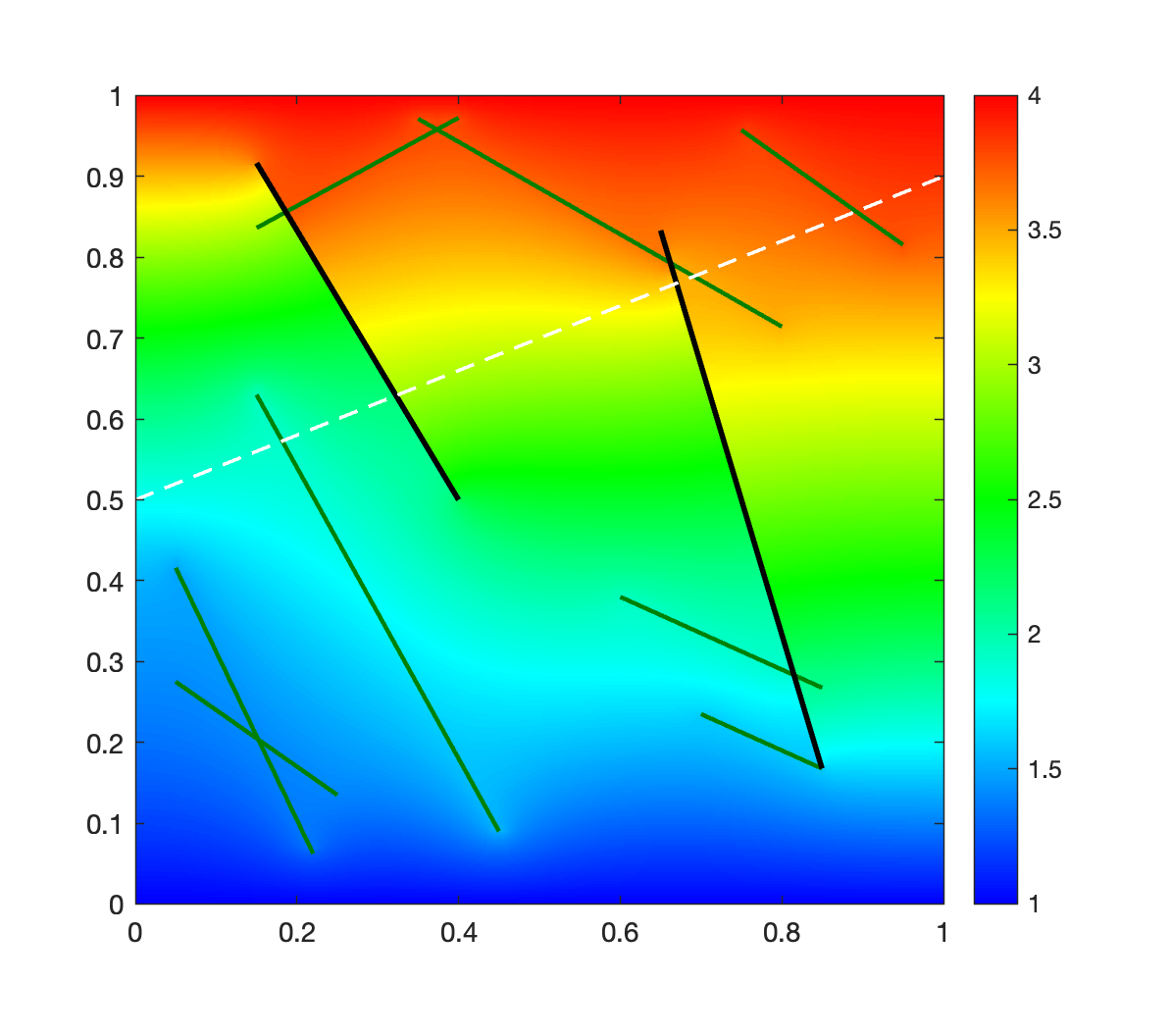}
\caption{vertical flow}
\end{subfigure}
\begin{subfigure}[b]{0.4\textwidth}
\includegraphics[width=\textwidth]{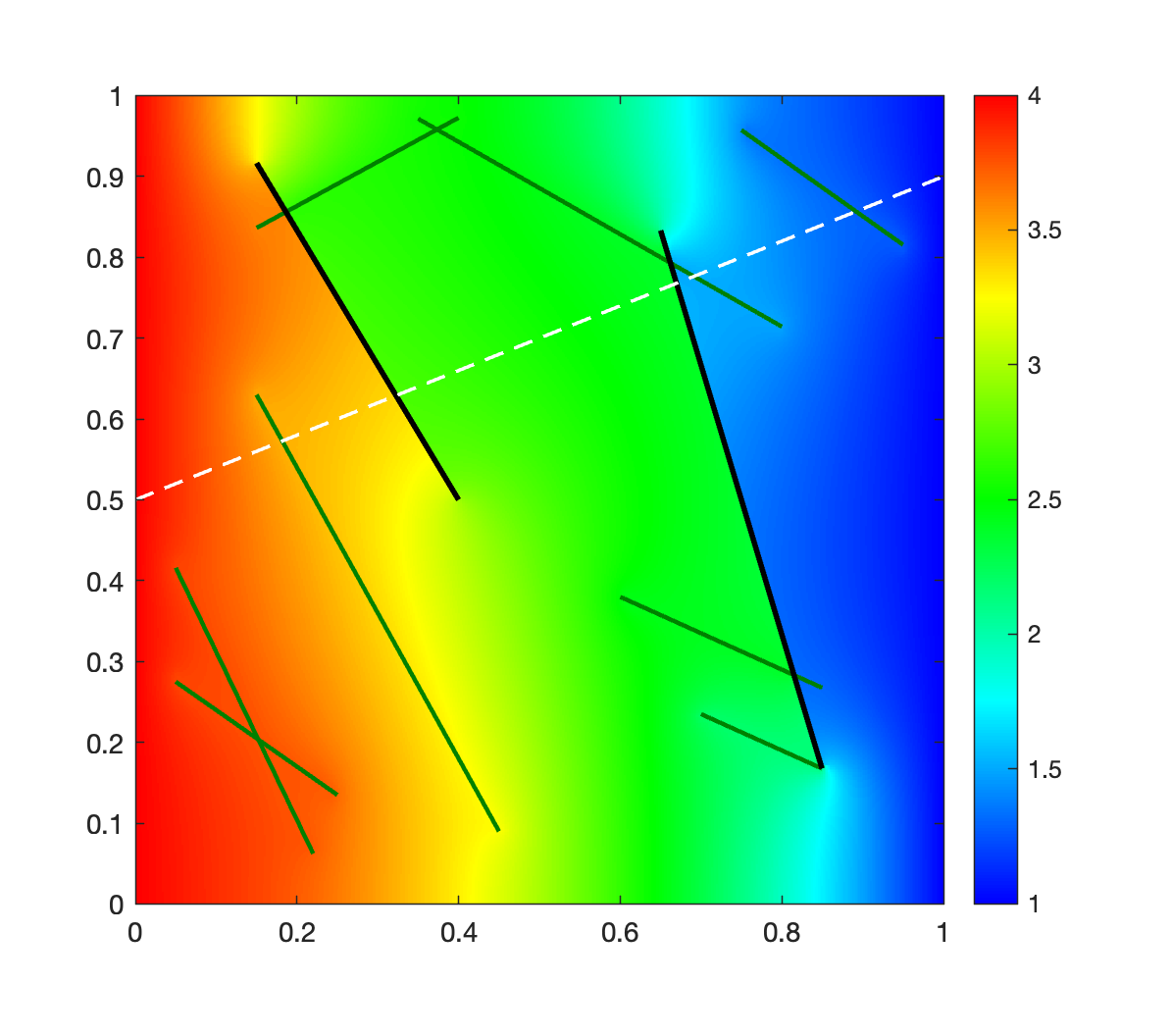}
\caption{horizontal flow}
\end{subfigure}
\caption{Example 4, pressure fields of the two flow configurations on the $N=99$ grid.
The dashed line marks the sampling slice from $(0,0.5)$ to $(1,0.9)$.}
\label{fig:example4_fields}
\end{figure}

\begin{figure}[htbp!]
\centering
\begin{subfigure}[b]{0.4\textwidth}
\includegraphics[width=\textwidth]{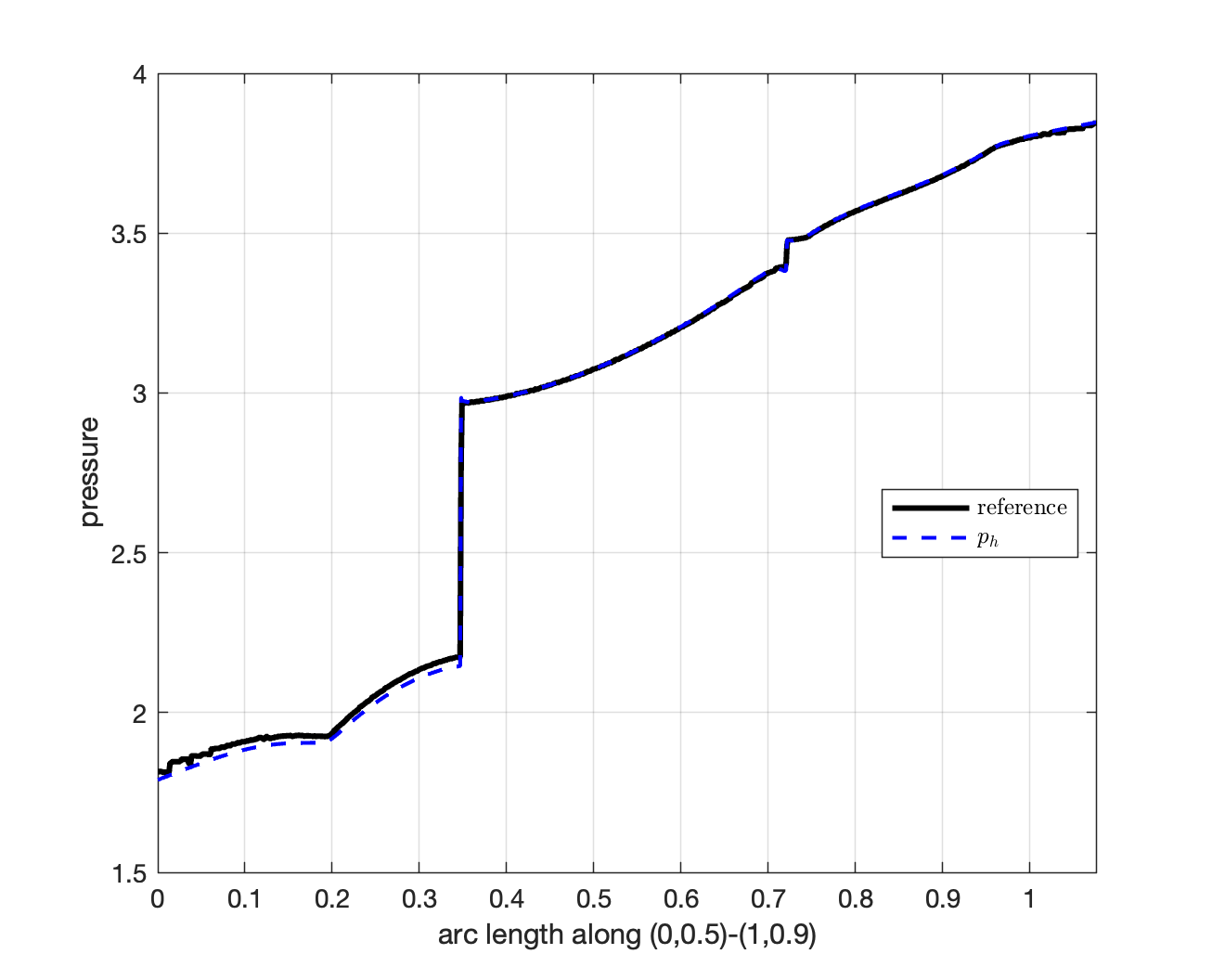}
\caption{vertical flow}
\end{subfigure}
\begin{subfigure}[b]{0.4\textwidth}
\includegraphics[width=\textwidth]{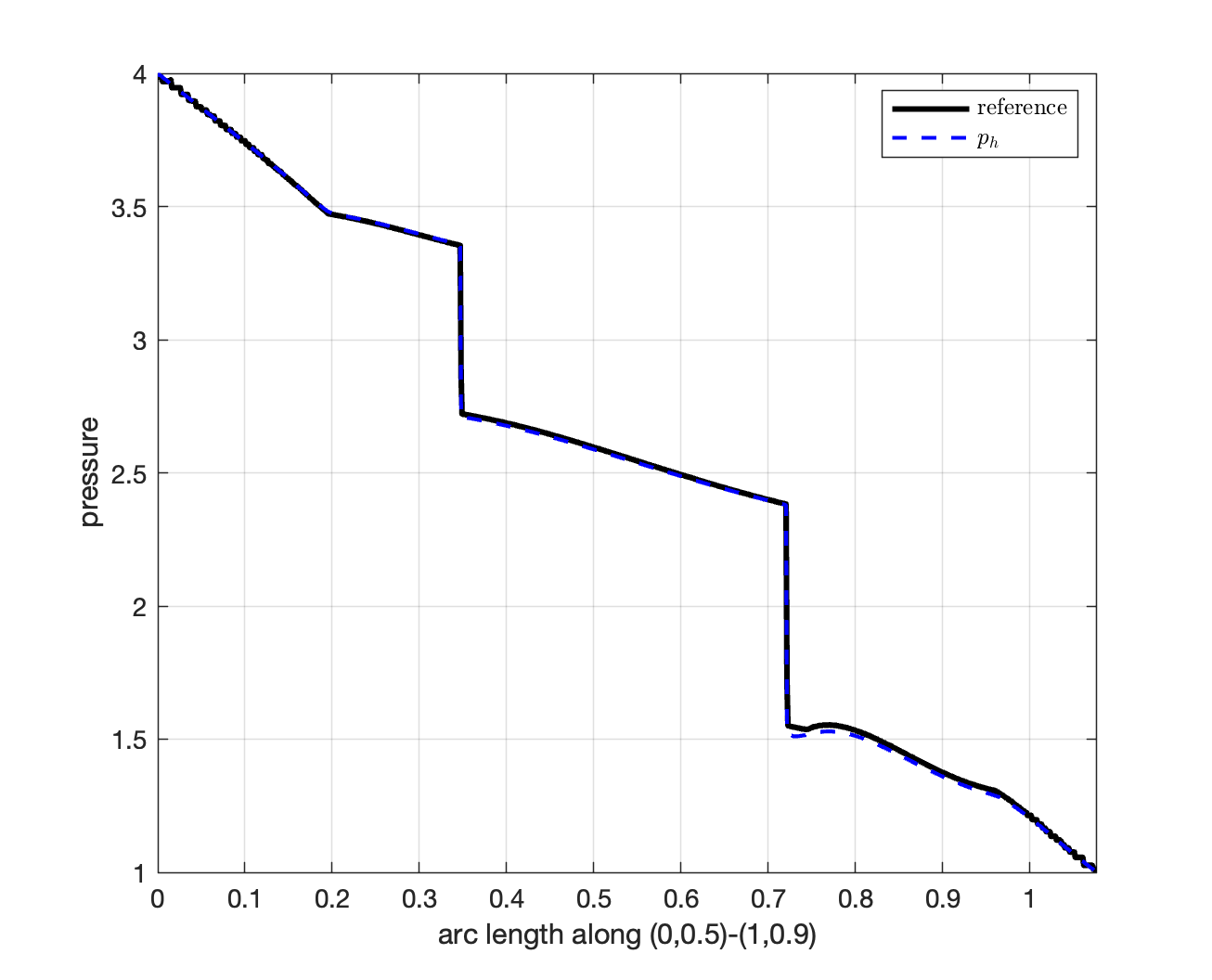}
\caption{horizontal flow}
\end{subfigure}
\caption{Example 4, pressure along the slice from $(0,0.5)$ to $(1,0.9)$ on
the $N=99$ grid, against the mimetic finite difference reference.}
\label{fig:example4_slices}
\end{figure}

Figure~\ref{fig:example4_fields} shows the pressure of both configurations.  
The two interface types shape the field in opposite ways: the fractures short-circuit the pressure differences along their length, while the barriers support visible discontinuities.  
Figure~\ref{fig:example4_slices} compares the pressure profile along the slice from $(0,0.5)$ to $(1,0.9)$, which crosses both barriers, with the reference solution, a mimetic finite difference computation provided by the authors of \cite{flemisch2018benchmarks}.  
In both configurations, the numerical profile reproduces the reference closely, including the sharp pressure jumps at the two barrier crossings and the gentle slope changes where the slice passes conductive fractures.

\subsection{Example 5: a realistic network}
\label{sec:example5}

The final two-dimensional experiment applies the method to a realistic geometry: the network of $63$ line segments digitized from an outcrop, adopted in the benchmark set of \cite{flemisch2018benchmarks}.  
The domain measures $700\times600$ meters and the segments range from a few meters to several hundred meters in length.
The matrix permeability is $k_m=10^{-14}\,\mathrm m^2$ and all features share the
aperture $a=10^{-2}\,\mathrm m$.  Flow is driven by the pressure $1{,}013{,}250$
Pa on the left boundary against zero on the right, with impermeable top and
bottom.  Two conductivities of the same geometry are
considered.  In case~(a), following \cite{flemisch2018benchmarks}, every
feature is a conductive fracture with $k_f=10^{-8}\,\mathrm m^2$, so that
$a_fk_f=10^{-10}\,\mathrm m^3$.  In case~(b), following the barrier
variant of the benchmark introduced in \cite{glaser2022comparison}, every
feature is a blocking barrier with $k_b=10^{-18}\,\mathrm m^2$, giving
$a_b/k_b=10^{16}\,\mathrm m^{-1}$.  The base grid consists of square cells of side $h=5$
meters, that is, $140\times120$ cells; the
slice comparisons additionally include a refined grid with $h=2.5$
($280\times240$ cells).
The two cases are assessed differently.  
For the fracture case the benchmark study
\cite{flemisch2018benchmarks} did not provide a reference solution for this geometry. 
We therefore compare our profiles along $x=625$ and
$y=500$ directly with the published solutions of the five methods that participated in that case taken from the benchmark's public repository.  
For the barrier case, a fine box-DFM reference \cite{xu2025extension} on about $3\times10^{5}$ conforming cells is available, sampled along the diagonal from $(0,0)$ to $(700,600)$ and along $x=625$.

\begin{figure}[htbp!]
\centering
\begin{subfigure}[b]{0.4\textwidth}
\includegraphics[width=\textwidth]{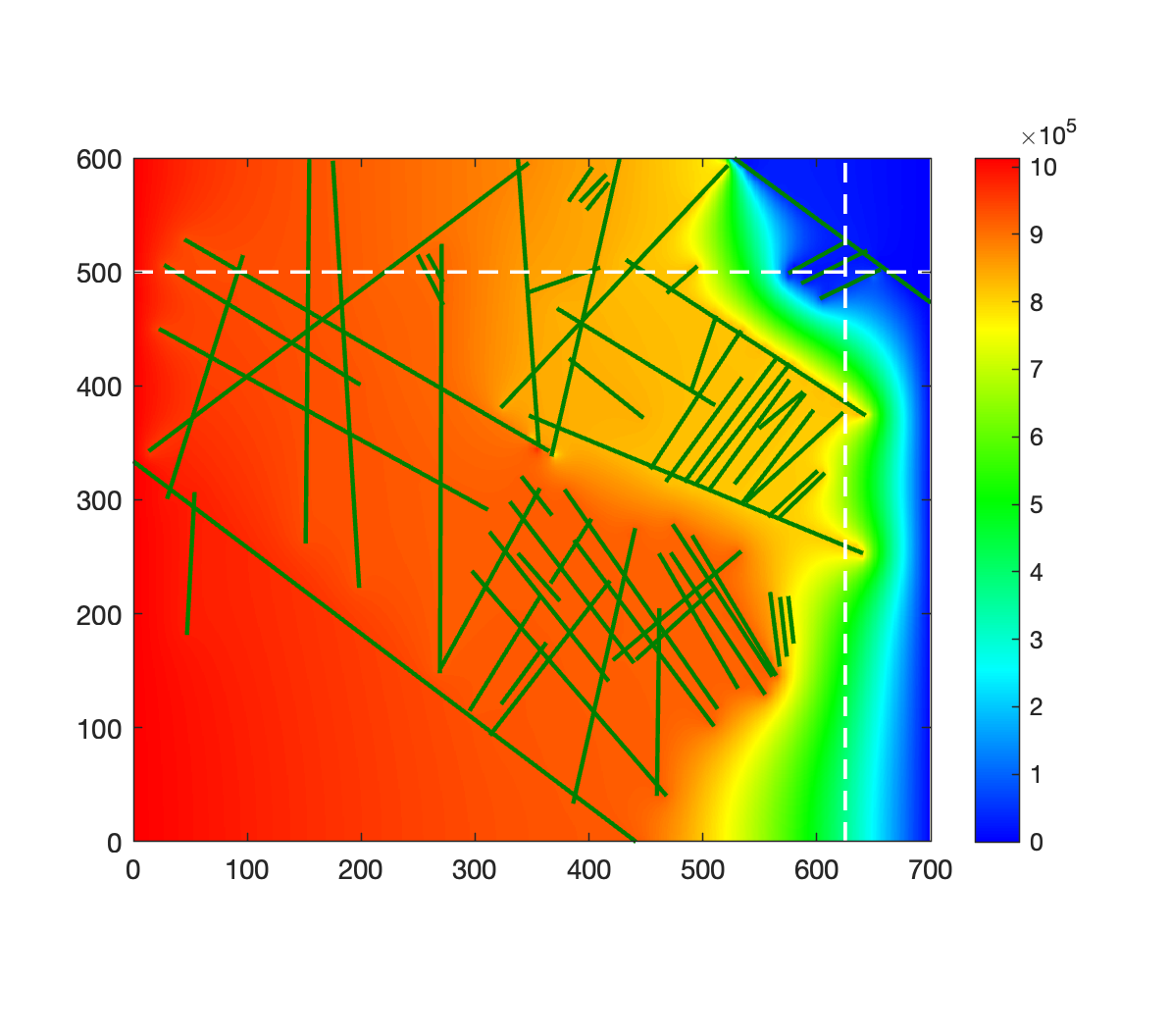}
\caption{conductive fractures}
\end{subfigure}
\begin{subfigure}[b]{0.4\textwidth}
\includegraphics[width=\textwidth]{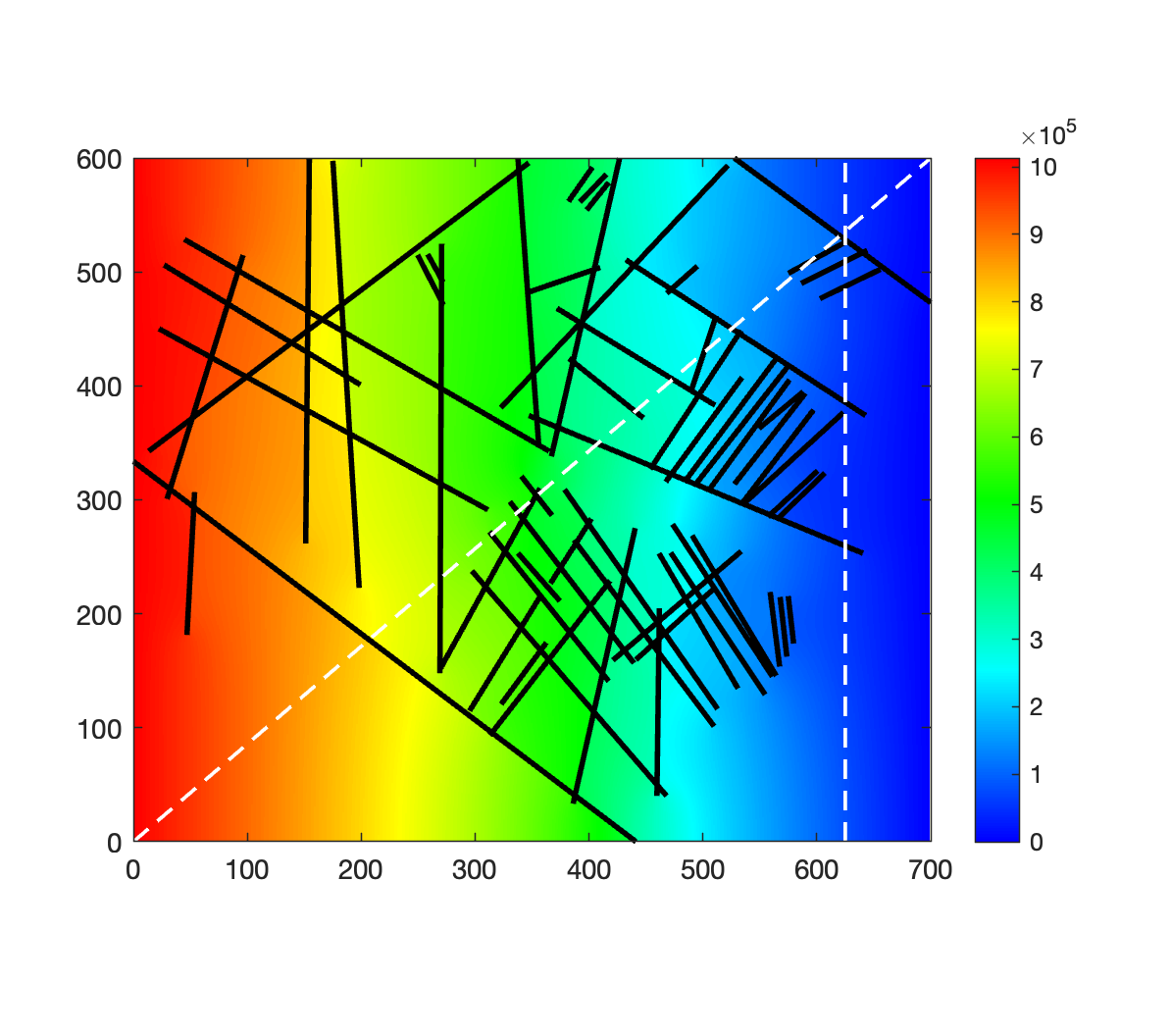}
\caption{blocking barriers}
\end{subfigure}
\caption{Example 5, pressure fields on the $140\times120$ grid ($h=5$).  
Dashed lines mark the sampling slices.}
\label{fig:example5_fields}
\end{figure}

\begin{figure}[htbp!]
\centering
\begin{subfigure}[b]{0.4\textwidth}
\includegraphics[width=\textwidth]{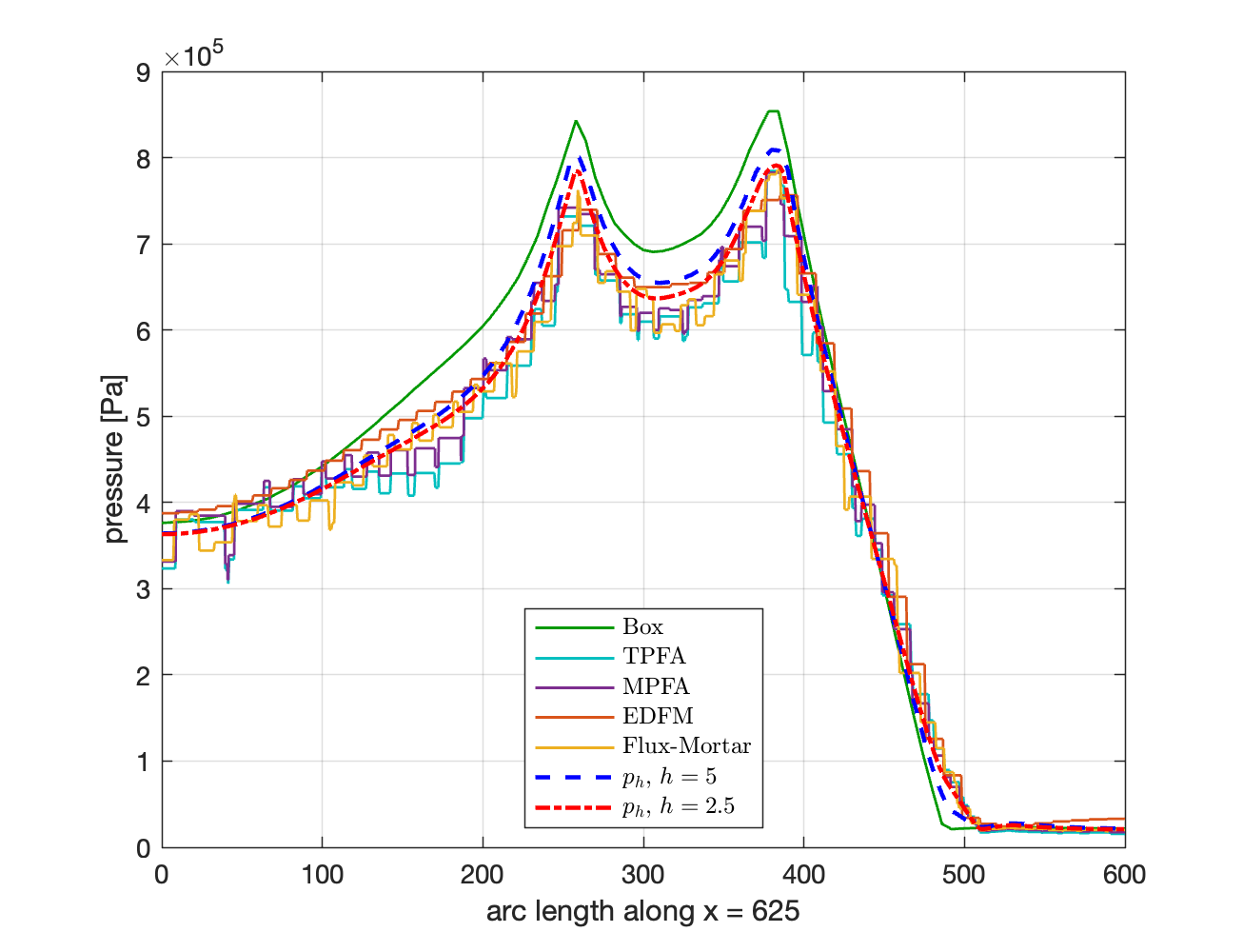}
\caption{fractures, $x=625$}
\end{subfigure}
\begin{subfigure}[b]{0.4\textwidth}
\includegraphics[width=\textwidth]{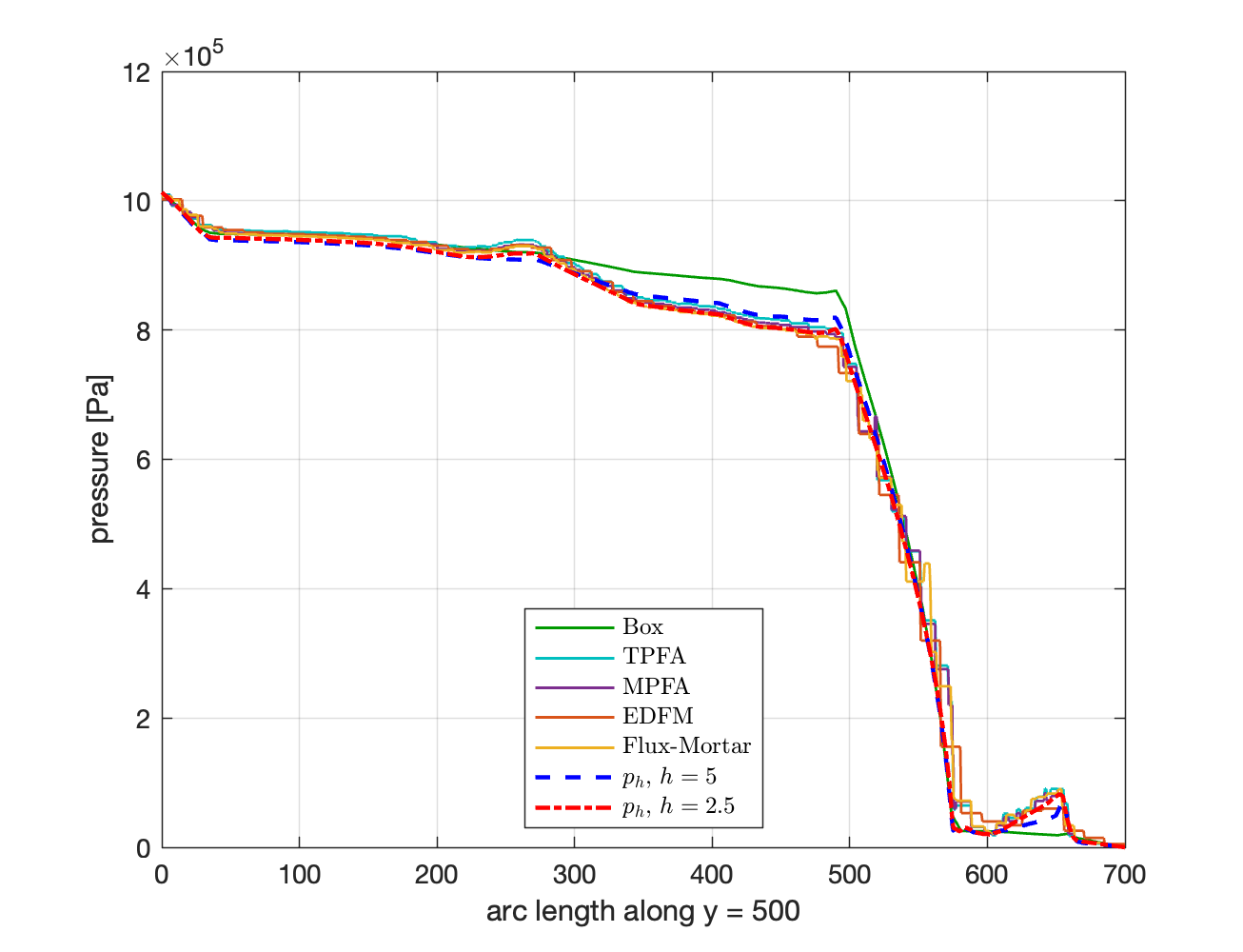}
\caption{fractures, $y=500$}
\end{subfigure}\\
\begin{subfigure}[b]{0.4\textwidth}
\includegraphics[width=\textwidth]{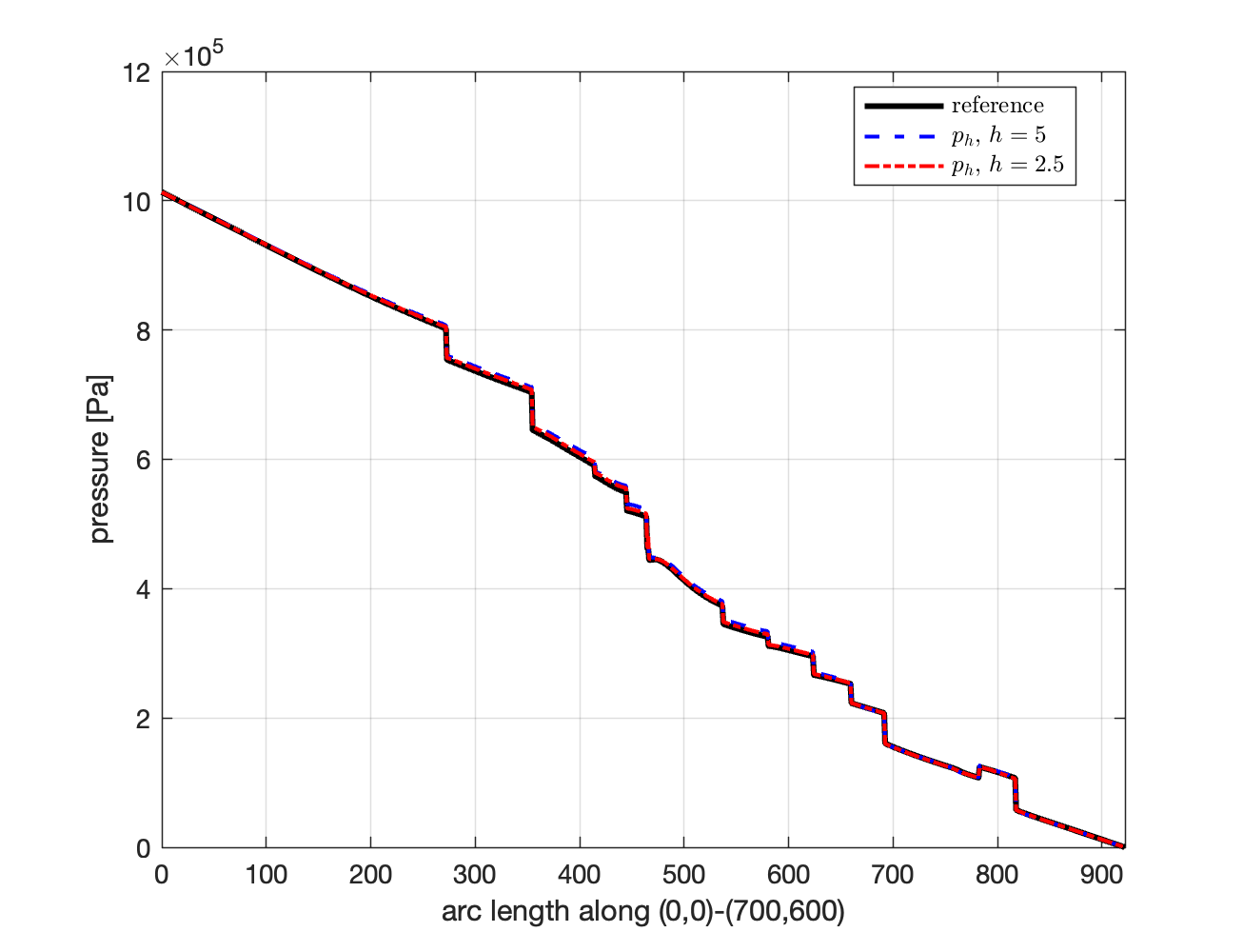}
\caption{barriers, diagonal}
\end{subfigure}
\begin{subfigure}[b]{0.4\textwidth}
\includegraphics[width=\textwidth]{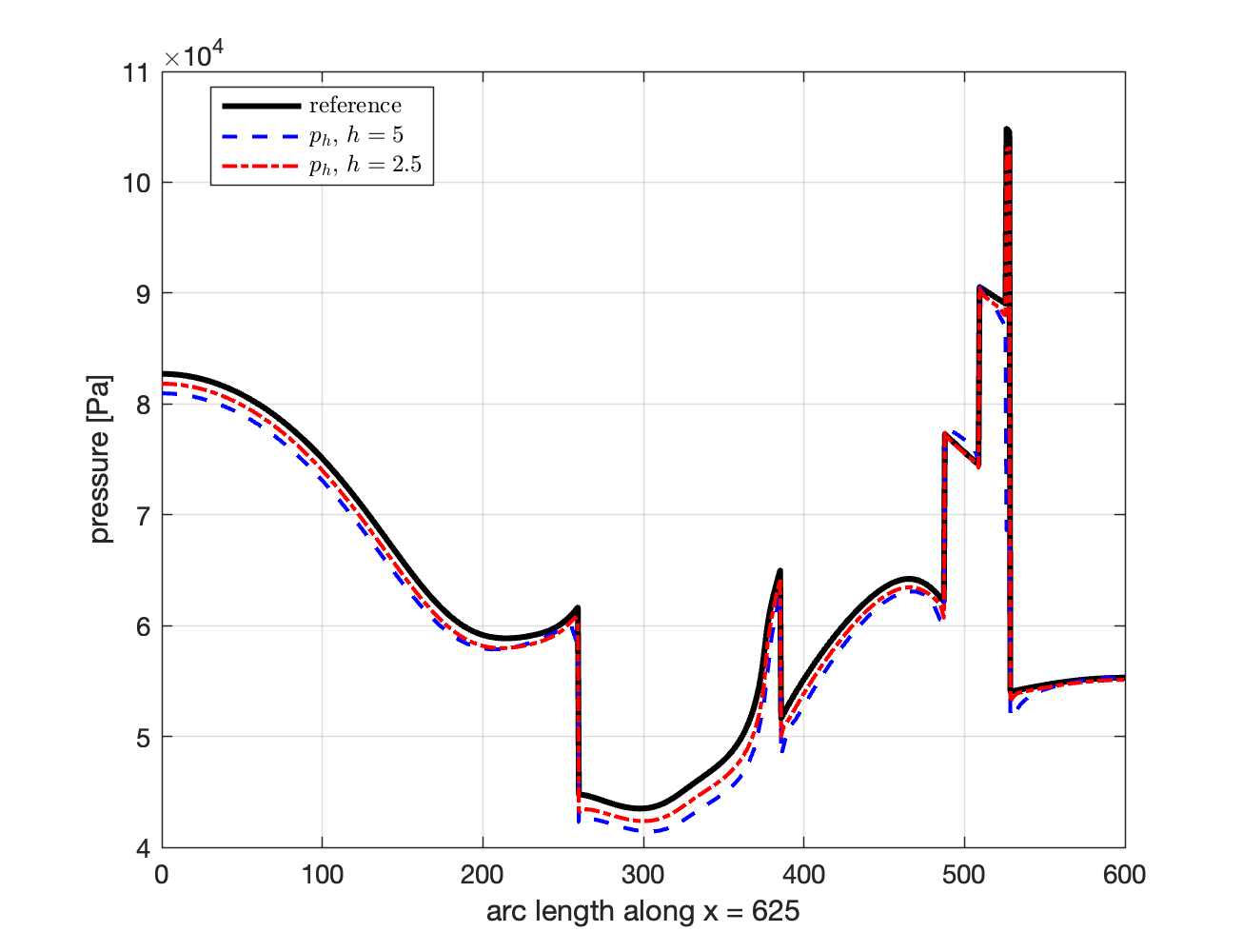}
\caption{barriers, $x=625$}
\end{subfigure}
\caption{Example 5, pressure profiles along the sampling slices on the $h=5$ and $h=2.5$ grids.  Top row (fractures): direct comparison with the
published profiles of the five methods that participated in this case of \cite{flemisch2018benchmarks}.  
Bottom row (barriers): comparison with the
fine conforming box-DFM reference of \cite{xu2025extension}.}
\label{fig:example5_slices}
\end{figure}

Figure~\ref{fig:example5_fields} shows the two pressure fields, and
Figure~\ref{fig:example5_slices} the four slice comparisons.  
In the fracture case, along both slices, our profiles lie within the family of the five published solutions over their entire length.
In the barrier case, the field remains globally graded from left to right along the diagonal slice.
The profile along $x=625$, which varies within a band of less than one tenth of the total pressure drop, resolves its sequence of six barrier crossings with moderate deviations.  
The only feature the coarse grid fails to resolve is the narrow spike near the top of that profile, where the slice crosses two barriers separated by only about two meters, and the flat crossing there is recovered on the fine grid.

\subsection{Example 6: a single fracture in three dimensions}
\label{sec:example6}

\begin{figure}[htbp!]
\centering
\includegraphics[width=0.5\textwidth]{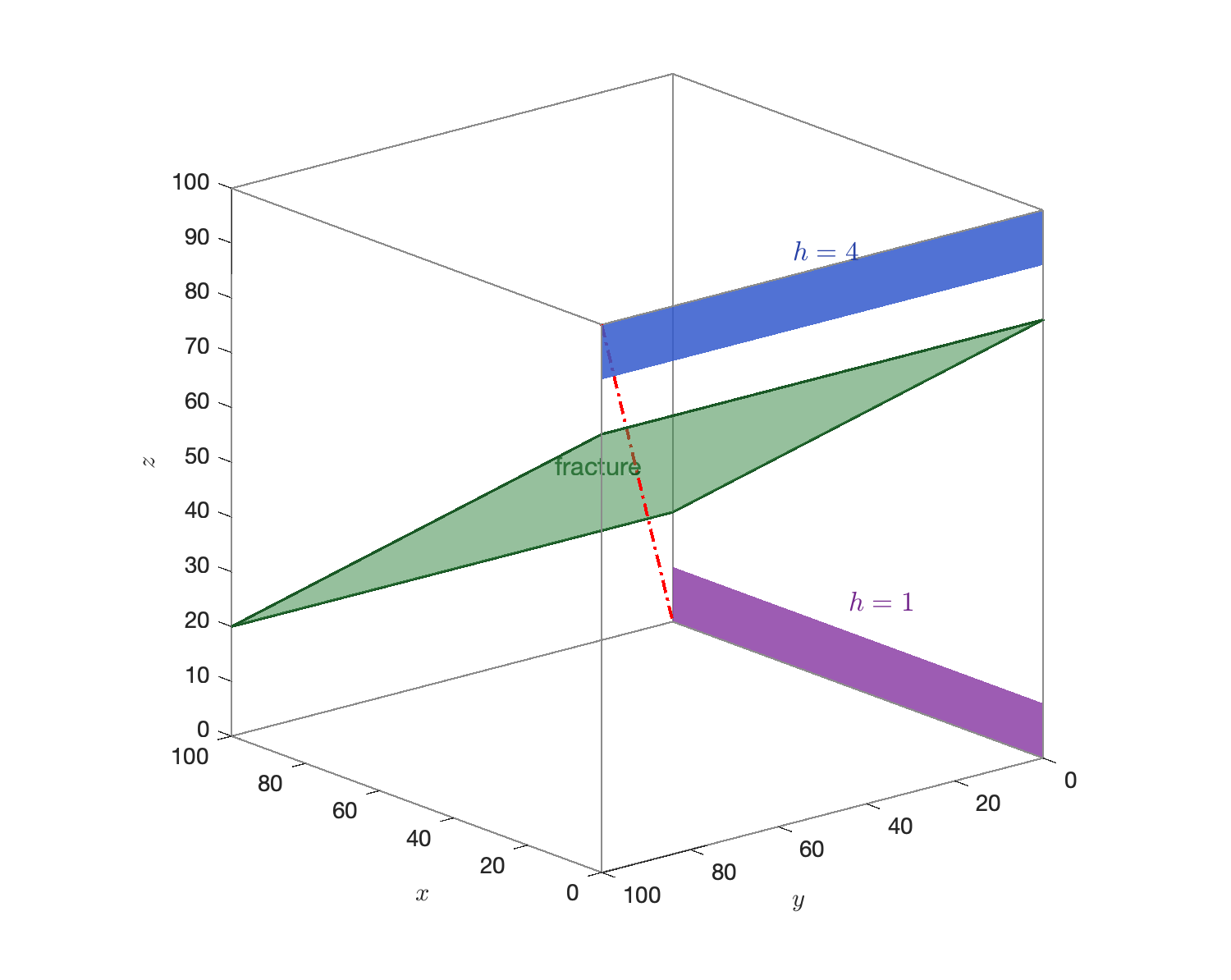}
\caption{Example 6, the computational domain. The inclined plane in green is the conductive fracture, and the colored strips are the inlet band ($h=4$, on the face $x=0$) and the outlet band ($h=1$, on the face $y=0$). The dash-dotted line is the sampling line of the head comparison.
}
\label{fig:example6_geometry}
\end{figure}

The remaining experiments move to three space dimensions and are taken from the verification benchmark set of \cite{berre2021verification}. 
The computations use the three-dimensional schemes summarized in Appendix~\ref{app:three_dimensional_extension}.  
The first three-dimensional test is the single, conductive fracture case.
The configuration is sketched in Figure~\ref{fig:example6_geometry}.
The domain is the cube $(0,100)^3$, crossed by the inclined planar fracture with corners $(0,0,80)$, $(100,0,20)$, $(100,100,20)$ and $(0,100,80)$ of aperture $a_f=10^{-2}$.  
The matrix conductivity is $10^{-6}\,\mathrm{m/s}$ for $z>10$ and raised to $10^{-5}\,\mathrm{m/s}$ in the layer $z<10$.
The fracture carries the effective tangential conductivity $a_fk_f=10^{-3}\,\mathrm{m^2/s}$. 
The flow is driven by the hydraulic head $h=4$ prescribed on the narrow band $\{0\}\times(0,100)\times(90,100)$ of one vertical face and $h=1$ on the band $(0,100)\times\{0\}\times(0,10)$ of an adjacent face; the remainder of the boundary is impermeable.  
We solve on uniform grids with $N=15$, $25$ and $45$ cells per direction, near the three refinement levels of the benchmark. 

\begin{figure}[htbp!]
\centering
\begin{subfigure}[b]{0.32\textwidth}
\includegraphics[width=\textwidth]{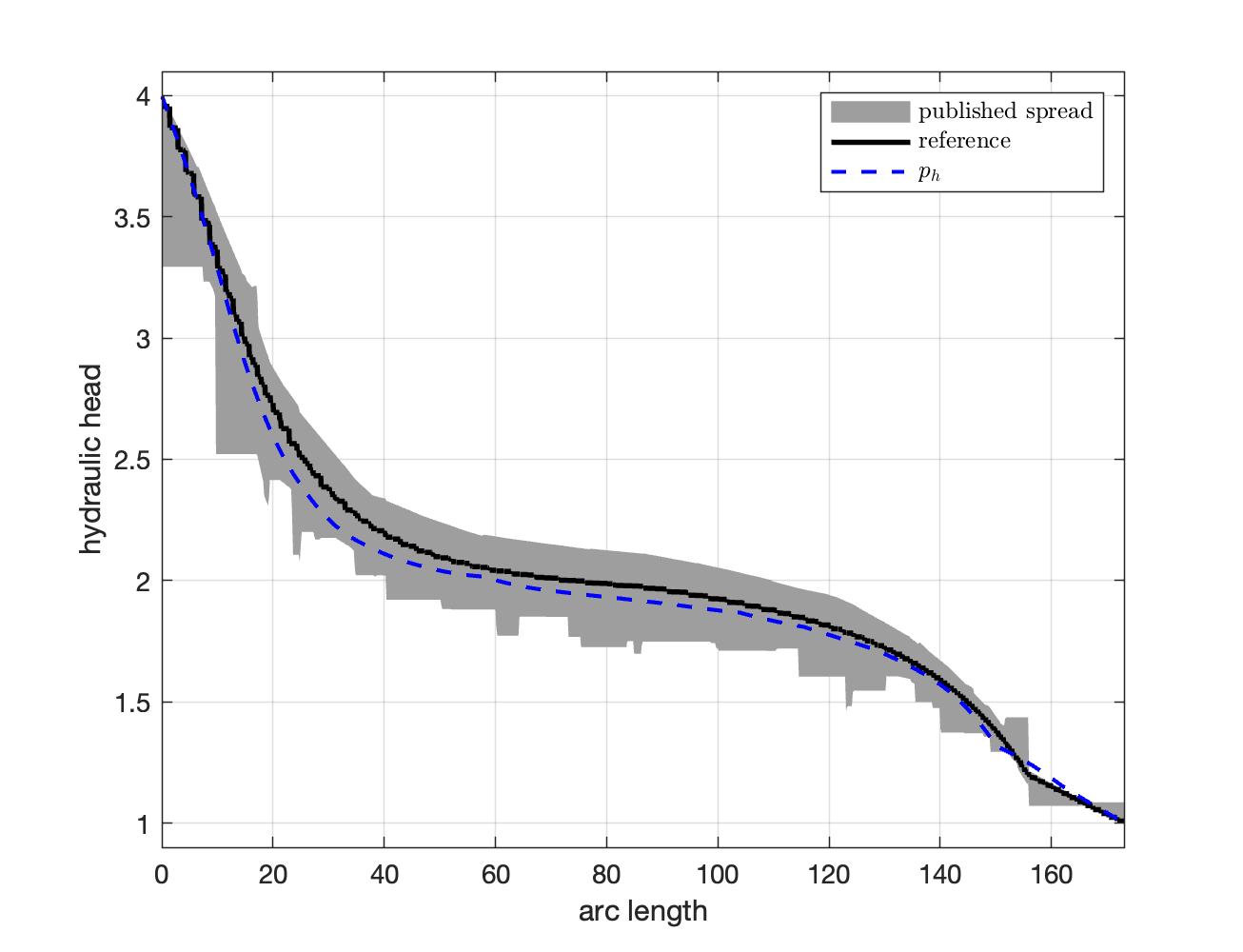}
\caption{$N=15$ ($3{,}375$ cells)}
\end{subfigure}
\hfill
\begin{subfigure}[b]{0.32\textwidth}
\includegraphics[width=\textwidth]{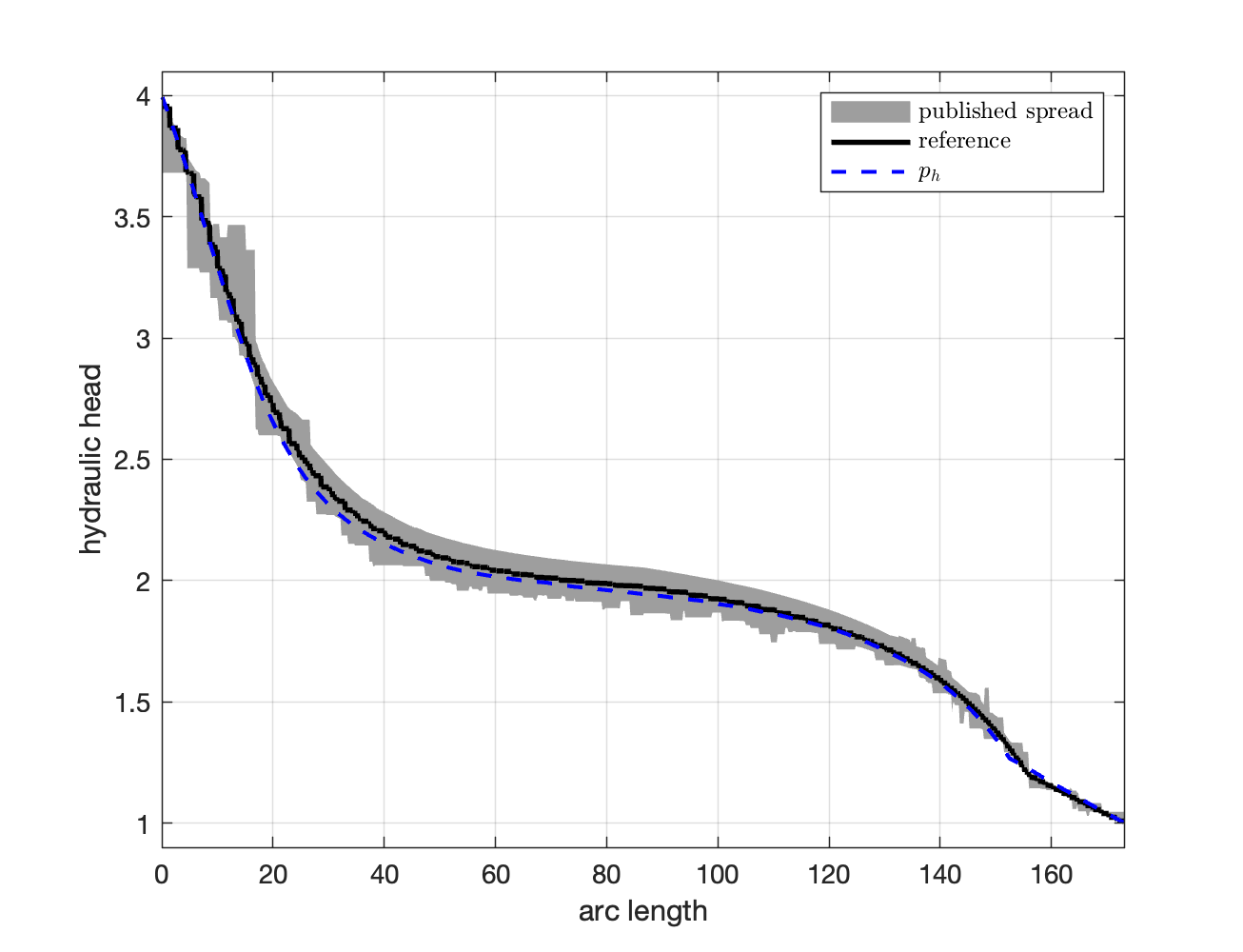}
\caption{$N=25$ ($15{,}625$ cells)}
\end{subfigure}
\hfill
\begin{subfigure}[b]{0.32\textwidth}
\includegraphics[width=\textwidth]{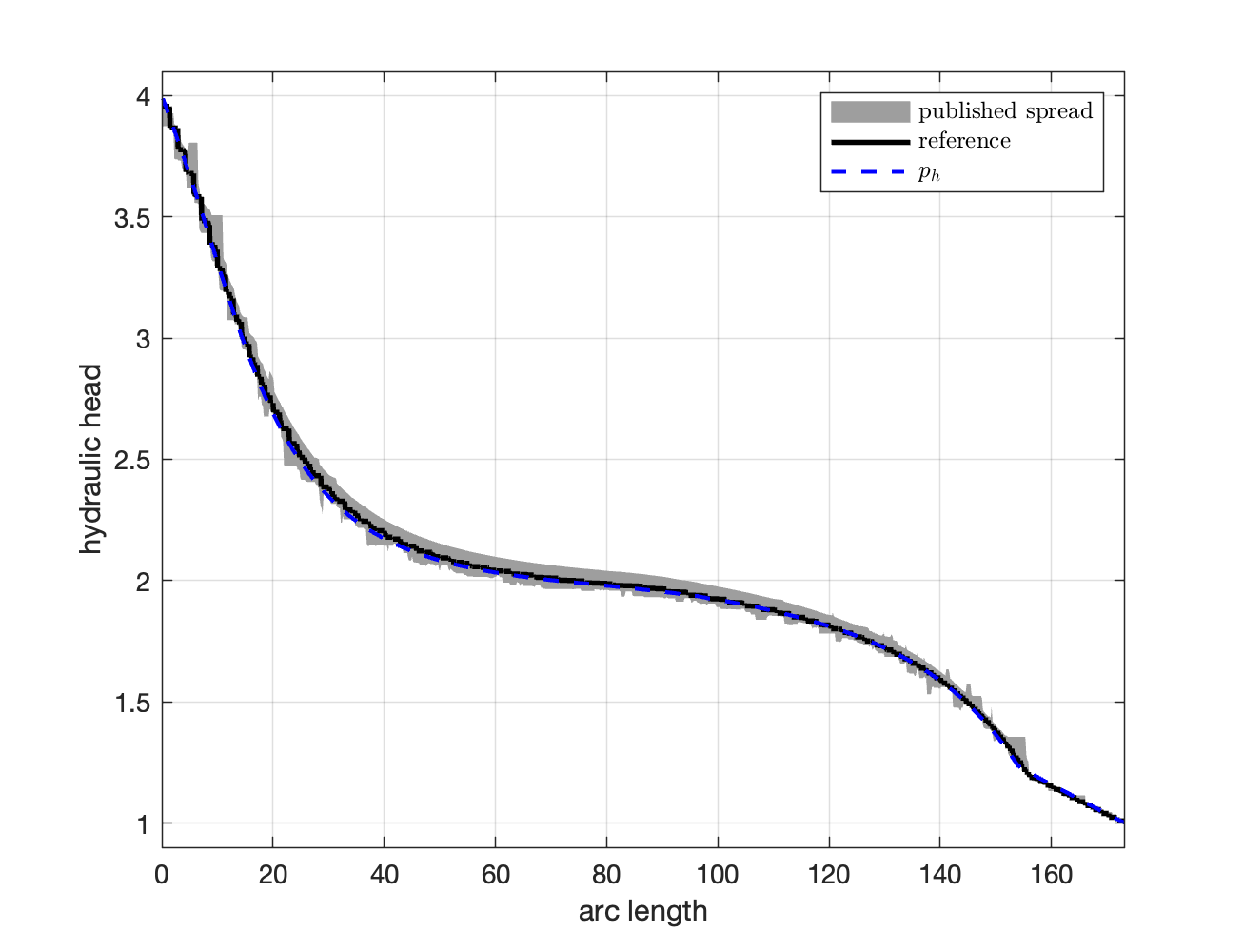}
\caption{$N=45$ ($91{,}125$ cells)}
\end{subfigure}
\caption{Example 6, hydraulic head along the line $(0,100,100)$--$(100,0,0)$ on the three levels of refinement.
The shaded region spans the $10$th-$90$th percentiles of the seventeen methods collected in \cite{berre2021verification} at refinement levels of $\sim 1$k, $\sim10$k, and $\sim100k$ cells. 
The black curve is the reference solution (USTUTT-MPFA).}
\label{fig:example6_pol}
\end{figure}

Figure~\ref{fig:example6_pol} compares the head along the line from $(0,100,100)$ to $(100,0,0)$ with the data collected in \cite{berre2021verification}. 
The solid black curve is the reference
solution, computed with the USTUTT-MPFA scheme on a refined grid, and the shaded region spans the 10th to 90th
percentiles of the seventeen participating methods at the matching
refinement level.  Already on the coarsest grid the computed profile runs
inside the published spread along the reference curve, and it approaches
the reference under refinement.

\subsection{Example 7: a regular fracture network in three dimensions}
\label{sec:example7}

\begin{figure}[htbp!]
\centering
\includegraphics[width=0.5\textwidth]{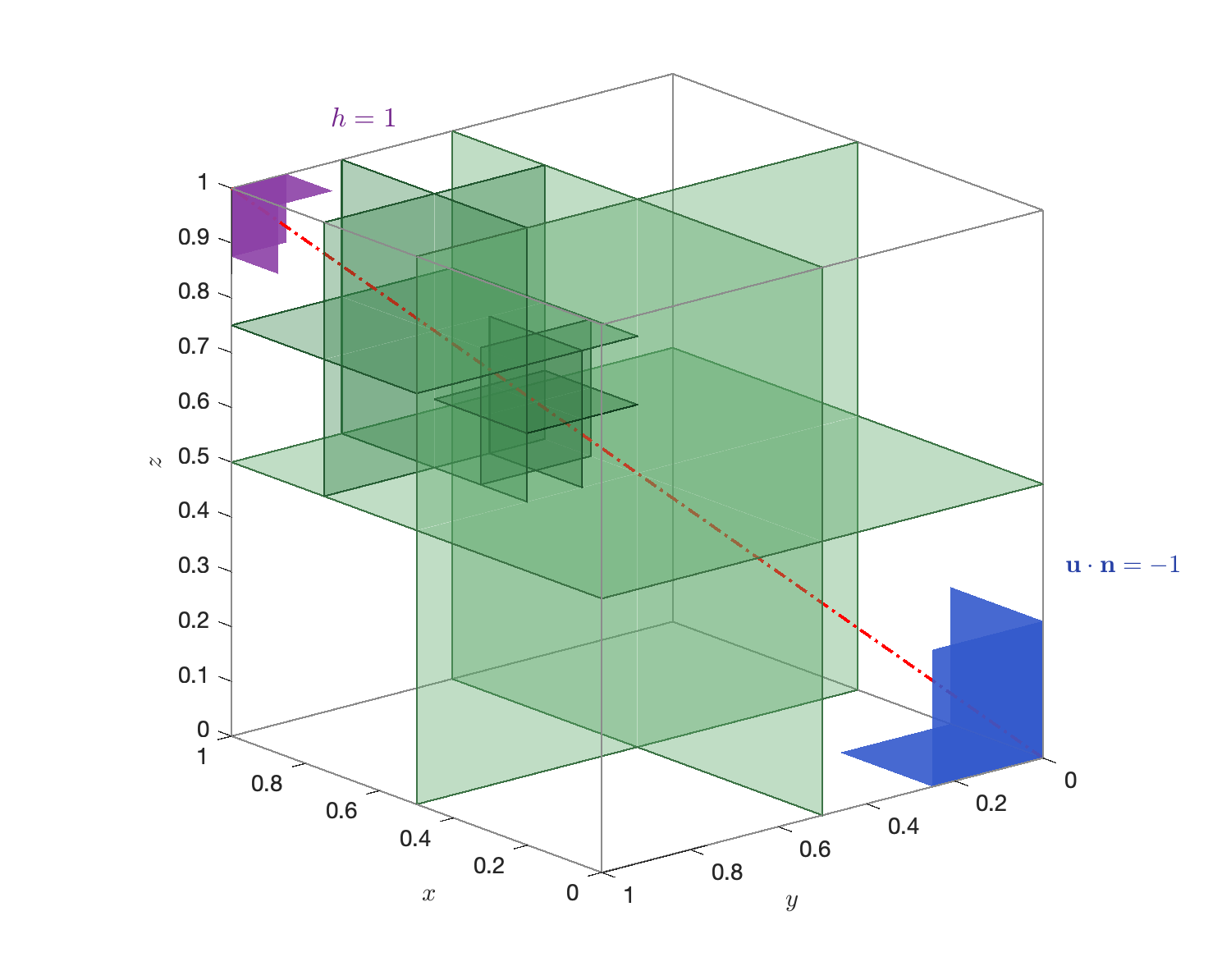}
\caption{Example 7, the computational domain. The nine axis-parallel fracture planes are shown in green (darker shades for the smaller planes), the inflow patches around the origin in blue, and the outflow patches around the corner $(1,1,1)$ in purple.  
The dash-dotted line is the sampling diagonal.}
\label{fig:example7_geometry}
\end{figure}

This experiment is the regular-network case of \cite{berre2021verification}, the three-dimensional analog of Example~3.  
The unit cube contains nine axis-parallel planar fractures of aperture $a=10^{-4}$: three planes spanning the whole cube through $0.5$, three quarter planes through $0.75$, and three eighth planes through $0.625$.
The configuration is sketched in Figure~\ref{fig:example7_geometry}.
The matrix conductivity is $0.1$ in the region \begin{equation*}
\begin{split}
\Omega_{1}=&\{(x,y,z)\in\Omega: x>0.5, y<0.5\}\cup
\{(x,y,z)\in\Omega: x>0.75, 0.5<y<0.75, z>0.5\}\\
&\cup
\{(x,y,z)\in\Omega: 0.625<x<0.75, 0.5<y<0.625, 0.5<z<0.75\},
\end{split}
\end{equation*}
and is $1$ in the rest of the bulk region.
In case (1) the fractures are conductive with $a_fk_f=1$ and negligible normal resistance. 
In case (2) they are blocking, with normal resistance $a_b/k_b=1$ and a negligible tangential effect.  
Flow is driven by a uniform influx $\mathbf u\cdot\mathbf n=-1$ across the three boundary patches with $x,y,z<0.25$ around the origin, against the head $h=1$ on the three patches with $x,y,z>0.875$ around the opposite corner.
The rest of the boundary is impermeable.  
We solve on grids with $N=11$, $19$ and $35$ cells per direction

\begin{figure}[htbp!]
\centering
\begin{subfigure}[b]{0.32\textwidth}
\includegraphics[width=\textwidth]{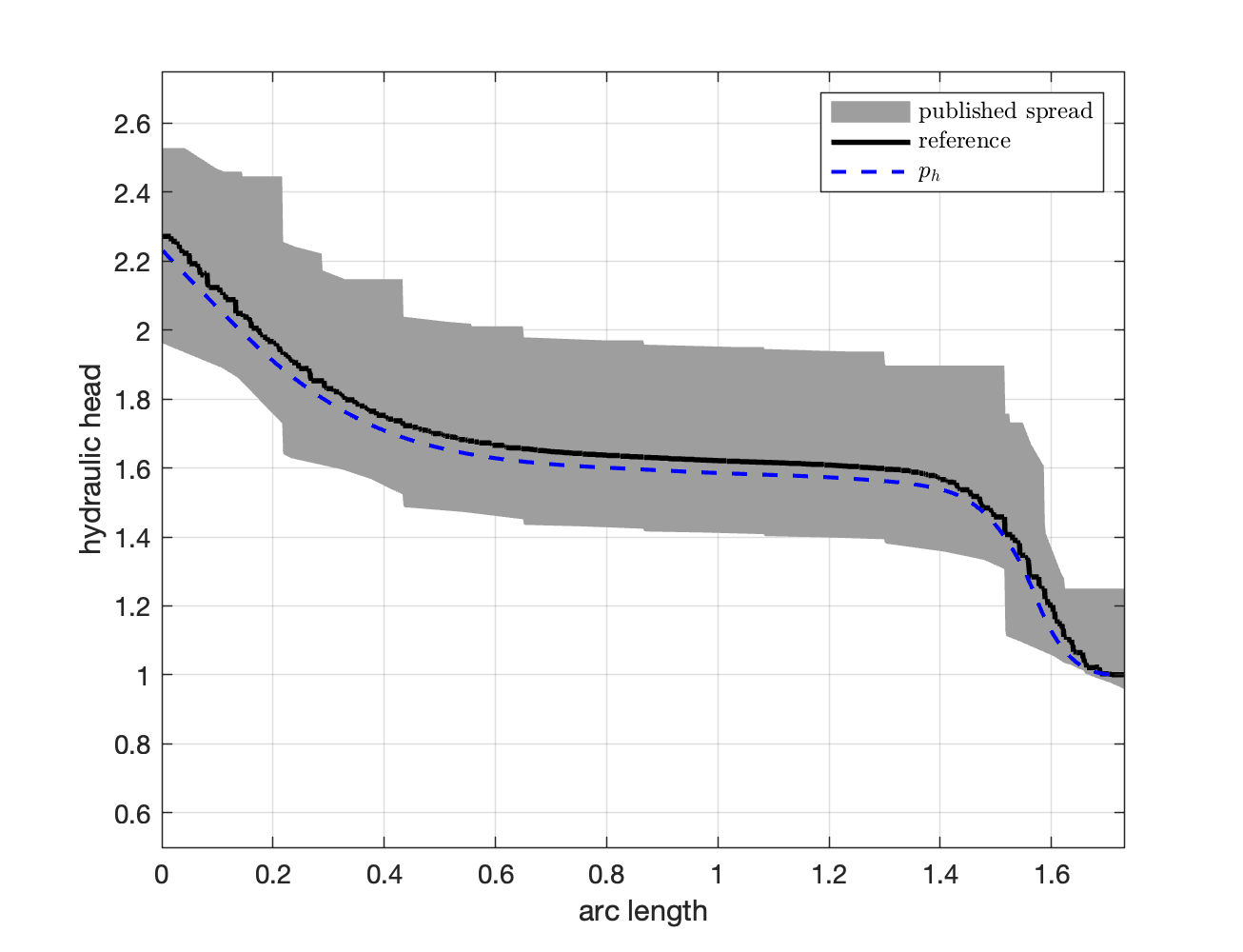}
\caption{$N=11$ ($1{,}331$ cells)}
\end{subfigure}
\hfill
\begin{subfigure}[b]{0.32\textwidth}
\includegraphics[width=\textwidth]{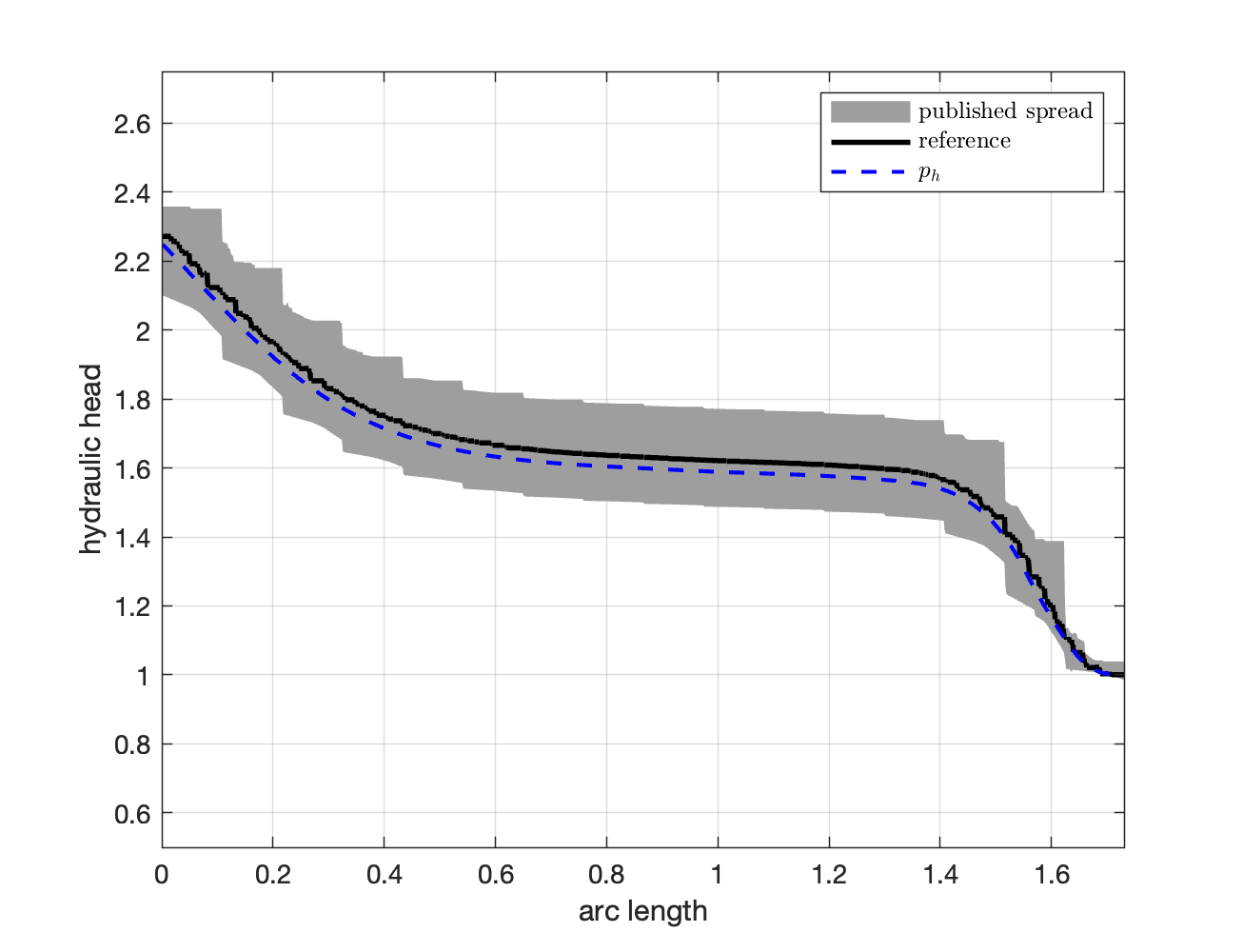}
\caption{$N=19$ ($6{,}859$ cells)}
\end{subfigure}
\hfill
\begin{subfigure}[b]{0.32\textwidth}
\includegraphics[width=\textwidth]{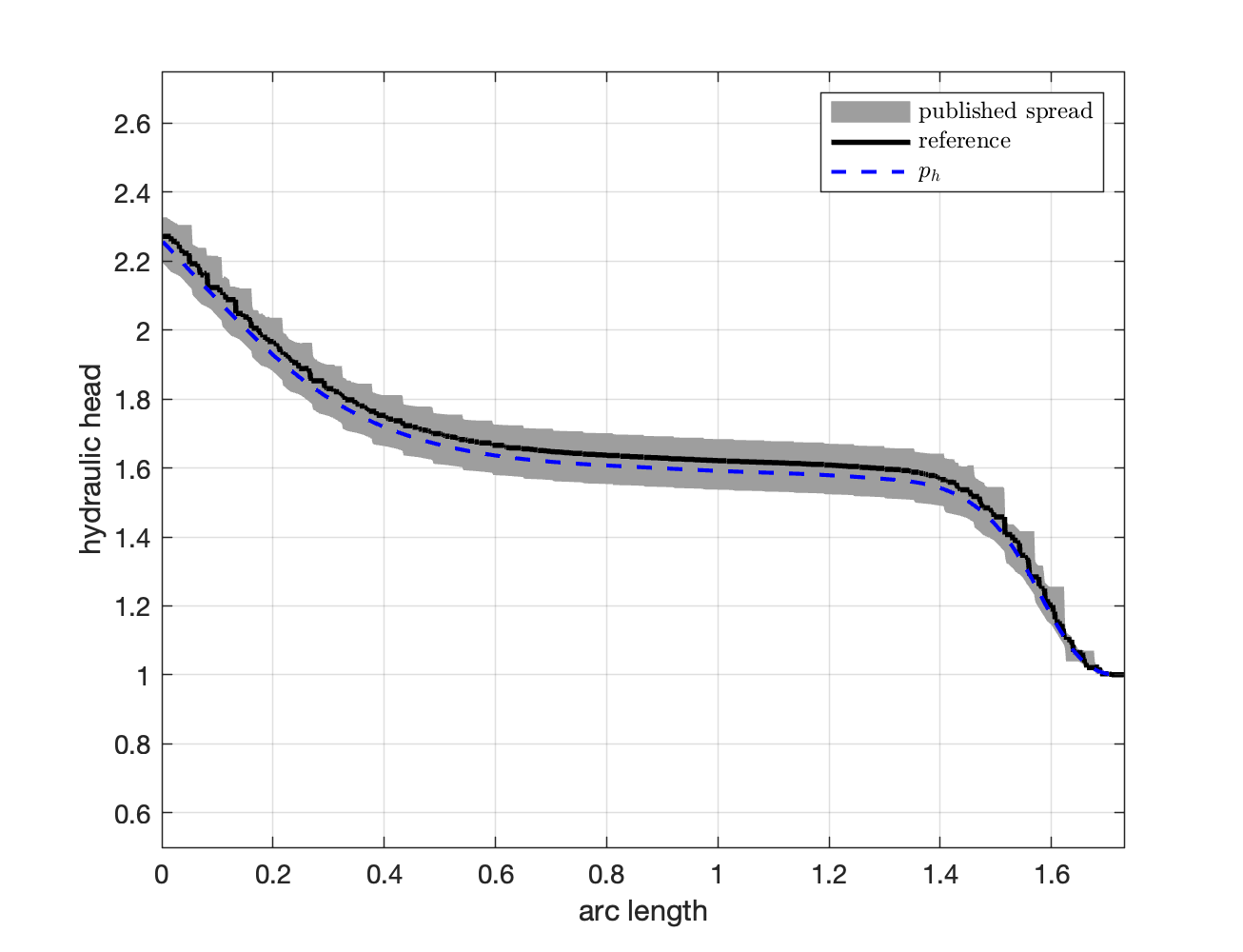}
\caption{$N=35$ ($42{,}875$ cells)}
\end{subfigure}
\caption{Example 7, conductive case, hydraulic head along the diagonal
$(0,0,0)$--$(1,1,1)$ on the three levels of refinement.
The shaded region spans the $10$th-$90$th percentiles of the fifteen methods collected in \cite{berre2021verification} at refinement levels of $\sim 500$, $\sim4$k, and $\sim32k$ cells. 
The black curve is the reference solution
(USTUTT-MPFA).}
\label{fig:example7_cond}
\end{figure}

\begin{figure}[htbp!]
\centering
\begin{subfigure}[b]{0.32\textwidth}
\includegraphics[width=\textwidth]{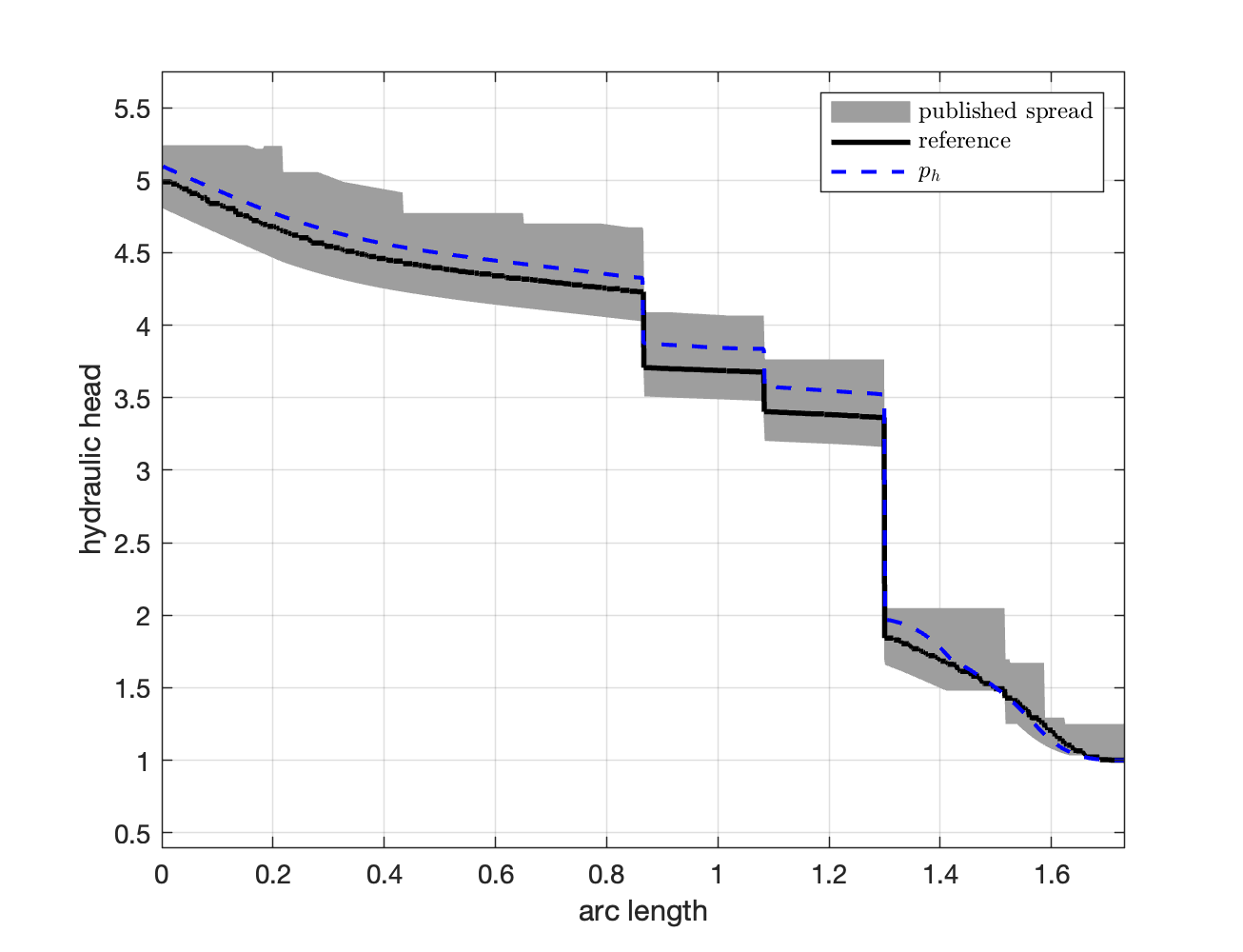}
\caption{$N=11$ ($1{,}331$ cells)}
\end{subfigure}
\hfill
\begin{subfigure}[b]{0.32\textwidth}
\includegraphics[width=\textwidth]{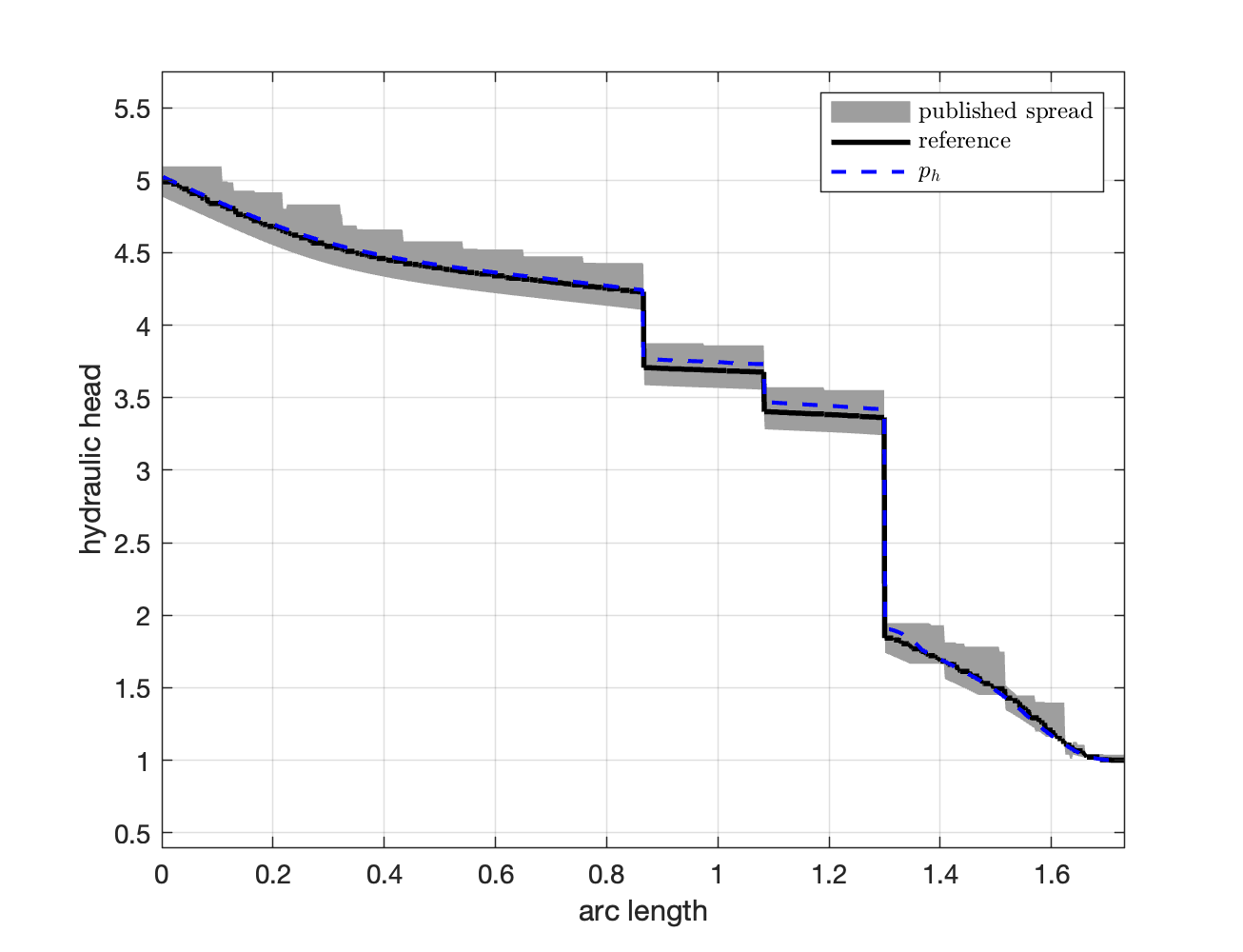}
\caption{$N=19$ ($6{,}859$ cells)}
\end{subfigure}
\hfill
\begin{subfigure}[b]{0.32\textwidth}
\includegraphics[width=\textwidth]{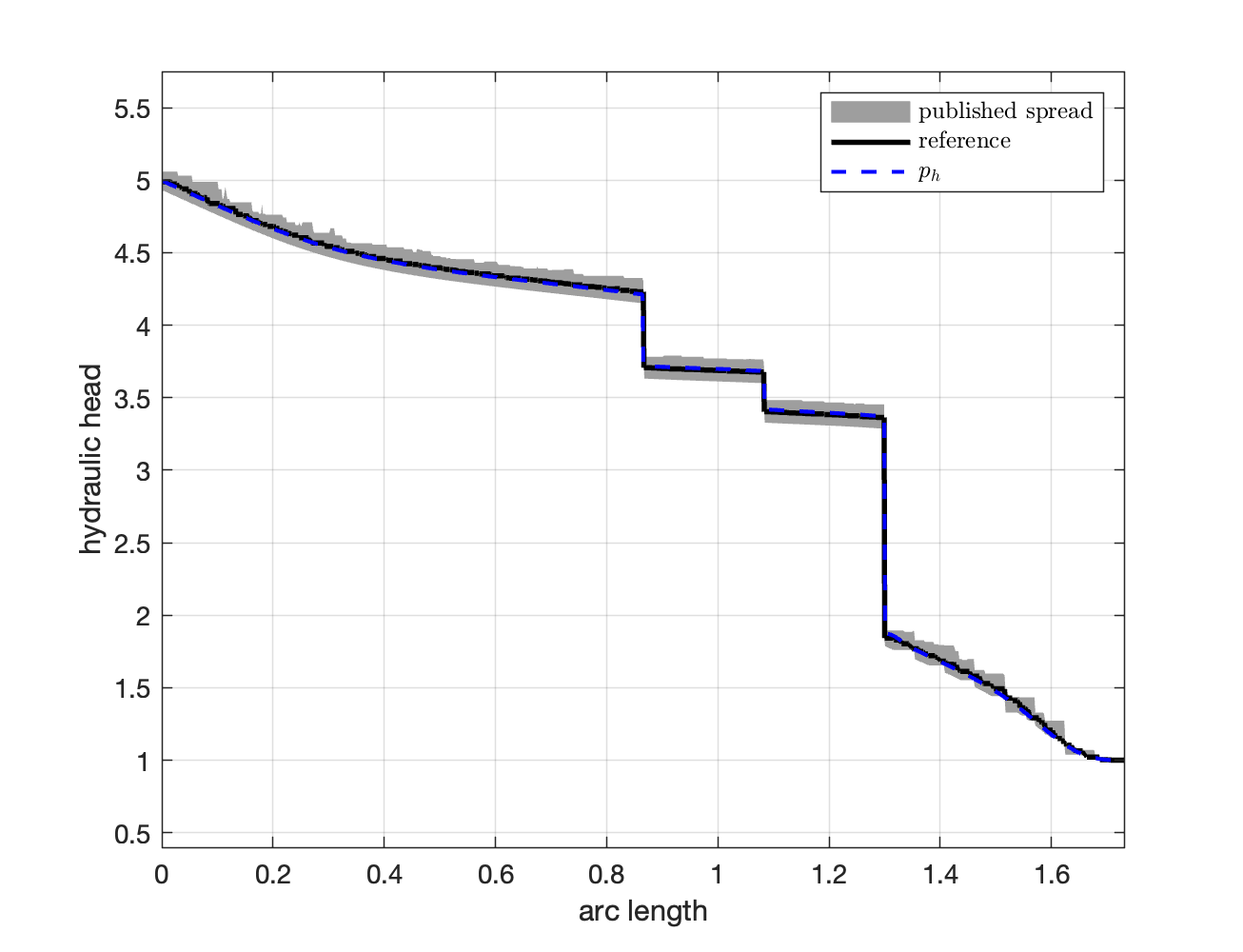}
\caption{$N=35$ ($42{,}875$ cells)}
\end{subfigure}
\caption{Example 7, blocking case, hydraulic head along the diagonal $(0,0,0)$--$(1,1,1)$ on the three levels of refinement.
The shaded region spans the $10$th--$90$th percentiles of the eleven methods that can represent a discontinuous head, at refinement levels of $\sim 500$, $\sim4$k, and $\sim32k$ cells.
The black curve is the reference solution
(USTUTT-MPFA).}
\label{fig:example7_block}
\end{figure}

Figures~\ref{fig:example7_cond} and~\ref{fig:example7_block} compare the head along the diagonal from $(0,0,0)$ to $(1,1,1)$ with the published data, in the same format as before: the band spans the $10$th--$90$th percentiles of the participating methods, and the reference is the USTUTT-MPFA solution on the benchmark's finest refinement.
In the conductive case, all fifteen participants enter the band, and the computed profiles reproduce the shape of the reference and lie within the band on all three levels of grid refinement.
In the blocking case, four of the participants (UNICE\_UNIGE-VAG\_Cont, UNICE\_UNIGE-HFV\_Cont, ETHZ\_USI-FEM\_LM and DTU-FEM\_COMSOL) keep the hydraulic head continuous across the fracture planes and cannot approach the reference under refinement.
These four schemes are excluded from the band of Figure~\ref{fig:example7_block}, which is built from the eleven remaining methods.
The reference profile is a staircase with sharp drops, and the our solution profiles capture every step at every resolution.
On the finest grid the computed profile is visually indistinguishable from the reference.

\subsection{Example 8: a network with small features in three dimensions}
\label{sec:example8}

\begin{figure}[htbp!]
\centering
\includegraphics[width=0.5\textwidth]{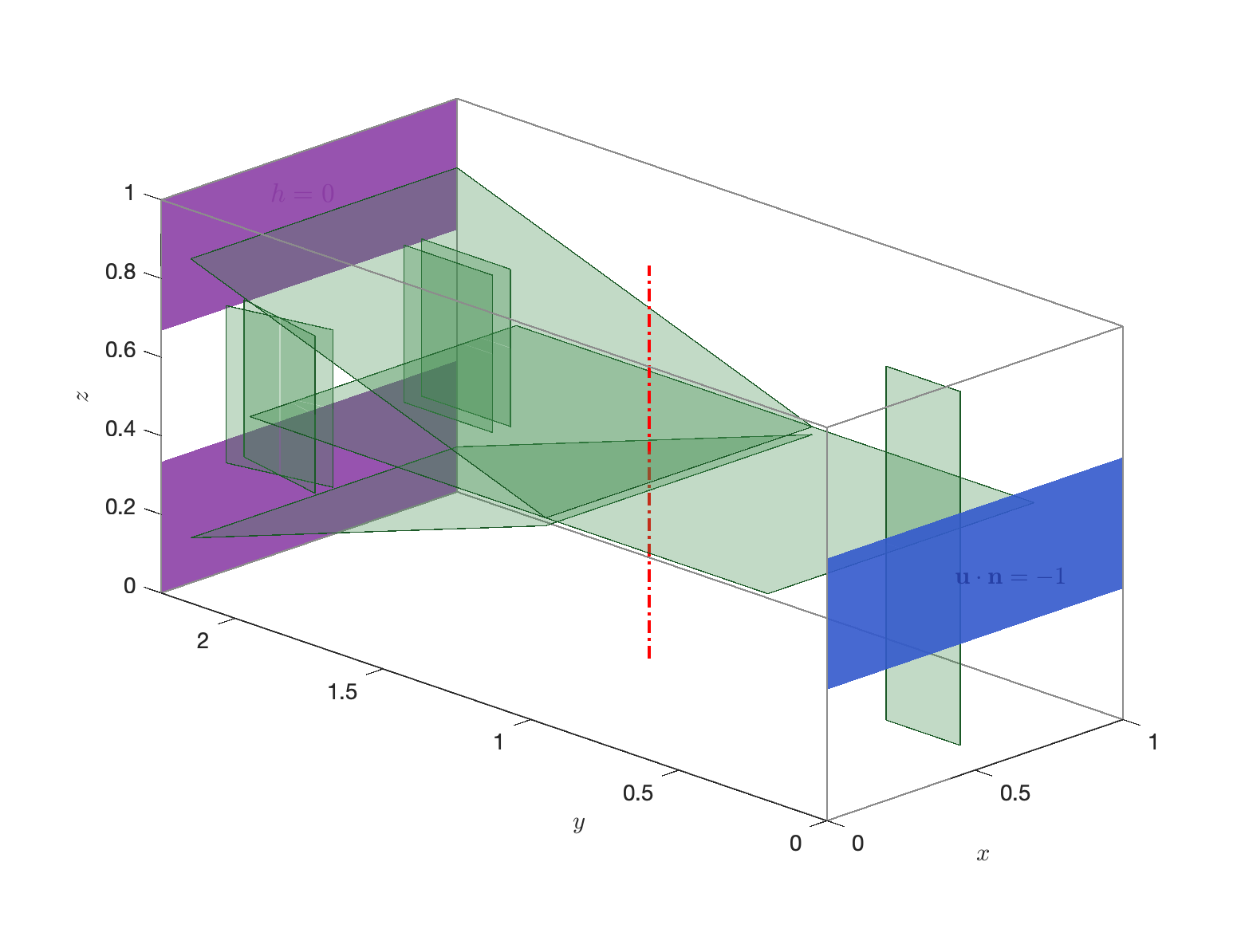}
\caption{Example 8, the computational domain with the eight fracture polygons, the inflow strip (blue, $\mathbf u\cdot\mathbf n=-1$), the two outlet strips (purple, $h=0$), and the sampling line (dash-dotted).}
\label{fig:example8_geometry}
\end{figure}

This experiment is the small-features case of
\cite{berre2021verification}, designed to probe how methods deal with geometric detail that troubles conforming mesh generation.  
The box $(0,1)\times(0,2.25)\times(0,1)$ contains eight conductive fractures combining axis-parallel and oblique planes with nearly touching members, small intersection segments, and features far below the domain scale.
The configuration is sketched in Figure~\ref{fig:example8_geometry}. 
All fractures have the tangential strength $a_fk_f=10^{2}$ and negligible normal resistance, and the matrix conductivity is $1$.
A uniform influx $\mathbf u\cdot\mathbf n=-1$ enters through the strip
$(0,1)\times\{0\}\times(1/3,2/3)$, the head $h=0$ is prescribed on the two strips $(0,1)\times\{2.25\}\times(0,\tfrac13)$ and $(0,1)\times\{2.25\}\times(\tfrac23,1)$, and
the remaining boundary is impermeable.  
We solve on grids of $25\times56\times25$ and $43\times97\times43$ cells.
The free-edge truncation rule of Appendix~\ref{app:three_dimensional_extension} acts on every polygon here, since none of the polygon edges lies on the domain boundary or on another
fracture.

\begin{figure}[htbp!]
\centering
\begin{subfigure}[b]{0.4\textwidth}
\includegraphics[width=\textwidth]{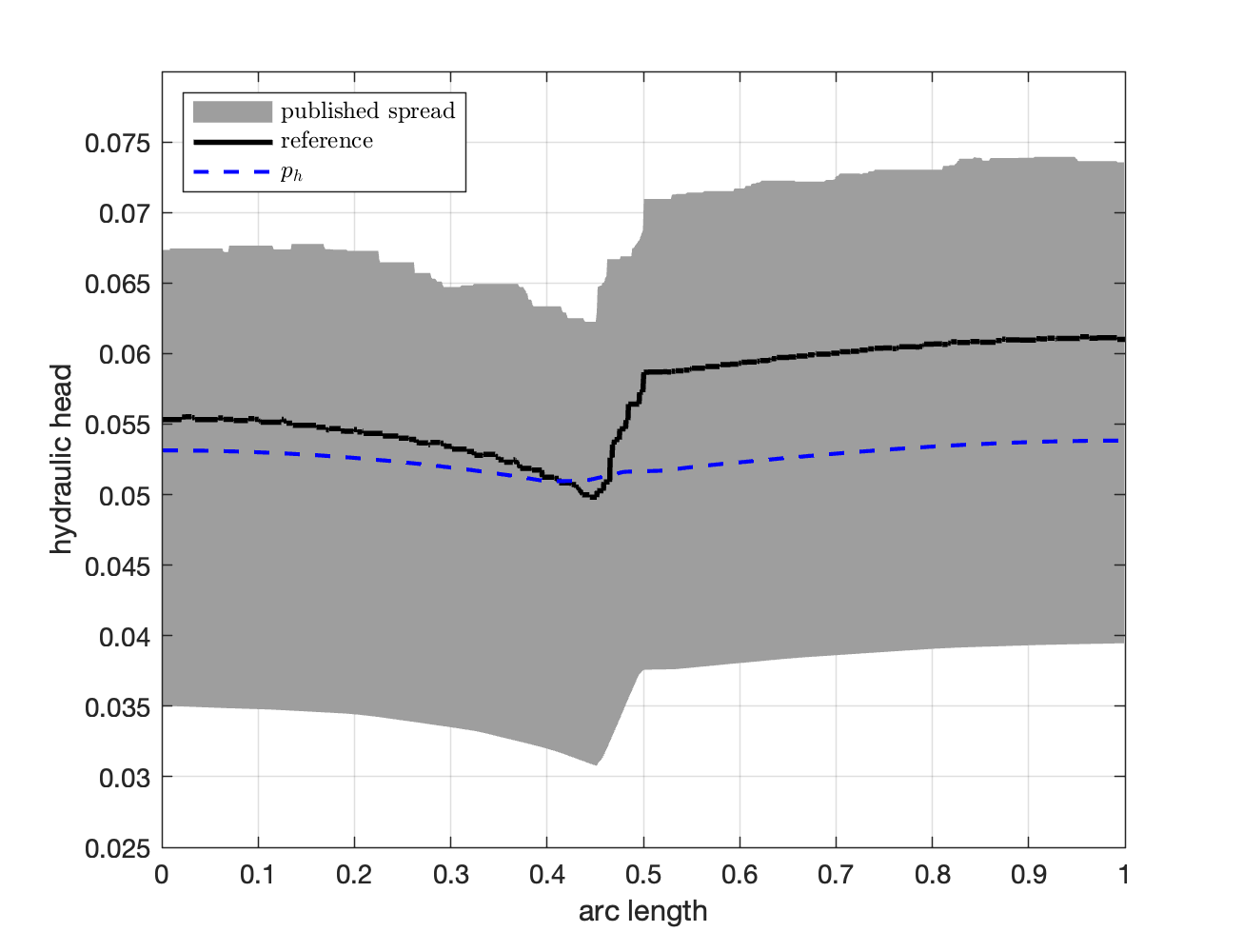}
\caption{$25\times56\times25$ ($35{,}000$ cells)}
\end{subfigure}
\hfill
\begin{subfigure}[b]{0.4\textwidth}
\includegraphics[width=\textwidth]{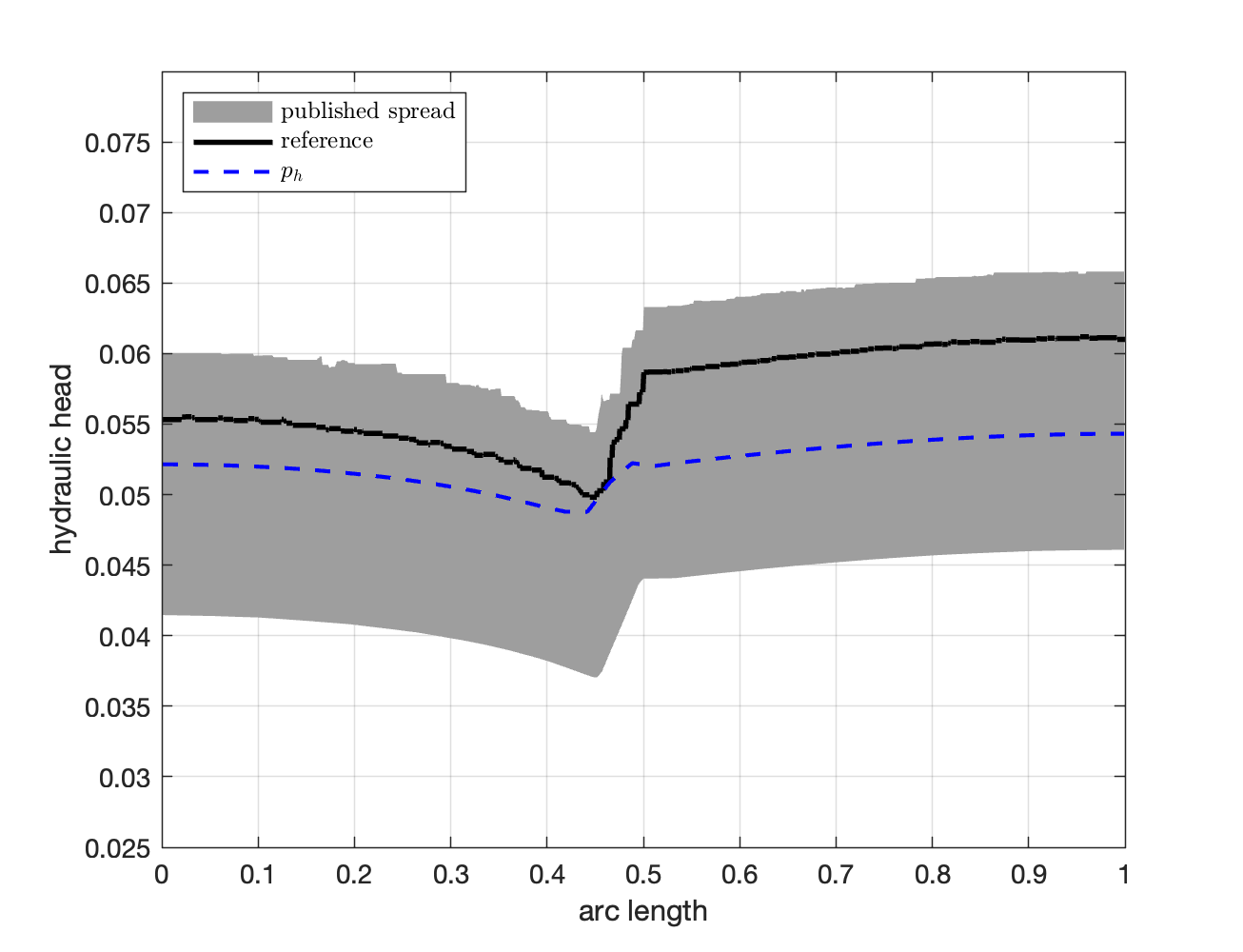}
\caption{$43\times97\times43$ ($179{,}353$ cells)}
\end{subfigure}
\caption{Example 8, hydraulic head along the line $(0.5,1.1,0)$--$(0.5,1.1,1)$ on the two levels of refinement.
The shaded region spans the 10th--90th
percentiles of the sixteen methods collected in \cite{berre2021verification} at refinement levels of $\sim 30$k and $\sim150$k cells.
The black curve is the reference solution (USTUTT-MPFA).}
\label{fig:example8_pol}
\end{figure}

Figure~\ref{fig:example8_pol} compares the head along the line $(0.5,1.1,0)$--$(0.5,1.1,1)$ with the published data.
The band spans the $10$th-$90$th percentiles of the sixteen participating methods, and the
reference is the USTUTT-MPFA solution a fine mesh.  
Our profiles lie in the interior of the wide
published spread.

\subsection{Example 9: a field case in three dimensions}
\label{sec:example9}

\begin{figure}[htbp!]
\centering
\includegraphics[width=0.5\textwidth]{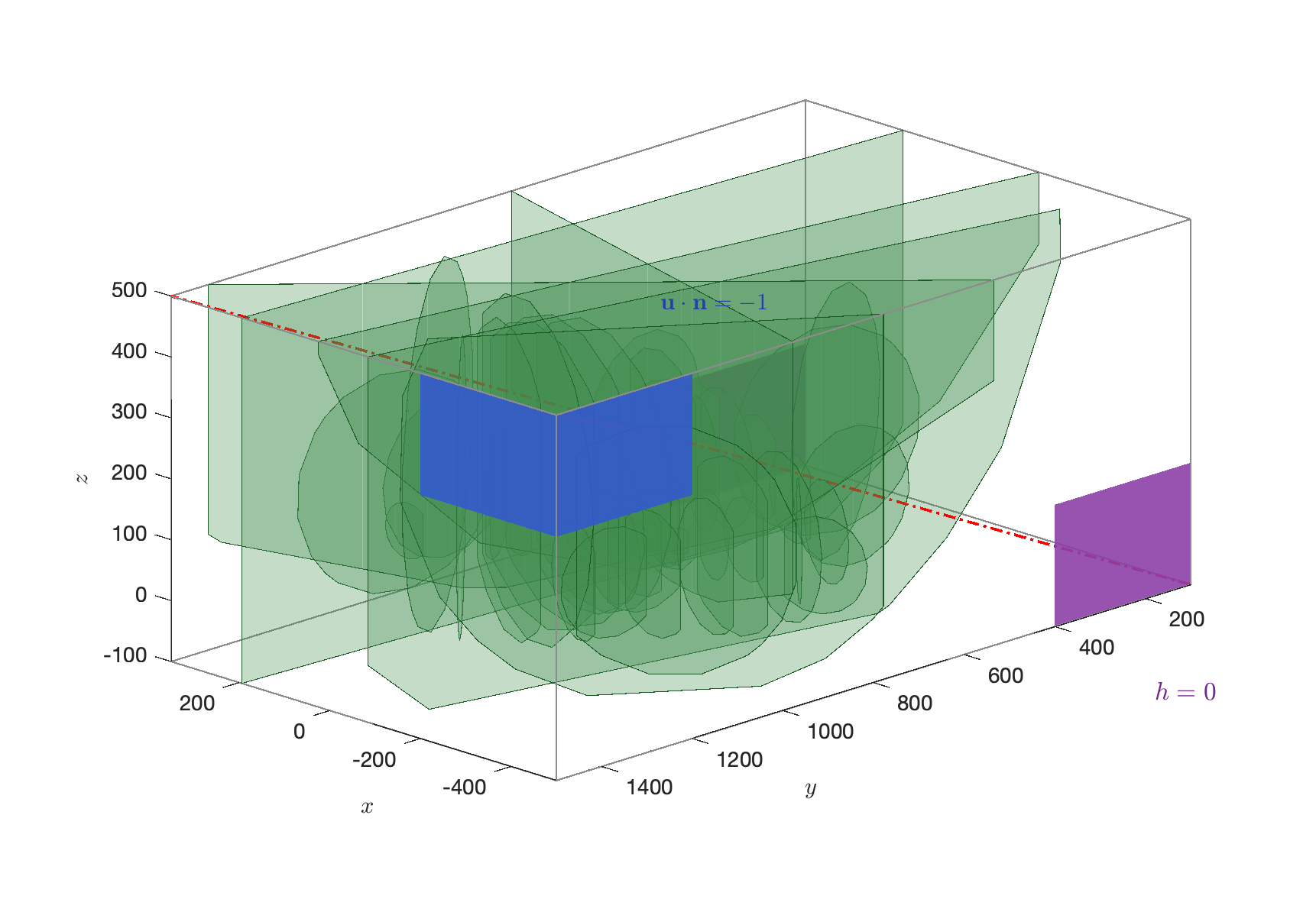}
\caption{Example 9, the computational domain with the 52 fracture polygons (green), the inlet patches (blue, $\mathbf u\cdot\mathbf n=-1$), the outlet patches (purple, $h=0$), and the sampling diagonals (dash-dotted).}
\label{fig:example9_geometry}
\end{figure}

The final experiment is the field case of \cite{berre2021verification}.  
$52$ fractures distribute in the box of $(-500,350)\times(100,1500)\times(-100,500)$.
The resulting network contains $106$ fracture intersections, and several fractures reach the domain boundary.  
Figure~\ref{fig:example9_geometry}
sketches the configuration.
The matrix conductivity is $1$, and all fractures are conductive with the tangential strength $a_fk_f=10^2$.
A uniform influx $\mathbf u\cdot\mathbf n=-1$ is prescribed on two rectangular inlet patches $(-500, -200)\times\{1500\}\times(300,500)$ and $\{-500\}\times(1200, 1500)\times(300, 500)$.
The head $h=0$ is prescribed on two outlet patches $\{-500\}\times(100, 400)\times(-100, 100)$ and $\{350\}\times(100, 400)\times(-100,100)$. 
The remaining boundary is impermeable.  
We use a grid of $61\times105\times43$ cells. 
The assembly and solution of the resulting $289$k unknowns take about a minute on a laptop.

\begin{figure}[htbp!]
\centering
\begin{subfigure}[b]{0.4\textwidth}
\includegraphics[width=\textwidth]{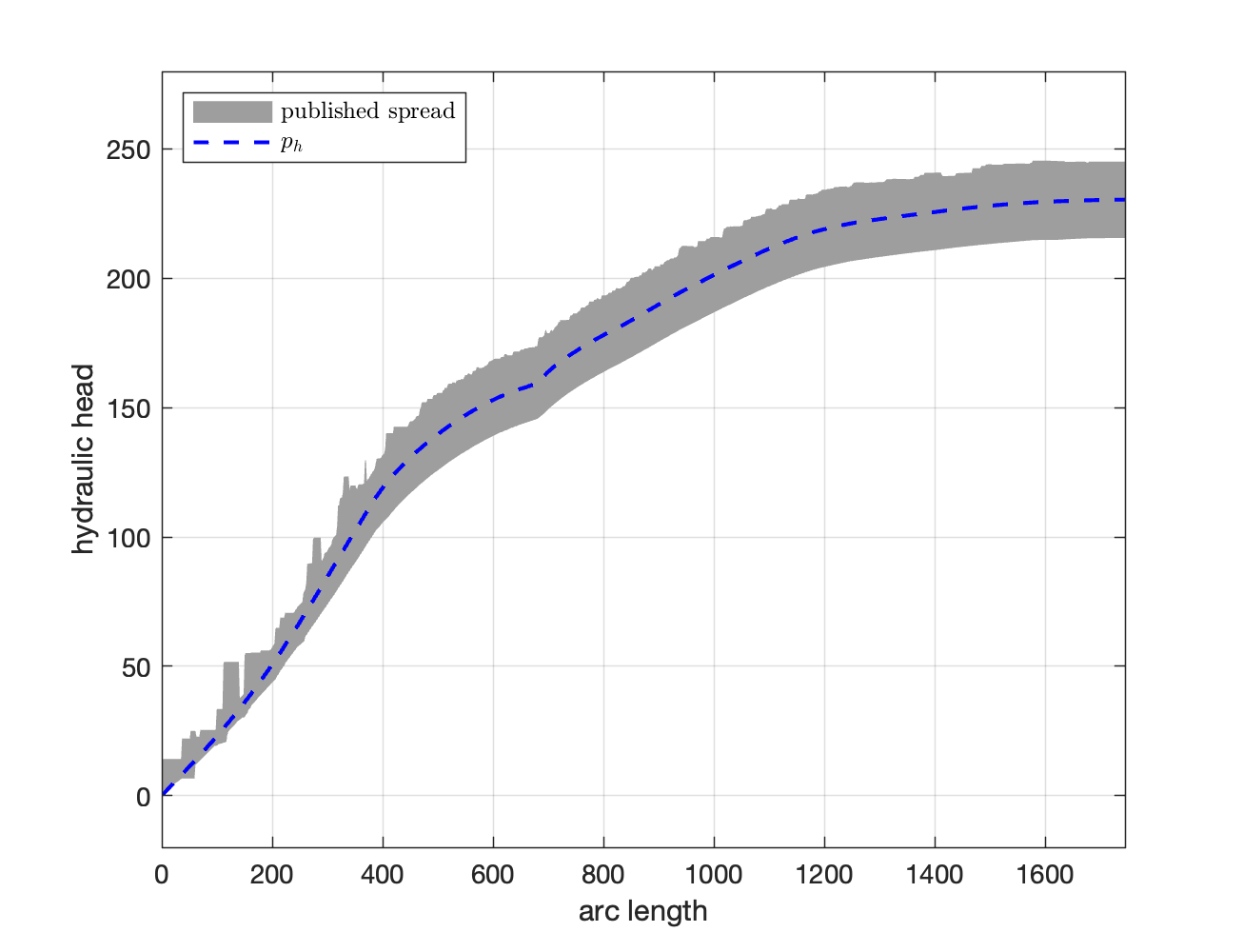}
\caption{slice 1, from $(-500,100,-100)$ to $(350,1500,500)$}
\end{subfigure}
\hfill
\begin{subfigure}[b]{0.4\textwidth}
\includegraphics[width=\textwidth]{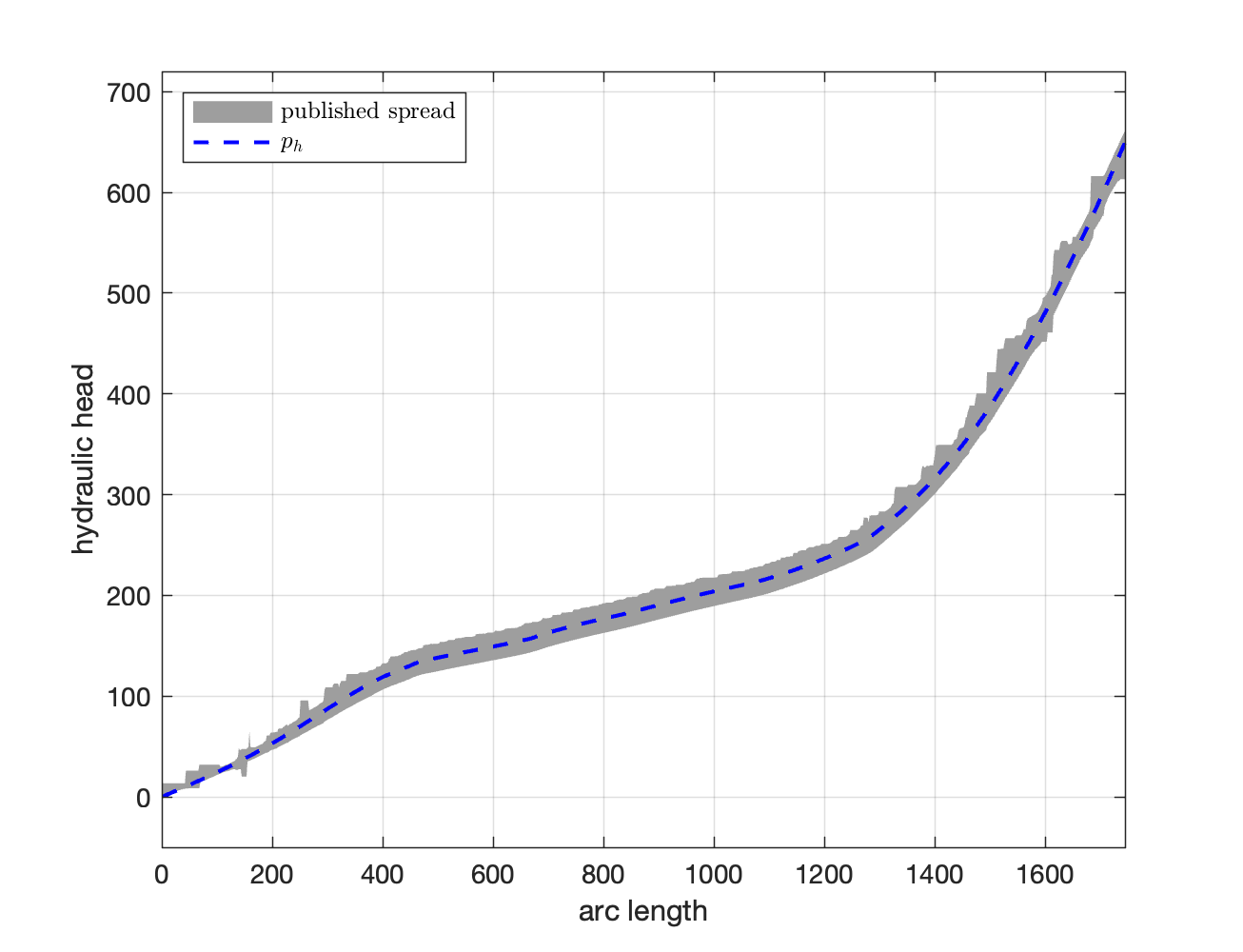}
\caption{slice 2, from $(350,100,-100)$ to $(-500,1500,500)$}
\end{subfigure}
\caption{Example 9, hydraulic head along the two sampling diagonals on the $61\times105\times43$ grid ($275{,}415$ cells).  
The shaded region spans the $10$th-$90$th percentiles of the fourteen methods collected in \cite{berre2021verification} at a refinement level of $\sim 260$k cells.}
\label{fig:example9_lines}
\end{figure}

Because a refinement study was considered
infeasible for this geometry, the benchmark took a single resolution of approximately $260$k cells.
Figure~\ref{fig:example9_lines} compares the head along the two box
diagonals with the fourteen published solutions with the shaded band spanning their $10$th to $90$th percentiles.
The first sampling line runs from $(-500,100,-100)$ to $(350,1500,500)$, and the second sampling line runs from $(350,100,-100)$ to $(-500,1500,500)$.
Along both lines the computed profile reproduces the shape of the published solutions and remains inside the
spread over nearly its entire length.

\section{Concluding remarks}
\label{sec:conclusion}

We have developed an energy-based finite difference discrete fracture model for single-phase flow in porous media containing conductive fractures and low-permeability barriers. 
The method is obtained by discretizing the matrix, barrier, and fracture contributions to the continuous energy directly on an unfitted Cartesian grid. 
Barrier effects are represented by local pressure-jump variables and eliminated before global assembly, whereas conductive fractures contribute local positive-semidefinite updates associated with tangential flow. 
The resulting system retains only the original nodal pressure unknowns, remains symmetric and compact, and coincides with the standard finite difference discretization away from the interfaces.
The eliminated barrier jumps and the tangential fracture fluxes can also be recovered through inexpensive local post-processing. 
The numerical results demonstrate that the method captures both enhanced fracture transport and barrier-induced pressure discontinuities for a range of isolated and network configurations.

Two directions appear to be natural extensions of the methodology developed in this paper. First, we plan to apply the energy principle to TPFA and MPFA discrete fracture models on non-conforming meshes. 
Second, an analogous complementary-energy formulation for mixed finite element discretizations should also be investigated.

\appendix
\section{Extension to three dimensional space}
\label{app:three_dimensional_extension}

The extension to three dimensions follows the same energy principle. 
We briefly outline the construction and state the resulting schemes.

Let $\Omega\subset\mathbb R^3$ be partitioned into Cartesian cells of side lengths $h_x$, $h_y$, and $h_z$. Barriers and conductive fractures are now represented by surfaces $\Gamma_b$ and $\Gamma_f$, respectively. 
The continuous energy has the same form as in two dimensions, with the line integrals replaced by surface integrals and the tangential derivative along a fracture replaced by the surface
gradient $\nabla_\tau$:
\[
\mathcal J(p)
=
\frac12\int_{\Omega_m}
\mathbf K_m\nabla p\cdot\nabla p\,d\mathbf x
+
\frac12\int_{\Gamma_b}
\frac{k_b}{a_b}\llbracket p\rrbracket^2\,dS
+
\frac12\int_{\Gamma_f}
a_fk_f|\nabla_\tau p|^2\,dS
-
\ell(p).
\]

For a Cartesian cell $C$, let $\mathbf p_C\in\mathbb R^8$ contain the pressures at its eight vertices. 
For each coordinate direction $i\in\{x,y,z\}$, denote by $\mathcal E_i(C)$ the four edges of $C$ parallel to that direction and define
\[
\left\langle D_i^2\right\rangle_C
=
\frac14\sum_{e\in\mathcal E_i(C)}(D_i^e p_C)^2,
\qquad
\overline D_i p_C
=
\frac14\sum_{e\in\mathcal E_i(C)}D_i^e p_C.
\]
For a representative full permeability tensor
$\mathbf K_C=(k_{ij,C})_{i,j=1}^3$, the three-dimensional matrix energy
may be discretized by
\begin{equation}
\label{eq:3d_tensor_cell_energy}
\mathcal E_{h,C}^{m}(\mathbf p_C)
=
\frac{|C|}{2}
\left[
\sum_{i=1}^3 k_{ii,C}\left\langle D_i^2\right\rangle_C
+
2\sum_{1\leq i<j\leq3}
k_{ij,C}\,
\overline D_i p_C\,\overline D_j p_C
\right]
=
\frac12\mathbf p_C^{\mathsf T}S_C\mathbf p_C,
\end{equation}
where $S_C\in\mathbb R^{8\times8}$ is symmetric positive semidefinite and satisfies $S_C\mathbf1=\mathbf0$. 
For scalar permeability, $\mathbf K_m=k_m\mathbf I$, the assembled energy reduces to
\[
\mathcal E_h^m(\mathbf p)
=
\frac12\sum_{e=(P,Q)\in\mathcal E_h}
c_e(p_P-p_Q)^2,
\]
with
\begin{equation}
\label{eq:3d_edge_conductance}
c_e
=
\begin{cases}
\displaystyle k_e\frac{h_yh_z}{h_x},
& e\parallel\mathbf e_x,\\[2mm]
\displaystyle k_e\frac{h_xh_z}{h_y},
& e\parallel\mathbf e_y,\\[2mm]
\displaystyle k_e\frac{h_xh_y}{h_z},
& e\parallel\mathbf e_z.
\end{cases}
\end{equation}
Its first variation is the standard seven-point finite difference operator.

For the scalar barrier treatment, let a barrier surface intersect an edge $e$ transversely at $\mathbf x_e$. The projected surface weight is
defined by
\begin{equation}
\label{eq:3d_projected_barrier_weight}
\omega_e
=
\begin{cases}
|\mathbf{n}_{\Gamma,x}(\mathbf x_e)|h_yh_z,
& e\parallel\mathbf e_x,\\
|\mathbf{n}_{\Gamma,y}(\mathbf x_e)|h_xh_z,
& e\parallel\mathbf e_y,\\
|\mathbf{n}_{\Gamma,z}(\mathbf x_e)|h_xh_y,
& e\parallel\mathbf e_z,
\end{cases}
\end{equation}
and
\[
d_e=\omega_e\frac{k_{b,e}}{a_{b,e}}.
\]
The same edgewise elimination as in two dimensions gives
\begin{equation}
\label{eq:3d_effective_edge_conductance}
\widehat c_e
=
\left(
\frac1{c_e}
+
\sum_{\alpha\in\mathcal A_e}\frac1{d_{e,\alpha}}
\right)^{-1}
\end{equation}
on a barrier-cut edge, while $\widehat c_e=c_e$ on an uncut edge.
Thus, the scalar barrier-only scheme retains the original seven-point
stencil.

For full-tensor permeability, suppose that the barrier surfaces
partition a cut cell into regions $C_0,\ldots,C_m$, with $C_0$ chosen
as the reference region. Let
$E_C\in\{0,1\}^{8\times m}$ be the vertex--region incidence matrix.
If a barrier patch $\Sigma_{C,\alpha}$ separates regions
$C_{i_\alpha}$ and $C_{j_\alpha}$, define
\[
d_{C,\alpha}
=
\int_{\Sigma_{C,\alpha}}\frac{k_b}{a_b}\,dS
\]
and let
\[
P_C
=
\sum_{\alpha}
d_{C,\alpha}
(\mathbf e_{i_\alpha}-\mathbf e_{j_\alpha})
(\mathbf e_{i_\alpha}-\mathbf e_{j_\alpha})^{\mathsf T},
\qquad
\mathbf e_0:=\mathbf0.
\]
Eliminating the region offsets gives the barrier-modified cell matrix
\begin{equation}
\label{eq:3d_condensed_barrier_matrix}
\widehat S_C
=
S_C
-
S_CE_C
\left(
E_C^{\mathsf T}S_CE_C+P_C
\right)^{-1}
E_C^{\mathsf T}S_C.
\end{equation}
On an uncut cell, $\widehat S_C=S_C$.

It remains to discretize tangential transport on the conductive
fracture surfaces. Let
$\boldsymbol\lambda_C=(\lambda_{C,1},\ldots,\lambda_{C,8})^{\mathsf T}$
be the trilinear basis on $C$ and define
\[
G_C(\mathbf x)
=
\begin{pmatrix}
\nabla\lambda_{C,1}(\mathbf x)&\cdots&
\nabla\lambda_{C,8}(\mathbf x)
\end{pmatrix}.
\]
For a fracture patch
$\Sigma_{C,\alpha}\subset\Gamma_f\cap C$ with unit normal
$\mathbf n_{\Gamma}$, set
\[
\Pi_\tau
=
I-\mathbf n_{\Gamma}\mathbf n_{\Gamma}^{\mathsf T}.
\]
Its discrete tangential energy is
\[
\mathcal E_{h,C,\alpha}^{f}(\mathbf p_C)
=
\frac12\mathbf p_C^{\mathsf T}
F_{C,\alpha}\mathbf p_C,
\]
where
\begin{equation}
\label{eq:3d_fracture_matrix}
F_{C,\alpha}
=
\int_{\Sigma_{C,\alpha}}
a_fk_f\,
G_C^{\mathsf T}\Pi_\tau G_C\,dS.
\end{equation}
The surface integral may be evaluated by triangulating the fracture patch and applying a standard quadrature rule. 
Summing over all fracture patches in the cell gives
\[
F_C=\sum_{\alpha\in\mathcal A_C^f}F_{C,\alpha}.
\]
The matrix $F_C$ is symmetric positive semidefinite and satisfies $F_C\mathbf1=\mathbf0$.
As in two dimensions, the free boundary of a fracture is treated by truncation. 

The final three-dimensional nodal schemes therefore take the same form as their two-dimensional counterparts. 
For scalar matrix permeability,
\begin{equation}
\label{eq:3d_scalar_final_scheme}
\sum_{Q:\,(P,Q)\in\mathcal E_h}
\widehat c_{PQ}(p_P-p_Q)
+
\sum_{\substack{C\in\mathcal T_h\\P\in C}}
(F_C\mathbf p_C)_P
=
(\mathbf b_h)_P,
\qquad
P\in\mathcal N_h\setminus\mathcal N_{h,D}.
\end{equation}
For full-tensor permeability,
\begin{equation}
\label{eq:3d_tensor_final_scheme}
\sum_{\substack{C\in\mathcal T_h\\P\in C}}
\left[
(\widehat S_C+F_C)\mathbf p_C
\right]_P
=
(\mathbf b_h)_P,
\qquad
P\in\mathcal N_h\setminus\mathcal N_{h,D}.
\end{equation}
Both schemes retain only the original Cartesian-grid nodal pressures.
The scalar barrier-only discretization has a seven-point stencil, whereas the full-tensor and fracture contributions may produce couplings within the surrounding $27$-point stencil.

The eliminated barrier jumps are recovered by the same local back-substitution formulas as in two dimensions. 
In particular,
\[
\mathbf t_C^\star
=
\left(
E_C^{\mathsf T}S_CE_C+P_C
\right)^{-1}
E_C^{\mathsf T}S_C\mathbf p_{h,C},
\]
while the tangential fracture flux is recovered from
$-a_fk_f\nabla_\tau p_h$ on each fracture patch. No additional global
solve is required.


\bibliographystyle{plain}

\bibliography{refs}

\end{document}